\documentclass{amsart}
\usepackage{amssymb}

\newtheorem{thm}{Theorem}[section]
\newtheorem{prop}[thm]{Proposition}
\newtheorem{cor}[thm]{Corollary}
\newtheorem{lem}[thm]{Lemma}
\newtheorem{rema}[thm]{Remark}
\newtheorem{defn}[thm]{Definition}

\numberwithin{equation}{section}

\newcommand{\nor}[0]{\Upsilon}
\newcommand{\sou}[0]{\Delta}

\newcommand{\Z}{\mathbb{Z}_+}

\begin{document}

\title[On uniformization of DeWitt super-Riemann surfaces]{On uniformization of N=2 superconformal and N=1 superanalytic DeWitt super-Riemann surfaces}

%Information for first author 
\author{Katrina Barron}
\address{Department of Mathematics, University of Notre Dame,
Notre Dame, IN 46556 and Max Planck Institute for Mathematics, 53111 Bonn, Germany}
\email{kbarron@nd.edu}


\date{February 11, 2011}


\thanks{Supported in part by NSA grant MSPF-07G-169, and by a research grant from the Max Planck Institute for Mathematics, Bonn, Germany.}

\begin{abstract}
We prove a general uniformization theorem for N=2 superconformal and N=1 superanalytic DeWitt super-Riemann surfaces, showing that in general an N=2 superconformal (resp. N=1 superanalytic) DeWitt super-Riemann surface is  N=2 superconformally (resp., N=1 superanalytically) equivalent to a manifold with transition functions containing no odd functions of the even variable if and only if  a certain cohomology group is trivial, namely the first \v  Cech cohomology group of the body Riemann surface with coefficients in the sheaf consisting of the reciprocal of a line bundle tensor the holomorphic vector fields over the body.  In particular, this gives a general criteria for when a DeWitt N=1 superanalytic super-Riemann surface is N=1 superanalytically equivalent to a ringed-space $(1,1)$-supermanifold, as studied in the algebro-geometric setting.   As a consequence of this general classification result, there is a countably infinite family of N=2 superconformal equivalence classes of N=2 superconformal DeWitt super-Riemann surfaces with genus-zero compact body, and N=2 superconformal DeWitt super-Riemann surfaces with simply connected body are classified up to N=2 superconformal equivalence by conformal equivalence classes of holomorphic line bundles over the underlying body Riemann surface.   In addition, we show that N=2 superconformal DeWitt super-Riemann surfaces with compact genus-one body and transition functions which correspond to the trivial cocycle in the first \v  Cech cohomology group of the body Riemann surface with coefficients in the reciprocal of a line bundle tensor the sheaf of holomorphic vector fields over the body are classified up to N=2 superconformal equivalence by theta functions associated to the underlying torus up to type modulo the trivial theta functions, or in other words, by holomorphic line bundles over the torus modulo conformal equivalence.   We also give the corresponding results for the uniformization of N=1 superanalytic DeWitt super-Riemann surfaces of genus zero or one.  
\end{abstract}

\maketitle

\section{Introduction}

In this paper, we prove a general uniformization theorem for N=2 superconformal and N=1 superanalytic DeWitt super-Riemann surfaces.   In N=2 superconformal field theory, the surfaces swept out by propagating strings with N=2 superconformal symmetry are N=2 superconformal DeWitt super-Riemann surfaces, \cite{Fd}, \cite{DPZ}, \cite{FMS}, \cite{Wa}, \cite{Ge}.   In order to construct an N=2 superconformal field theory, one needs, in particular, a precise description of the moduli space of N=2 superconformal DeWitt super-Riemann surfaces, under N=2 superconformal equivalence.   

There are two main approaches to supermanifolds, the ``concrete" or ``DeWitt" approach \cite{Fd}, \cite{D}, \cite{Ro} and the ``ringed-space" approach \cite{Leites}, \cite{Manin1}, \cite{Rothstein}.    The DeWitt approach and the ringed-space approach to supermanifolds are equivalent if one restricts the supermanifolds in the DeWitt approach to only allow for transition functions which do not include components that are odd functions of an even variable \cite{Batchelor}, \cite{Ro}; see Remark \ref{DeWitt-v-ringed-space-remark}.  However, to study the worldsheet geometry underlying superconformal field theory, one needs to include these more general transition functions, which are naturally incorporated using the DeWitt approach.   Only if one allows the worldsheets swept out by particles modeled as superstrings to include these more general supermanifolds does one have boson-fermion supersymmetry in the algebra of correlation functions governed by these worldsheets (cf. \cite{B-memoir}, \cite{B-iso}, \cite{B-n2moduli}, \cite{B-axiomatic}); that is, it is only then that one has the possibility of a nontrivial representation of the Neveu-Schwarz algebra of infinitesimal superconformal symmetries imposed on the space of particle states rather than merely a representation of the Virasoro algebra of infinitesimal conformal symmetries.  Using the ringed-space approach, in order to incorporate the more general transition functions allowed in the DeWitt approach and in superconformal field theories, one must consider, for instance, families of ringed-space supermanifolds over a given supermanifold; see for instance \cite{LR}, \cite{Falqui-Reina-1}, \cite{Manin2},  \cite{BR}.

The main result we prove in this paper states, in particular, that a  general N=2 superconformal DeWitt super-Riemann surface with body Riemann surface $M_B$ is N=2 superconformally equivalent to an N=2 superconformal DeWitt super-Riemann surface with transition functions which do not include components that are odd functions of an even variable if and only if a certain cohomology vanishes, namely the first \v Cech cohomology of $M_B$ with coefficients in the space of holomorphic vector fields over $M_B$ tensor the reciprocal of a line bundle over $M_B$, denoted $\check{H}^1(M_B, \mathcal{L}^{-1} \otimes TM_B)$.   In other words, an N=2 superconformal DeWitt super-Riemann surface is, in general, N=2 superconformally equivalent to a supermanifold as defined in the ringed-space approach if and only if   $\check{H}^1(M_B, \mathcal{L}^{-1} \otimes TM_B) = 0$.

%Thus, for instance, if $M_B$ is a simply connected Riemann surface, the moduli space of N=2 
%superconformal DeWitt super-Riemann surfaces with body $M_B$ under N=2 superconformal 
%equivalence is given by the moduli space of N=2 superconformal ringed-space supermanifolds.  
%However, if, for instance $M_B$ is genus one, the moduli space is considerably more complicated for 
%N=2 superconformal DeWitt super-Riemann surfaces than for N=2 superconformal ringed-space 
%supermanifolds. 

As shown in \cite{DRS}, N=2 superconformal super-Riemann surfaces are equivalent to N=1 superanalytic super-Riemann surfaces, and thus our results apply to these surfaces as well, under this equivalence.   In fact it is interesting to note, that our main uniformization theorem is proved in the N=2 superconformal ``homogeneous coordinate setting" where the dependency on the \v Cech cohomology associated to the underlying body Riemann surface is transparent when one considers consistency conditions of local coordinate transformations on triple overlaps.  This dependency is not transparent if one looks at these consistency conditions for N=1 superanalytic super-Riemann surfaces, as discussed in Remark \ref{setting-remark}, or for N=2 superconformal super-Riemann surfaces in the ``nonhomogeneous coordinate setting" as discussed in Section \ref{nonhomo-section}.

The problem of the classification of N=1 superanalytic ringed-space manifolds has been studied in for instance \cite{Manin1}, \cite{Rothstein}, \cite{Vaintrob}, \cite{Manin2}.  Our main uniformization theorem, Theorem 4.1,  (along with Corollary 4.2) provides criteria for when an N=2 superconformal (resp. N=1 superanalytic) DeWitt super-Riemann surface is N=2 superconformally (resp. N=1 superanalytically) equivalent to a ringed-space supermanifold for which these prior results on classifying N=1 superanalytic ringed-space manifolds apply.   

Thus, for instance, when the body of an underlying N=2 superconformal (resp. N=1 superanalytic) super-Riemann surface is simply connected, the vanishing of the cohomologies that appear in Theorem 4.1, imply that any of this simply connected DeWitt supermanifolds is equivalent to a simply connected ringed-space supermanifold, and previous classification results regarding these ringed-space manifolds, as for instance in \cite{Manin1}, \cite{Manin2}, can be applied.  
In the present paper, using the DeWitt setting, we give concrete realizations of the further classification of these simply connected supermanifolds as already obtained in the ringed-space approach once one applies our Corollary 4.2.  That is, using our Theorem 4.1 and Corollary 4.2, we classify, up to N=2 superconformal (resp. N=1 superanalytic) equivalence, N=2 superconformal (resp. N=1 superanalytic) DeWitt super-Riemann surfaces with simply connected body.  We show that there are unique, up to N=2 superconformal (resp. N=1 superanalytic) equivalence, N=2 superconformal (resp. N=1 superanalytic) DeWitt structures over the complex plane and complex upper-half plane, and a countably infinite number of inequivalent  N=2 superconformal (resp. N=1 superanalytic) DeWitt structures over the Riemann sphere.

In addition, we give the classification of N=2 superconformal (resp. N=1 superanalytic) DeWitt super-Riemann surfaces with compact genus-one body and with transition functions restricted to those that do not contain components involving odd functions of an even variable, i.e., those found already in the ringed-space approach.  For N=2 superconformal  (resp. N=1 superanalytic) DeWitt structures over a complex torus,  we show that if we restrict to supermanifolds where the transition functions, in particular, do not contain components involving odd functions of an even variable, then there is a doubly infinite family of N=2 superconformal (resp. N=1 superanalytic) equivalence classes of genus-one N=2 superconformal (resp. N=1 superanalytic) DeWitt super-Riemann surfaces over a given torus, and that the N=2 superconformal (resp. N=1 superanalytic) equivalence classes are determined by theta functions over the underlying lattice defining the complex torus, up to type and up to equivalence with respect to the trivial theta functions.  

The classification of simply connected N=2 superconformal (resp. N=1 superanalytic) DeWitt super-Riemann surfaces and this subclass of genus-one N=2 superconformal (resp. N=1 superanalytic)  DeWitt super-Riemann surfaces can be restated as follows: N=2 superconformal (resp. N=1 superanalytic) DeWitt super-Riemann surfaces with simply connected body or with compact genus-one body and transition functions restricted, in particular,  to contain no odd functions of an even variable are classified up to N=2 superconformal (resp. N=1 superanalytic) equivalence by holomorphic line bundles over the underlying Riemann surface up to conformal equivalence.

This work has important implications for N=2 superconformal field theory \cite{BPZ}, \cite{Fd}, \cite{DPZ}, \cite{FMS}, \cite{Wa}, \cite{Ge}, and its connections to mathematics through, for instance, the phenomenon of mirror symmetry \cite{GP}, \cite{CdGP}, \cite{Gi}, \cite{CK}, \cite{HKKP},  and the theory of vertex operator superalgebras \cite{FLM}, \cite{KT}, \cite{DPZ}, \cite{YZ}, \cite{LVW}, \cite{FFR}, \cite{DL}, \cite{B-axiomatic}, in addition to having applications, through the theory of vertex operator superalgebras, to aspects of number theory \cite{Milas}, \cite{STT} and integrable systems \cite{Manin-Radul}, \cite{BR}.   

This uniformization result is a crucial step in constructing an N=2 superconformal field theory following the work of Huang  \cite{H-book}-\cite{H2005-2} and Huang and Lepowsky \cite{HL1}-\cite{HL7} in the nonsuper case, and the author \cite{B-announce}-\cite{B-change} in the N=1 superconformal case.   The extension of this program to the N=2 superconformal case was initiated in \cite{B-n2moduli}.  The present paper allows the results of \cite{B-n2moduli} to be put into the context of the description of the entire moduli space of genus-zero N=2 superconformal worldsheets modeled as genus-zero N=2 superconformal DeWitt super-Riemann surfaces with half-infinite tubes.  In \cite{B-autogroups}, we have determined the Lie supergroups of automorphisms of each N=2 superconformal equivalence class of genus-zero N=2 superconformal super-Riemann surface, and this work along with that of the present paper provides the necessary results to proceed in defining the moduli space of N=2 superspheres with tubes and a sewing operation modeling the worldsheet approach to genus-zero two-dimensional N=2 superconformal field theory.  

In particular, we note here that in both the nonsuper case and N=1 super case there is a unique genus-zero surface up to conformal or N=1 superconformal equivalence, respectively.  However, as proved in this paper, there is a countably infinite number of equivalence classes of N=2 superconformal genus-zero super-Riemann surfaces.  Thus, in the N=2 superconformal case, the operad structure on the moduli space of genus-zero worldsheets will have considerably more structure in comparison to the nonsuper and N=1 super cases.  Consequently the algebraic structure imposed on the space of particle states by the worldsheet operad structure has considerably more structure in the N=2 superconformal case than in the analogous conformal and N=1 superconformal cases, where the structure of a vertex operator algebra \cite{H-book} or N=1 Neveu-Schwarz vertex operator superalgebra \cite{B-vosas}, \cite{B-memoir}, \cite{B-iso}, respectively, is imposed.  In fact, an algebra over the operadic structure of genus-zero N=2 superconformal worldsheets will  necessarily carry much more structure than that of an N=2 Neveu-Schwarz vertex operator superalgebra \cite{B-axiomatic}.

This paper is organized as follows:  In Section \ref{preliminaries-section}, the basic notions of superfunctions, superanalyticity and N=1 and  N=2 superconformality are presented as well as the notion of DeWitt supermanifold and N=2 superconformal super-Riemann surface.   As subclasses of these notions, we also present notation for superfunctions and supermanifolds found in the more restricted ringed-space approach.  In Section \ref{equivalence-section}, we recall the construction of an equivalence between N=2 superconformal DeWitt super-Riemann surfaces and N=1 superanalytic DeWitt super-Riemann surfaces following \cite{DRS}.

In Section \ref{uniformization-section}, we prove the main theorem of this paper, Theorem 4.1.  We prove that in general, an N=2 superconformal DeWitt super-Riemann surface with body Riemann surface $M_B$ is N=2 superconformally equivalent to an N=2 superconformal DeWitt super-Riemann surface with transition functions which do not include components that are odd functions of an even variable if and only if the first \v Cech cohomology of $M_B$ with coefficients in the space of holomorphic vector fields over $M_B$ tensor the reciprocal of a line bundle over $M_B$ is trivial.  Furthermore, the remaining even transition functions can be further reduced to have no soul components in general, if and only if the  first \v Cech cohomology of $M_B$ with coefficients in the space of holomorphic vector fields over $M_B$ is trivial.  We then use the equivalence between N=2 superconformal and N=1 superanalytic DeWitt super-Riemann surfaces to state an analogous theorem for N=1 superanalytic DeWitt super-Riemann surfaces.

In Section \ref{sphere-section}, we study N=2 superconformal DeWitt super-Riemann surfaces over the Riemann sphere with certain restricted N=2 superconformal coordinate transformations.  We then use the uniformization theorem, Theorem 4.1, to classify genus-zero N=2 superconformal DeWitt super-Riemann surfaces.  In particular, we show that compact genus-zero N=2 superconformal DeWitt super-Riemann surfaces are classified up to N=2 superconformal equivalence by a countably infinite family of superspheres, and there is a bijection between these equivalence classes and conformal equivalence classes of holomorphic line bundles over the Riemann sphere.  We then state the analogous results for genus-zero N=1 superanalytic DeWitt super-Riemann surfaces, which also essentially follow from Corollary \ref{uniformization-cor} on uniformization for N=1 superanalytic DeWitt super-Riemann surfaces, along with results on ringed-space supermanifolds as in \cite{Manin1}.

In Section \ref{torus-section}, we first recall some preliminary results about complex tori and theta functions.  Then in Section \ref{torus-superconformal-section}, we study N=2 superconformal DeWitt super-Riemann surfaces over a fixed complex torus with certain restricted N=2 superconformal coordinate transformations.  In particular, we show that there is a doubly infinite family of distinct N=2 superconformal equivalence classes of compact genus-one N=2 superconformal DeWitt super-Riemann surfaces and that these distinct equivalence classes are characterized by theta functions associated to the torus up to type modulo the trivial theta functions, or equivalently by holomorphic line bundles over the underlying complex torus up to conformal equivalence.  Aspects of genus-zero supermanifolds have been studied previously in, for instance, \cite{Manin2}.

In Section \ref{nonhomo-section}, we reformulate our results in the ``nonhomogeneous coordinate system" as opposed to the ``homogeneous coordinate system" we had been using.  We use this reformulation to give some intuition as to why our uniformization theorems hold.  In particular, we discuss the infinitesimal N=1 and N=2 superconformal coordinate transformations, the presence of a representation of the affine unitary Lie algebra and an action of the $GL(1)$ loop group on the moduli spaces of genus-zero and genus-one N=2 superconformal super-Riemann surfaces.

Acknowledgments: The author thanks Yi-Zhi Huang, Liviu Nicolaescu, Jeffrey Rabin and Stephan Stolz  for insightful comments and discussions.   The author thanks Jeffrey Rabin as well as an anonymous referee for pointing out mistakes in preliminary versions of this paper.   The author also thanks the kind support and hospitality of the Max Planck Institute for Mathematics in Bonn, Germany. 

\section{Preliminaries: Superanalytic functions, N=2 superconformal functions and supermanifolds}\label{preliminaries-section}

\subsection{Superanalytic functions}

In this section, we recall the notion of Grassmann algebra and superanalytic function following, for instance \cite{B-n2moduli}, \cite{D}, \cite{Ro}.  

Let $\mathbb{C}$ denote the complex numbers, let $\mathbb{Z}$ denote the integers, and let $\mathbb{Z}_2$ denote the integers modulo 2.  For a $\mathbb{Z}_2$-graded vector space $V = V^0 \oplus V^1$, over $\mathbb{C}$, define the {\it sign function} $| \cdot |$ on the homogeneous subspaces of $V$ by $|v| = j$, for $v \in V^j$ and $j \in \mathbb{Z}_2$.  If $|v| = 0$, we say that $v$ is {\it even}, and if $|v| = 1$, we say that $v$ is {\it odd}.    A {\it superalgebra} is an (associative) algebra $A$ (with identity $1 \in A$), such that: (i) $A$ is a $\mathbb{Z}_2$-graded algebra; (ii) $ab = (-1)^{|a| \, |b|} ba$ for $a,b$ homogeneous in $A$.

Let $V$ be a complex vector space.  The exterior algebra generated by $V$, denoted  $\bigwedge (V)$, has the structure of a superalgebra.  Let $\mathbb{N}$ denote the nonnegative integers.  For $L \in \mathbb{N}$, fix $V_L$ to be  an $L$-dimensional vector space over $\mathbb{C}$ with basis $\{\zeta_1,\zeta_2, \ldots, \zeta_L\}$ such that $V_L \subset V_{L+1}$.  We denote $\bigwedge(V_L)$ by $\bigwedge_L$ and call this the {\it complex Grassmann algebra on $L$  generators}.  In other words, {}from now on we will consider the Grassmann algebras to have a fixed sequence of generators. Then $\bigwedge_L$ is the associative algebra over $\mathbb{C}$ with generators $1, \zeta_j$, for $j = 1, 2, \dots, L$, and with relations $1 \cdot \zeta_j = \zeta_j \cdot 1$, $\zeta_j \zeta_k = - \zeta_k \zeta_j$ and $\zeta_j^2 = 0$, for $j,k = 1,2, \dots, L$.    

Let 
\begin{eqnarray*}
J^0_L \! \!  &=&  \! \! \bigl\{ (j) = (j_1, j_2, \ldots, j_{2n}) \; | \; j_1 < j_2 < \cdots < j_{2n}, \; j_l \in \{1, 2, \dots, L\}, \; n \in \mathbb{N} \bigr\}, \\ 
J^1_L  \! \! &=&  \! \! \bigl\{(j) = (j_1, j_2, \ldots, j_{2n + 1}) \; | \; j_1 < j_2 < \cdots < j_{2n + 1}, \; j_l \in \{1, 2, \dots, L\}, \;  n \in \mathbb{N} \bigr\},
\end{eqnarray*}
and $J_L = J^0_L \cup J^1_L$.  Let $\mathbb{Z}_+$ denote the positive integers, and let
\begin{eqnarray*}
J^0_\infty \! \! &=& \! \! \bigl\{(j) = (j_1, j_2, \ldots, j_{2n})\; | \; j_1 < j_2 < \cdots < j_{2n}, \; j_l \in \mathbb{Z}_+, \; n \in \mathbb{N} \bigr\}, \\
J^1_\infty \! \! &=& \! \! \bigl\{(j) = (j_1, j_2, \ldots, j_{2n + 1})\; | \; j_1 < j_2 < \cdots < j_{2n + 1}, \; j_l \in \mathbb{Z}_+, \; n \in \mathbb{N} \bigr\},
\end{eqnarray*}
and $J_\infty = J^0_\infty \cup J^1_\infty$.  We use $J^0_*$, $J^1_*$, and $J_*$ to denote $J^0_L$ or $J^0_\infty$, $J^1_L$ or $J^1_\infty$, and $J_L$ or $J_\infty$, respectively.  Note that $(j) = (j_1,\dots,j_{2n})$ for $n = 0$ is in $J^0_*$, and we denote this element by $(\emptyset)$.  

The {\it infinite Grassmann algebra}, denoted  $\bigwedge_\infty$, is the superalgebra over $\mathbb{C}$ given by
\[\mbox{$\bigwedge_\infty$} = \Bigl\{ \sum_{(j) \in J_\infty} a_{(j)}\zeta_{j_{1}}\zeta_{j_{2}} \cdots \zeta_{j_{n}} \; \big\vert \; a_{(j)} \in \mathbb{C}, \; n \in \mathbb{N} \Bigr\}\] 
and with relations $1 \cdot \zeta_j = \zeta_j \cdot 1$, $\zeta_j \zeta_k = - \zeta_k \zeta_j$ and $\zeta_j^2 = 0$, for $j,k = 1, 2, \dots$.    We use the notation $\bigwedge_*$ to denote a Grassmann algebra, finite or infinite.  
The reason we take $\bigwedge_*$ to be over $\mathbb{C}$ is that we will be interested in complex supergeometry.  However, formally, we could just as well have taken $\mathbb{C}$ to be any field of characteristic zero.

The $\mathbb{Z}_2$-grading of $\bigwedge_*$ is given explicitly by
\begin{eqnarray*}
\mbox{$\bigwedge_*^0$} \! &=& \! \Bigl\{ \sum_{(j) \in J^0_*} a_{(j)}\zeta_{j_{1}}\zeta_{j_{2}} \cdots \zeta_{j_{2n}} \; \big\vert \; a_{(j)} \in \mathbb{C}, \; n \in \mathbb{N} \Bigr\}\\ 
\mbox{$\bigwedge_*^1$} \! &=& \! \Bigl\{ \sum_{(j) \in J^1_*} a_{(j)}\zeta_{j_{1}}\zeta_{j_{2}} \cdots \zeta_{j_{2n + 1}} \; \big\vert  \; a_{(j)} \in \mathbb{C}, \; n \in \mathbb{N}
\Bigr\} . 
\end{eqnarray*} 

We can also decompose $\bigwedge_*$ into {\it body}, $(\bigwedge_*)_B = \{ a_{(\emptyset)} \in \mathbb{C} \}$,  and {\it soul} 
\[(\mbox{$\bigwedge_*$})_S \; = \; \Bigl\{\sum_{(j) \in J_* \smallsetminus \{(\emptyset) \}} \! 
a_{(j)} \zeta_{j_1} \zeta_{j_2} \cdots \zeta_{j_n} \; \big\vert \; a_{(j)} \in \mathbb{C}  \Bigr\}\] 
subspaces such that $\bigwedge_* = (\bigwedge_*)_B \oplus (\bigwedge_*)_S$.  For $a \in \bigwedge_*$, we write $a = a_B + a_S$ for its body and soul decomposition.  We will use both notations $a_B$ and $a_{(\emptyset)}$ for the body of a supernumber $a \in \bigwedge_*$ interchangeably.  

For $n \in \mathbb{N}$, we introduce the notation $\bigwedge_{*>n}$ to denote a finite Grassmann algebra $\bigwedge_L$ with $L > n$ or an infinite Grassmann algebra.  We will use the corresponding index notations for the corresponding indexing sets $J^0_{*>n}, J^1_{*>n}$ and $J_{*>n}$.

Let $m,n \in \mathbb{N}$, and let $U$ be a subset of $(\bigwedge_*^0)^m \oplus (\bigwedge_*^1)^n$.  A $\bigwedge_*$-superfunction $H$ on $U$ in $(m,n)$-variables is given by
\begin{eqnarray*}
H: U &\longrightarrow& \mbox{$\bigwedge_*$} \\ 
(z_1,z_2,\dots,z_m ,\theta_1,\theta_2,\dots,\theta_n) &\mapsto& H(z_1,z_2,\dots,z_m ,\theta_1,\theta_2,\dots,\theta_n)
\end{eqnarray*}
where $z_k$, for $k = 1,\dots,m$, are even variables in $\bigwedge_*^0$ and $\theta_k$, for $k = 1,\dots,n$, are odd variables in $\bigwedge_*^1$.  If $H$ takes values only in $\bigwedge_*^0$, respectively in $\bigwedge_*^1$, we say that $H$ is an {\it even}, respectively {\it odd}, superfunction.  Let $f((z_1)_B,(z_2)_B,\dots,(z_m)_B)$ be a complex analytic function in $(z_k)_B$, for $k = 1,\dots,m$.  For $z_k \in \bigwedge_*^0$, and $k = 1,\dots,m$, define
\begin{multline}\label{more-than-one-variable}
f(z_1,z_2,\dots,z_m) = \! \sum_{l_1,\dots,l_m \in \mathbb{N}} \! \! \frac{(z_1)_S^{l_1} (z_2)_S^{l_2} \cdots
(z_m)_S^{l_m}}{l_1 ! l_2 ! \cdots l_m !} \biggl(\frac{\partial \; \;}{\partial (z_1)_B}\biggr)^{l_1} \biggl(\frac{\partial \; \;}{\partial (z_2)_B}\biggr)^{l_2} \\
  \cdots \biggl(\frac{\partial \; \;}{\partial (z_m)_B}\biggr)^{l_m}  f((z_1)_B,(z_2)_B,\dots,(z_m)_B) .
\end{multline}
Note that if $\bigwedge_*$ is finite, then Eq. (\ref{more-than-one-variable}) is a finite sum since $z^{L+1} = 0$ for $z \in \bigwedge_L$.   If $\bigwedge_*$ is infinite, then Eq. (\ref{more-than-one-variable}) is an infinite sum but is well defined since the coefficient of each basis element $\zeta_{j_1} \cdots \zeta_{j_m}$ is finite. 

Consider the projection
\begin{eqnarray}
\pi^{(m,n)}_B : (\mbox{$\bigwedge_{* > n-1}^0$})^m \oplus (\mbox{$\bigwedge_{* > n-1}^1$})^n & \longrightarrow & \mathbb{C}^m  \label{projection-onto-body}\\ 
(z_1,\dots,z_m, \theta_1,\dots,\theta_n) & \mapsto & ((z_1)_B,(z_2)_B,\dots,(z_m)_B) . \nonumber
\end{eqnarray}

\begin{defn}\label{superanal-definition}
Let $m, n \in \mathbb{N}$.  Let $U \subseteq (\bigwedge_{*>n-1}^0)^m \oplus (\bigwedge_{*>n-1}^1)^n$, and let $H$ be a $\bigwedge_{*>n-1}$-superfunction in $(m,n)$-variables defined on $U$. Then $H$ is said to be {\em superanalytic} if $H$ is of the form
\begin{equation}
H(z_1,z_2,\dots,z_m ,\theta_1,\theta_2,\dots,\theta_n) = \sum_{ (j) \in J_n} \theta_{j_1} \cdots \theta_{j_{l}} f_{(j)}(z_1,z_2,\dots,z_m) , 
\end{equation} 
where each $f_{(j)}$ is of the form 
\begin{equation}\label{restrict-coefficients}
f_{(j)}(z_1,z_2,\dots,z_m) = \sum_{(k) \in J_{ * - n}} f_{(j),(k)}(z_1,z_2,\dots,z_m) \zeta_{k_1}\zeta_{k_2} \cdots \zeta_{k_{s}}, 
\end{equation} 
and each $f_{(j),(k)}((z_1)_B,(z_2)_B,\dots,(z_m)_B)$ is analytic  in $(z_l)_B$, for $l = 1,\dots,m$ and $((z_1)_B,(z_2)_B,\dots,(z_m)_B) \in  U_B = \pi_B^{(m,n)} (U) \subseteq \mathbb{C}^m$.   
\end{defn}

We require the even and odd variables to be in $\bigwedge_{* > n-1}$, and we restrict the coefficients of the $f_{(j),(k)}$'s to be in $\bigwedge_{* - n} \subseteq \bigwedge_{*>n-1}$ in order for the partial derivatives with respect to each of the $n$ odd variables to be well defined and for multiple partials to be well defined (cf. \cite{D}, \cite{B-memoir}, \cite{Ro}, \cite{B-n2moduli}).    

\begin{rema}{\em   In the language of \cite{Ro}, these superanalytic $\bigwedge_{*>n-1}$-superfunctions just defined are called {\it $GC^\omega$ functions on $\mathbb{C}_S^{m,n}$} if $\bigwedge_* = \bigwedge_\infty$, and are called {\it $GHC^\omega$ functions on 
$\mathbb{C}_{S[L]}^{m,n}$} if $\bigwedge_* = \bigwedge_L$.   In the language of \cite{Ro}, the class of {\it $HC^\omega$ functions on $\mathbb{C}_S^{m,n}$} are those superanalytic $\bigwedge_\infty$-superfunctions in $(m,n)$-variables for which the coefficient functions $f_{(j)}$ restricted to $\mathbb{C}^m$ take values in $\mathbb{C}$ rather than more generally in $\bigwedge_\infty$.  In other words, the  $HC^\omega$ functions on $\mathbb{C}_S^{m,n}$ are the subclass of $GC^\omega$ functions on $\mathbb{C}_S^{m,n}$  for which $f_{(j), (k)} = 0$ if $(k) \neq (\emptyset)$, and thus are a proper subclass of the functions we call superanalytic $\bigwedge_{*>n-1}$-superfunctions in $(m,n)$-variables. 
Similarly, the class of {\it $HC^\omega$ functions on $\mathbb{C}_{S[L]}^{m,n}$} in \cite{Ro} are those superanalytic $\bigwedge_{L>n-1}$-superfunctions in $(m,n)$-variables for which the coefficient functions $f_{(j)}$ restricted to $\mathbb{C}^m$ take values in $\mathbb{C}$ rather than more generally in $\bigwedge_{L-n}$, and thus are also a proper subclass of the functions we call superanalytic $\bigwedge_{*>n-1}$-superfunctions in $(m,n)$-variables.   Often the ringed-space approach to supermanifolds restricts to functions in this subclass $HC^\omega$ of the more general class $GC^\omega$.}
\end{rema}

We define the {\it DeWitt topology on $(\bigwedge_{* > n-1}^0)^m \oplus (\bigwedge_{* > n-1}^1)^n$} by letting a subset $U$ of $(\mbox{$\bigwedge_{*> n-1}^0$})^m \oplus (\mbox{$\bigwedge_{* > n-1}^1$})^n$ be an open set in the DeWitt topology if and only if $U = (\pi^{(m,n)}_B)^{-1} (V)$ for some open set $V \subseteq \mathbb{C}^m$.  Note that the natural domain of a superanalytic $\bigwedge_{* > n-1}$-superfunction in $(m,n)$-variables is an open set in the DeWitt topology.  

Let $(\bigwedge_*)^\times$ denote the set of invertible elements in $\bigwedge_*$.  Then $(\mbox{$\bigwedge_*$})^\times = \{a \in \mbox{$\bigwedge_*$} \; | \; a_B \neq 0 \},$
since $\frac{1}{a} = \frac{1}{a_B + a_S} = \sum_{n \in \mathbb{N}} \frac{(-1)^n a_S^n}{a_B^{n + 1}}$
is well defined if and only if $a_B \neq 0$. 

\begin{rema}\label{H_L-remark}
{\em Recall that $\bigwedge_L \subset \bigwedge_{L+1}$ for $L \in \mathbb{N}$, and note that {}from (\ref{more-than-one-variable}), any superanalytic $\bigwedge_L$-superfunction, $H_L$, in $(m,n)$-variables for $L \geq n$ can naturally be extended to a superanalytic $\bigwedge_{L'}$-superfunction in $(m,n)$-variables for $L'>L$ and hence to a superanalytic $\bigwedge_\infty$-superfunction.  Conversely, if $H_{L'}$ is a superanalytic $\bigwedge_{L'}$-superfunction (or $\bigwedge_\infty$-superfunction) in $(m,n)$-variables for $L' >n$, then we can restrict $H_{L'}$ to a superanalytic $\bigwedge_L$-superfunction for $L'> L\geq n$ by restricting $(z_1,\dots,z_m,\theta_1,\dots,\theta_n) \in 
(\bigwedge_L^0)^m \oplus (\bigwedge_L^1)^n$ and setting $f_{(j),(k)} \equiv 0$ if $(k) \notin J_{L-n}$. 
If $H_{L'}$ satisfies $f_{(j), (k)} \equiv 0$ for $(k) \notin J_{L-n}$, then restriction to $\bigwedge_L$ and then extension to $\bigwedge_{L'}$ results in the identity mapping, i.e., leaves $H_{L'}$ unchanged.  Thus any superanalytic function over $\bigwedge_{L'}$ in $(m,n)$-variables with coefficient functions $f_{(j), (k)} = 0$ for $(k) \notin J_{L-n}$,  for $L\leq L'$, can be thought of as a functor from the category of Grassmann algebras $\bigwedge_*$ with $*\geq L+n$ to superanalytic functions over $\bigwedge_*$ in $(m,n)$-variables (cf. \cite{Schwarz1984}, \cite{KS}).   
}
\end{rema}

\subsection{Subclasses of $(1,1)$- and $(1, 2)$-superfunctions and N=1 and N=2 superconformal functions}\label{superconformal-section}

For the purposes of this paper, our focus will be on superanalytic $\bigwedge_{*>n-1}$-superfunctions in $(1,n)$-variables for $n \leq 2$, i.e., the case of one even variable and one or two odd variables.  Here we introduce some notation for some of the subclasses of $(1,n)$-superfunctions for $n = 1,2$.  In particular, distinguishing these subclasses will be useful both for stating our results, and for relating these results to some of the results from the ringed-space approach which is more restrictive than the general approach we take throughout this paper.   In addition, we recall the notions of N=1 and N=2 superconformal functions following, for instance,  \cite{CR}, \cite{DRS}, \cite{B-thesis}, \cite{B-memoir}, \cite{B-n2moduli}.

A superanalytic  $(1,1)$-superfunction $H(z,\theta) = (\tilde{z}, \tilde{\theta})$ from a DeWitt open neighborhood in $\bigwedge_{*>0}^0 \oplus \bigwedge_{*>0}^1$ to $\bigwedge_{*>0}^0 \oplus \bigwedge_{*>0}^1$ is of the form
\begin{eqnarray}
\tilde{z} &=& f(z) + \theta \xi(z)  \label{N=1-superanalytic1} \\
\tilde{\theta} &=& \psi (z) + \theta g(z) \label{N=1-superanalytic2}
\end{eqnarray}
for $f,g$ even and $\xi, \psi$ odd superanalytic $(1,0)$-superfunctions in $z$.  We will call this class of functions {\it $N=1$ superanalytic functions}, and denote the class of such functions by $\mathcal{G}_{*>0} (1)$.  In the ringed-space approach, one does not consider odd functions of an even variable such as $\xi$ and $\psi$.  Let $\mathcal{H}_{*>0}(1)$, be the subclass of $\mathcal{G}_{*>0} (1)$ consisting of $(1,1)$-superfunctions $H(z,\theta) = (\tilde{z}, \tilde{\theta})$ of the form (\ref{N=1-superanalytic1})-(\ref{N=1-superanalytic2}), where $\psi(z)$ and $\xi(z)$ and the soul portion of $f(z)$ are identically zero.  That is, $f(z_B)$ is a complex analytic function on an open subset of $\mathbb{C}$.  Finally, let $\mathcal{C}_{*>0}(1)$ be the subclass of $\mathcal{H}_{*>0}(1)$ consisting of $\mathcal{H}_{*>0}(1)$ functions such that $g(z_B)$ is a complex analytic functions on an open subset of $\mathbb{C}$.   Thus we have
$\mathcal{C}_{*>0}(1) \subseteq \mathcal{H}_{*>0}(1) \subseteq \mathcal{G}_{*>0}(1)$, and the only time equality holds is in the case when $* = 1$,  i.e., when we are working over the Grassmann algebra $\bigwedge_1$.  In this case, $\mathcal{C}_1 (1) = \mathcal{H}_1(1) = \mathcal {G}_1(1)$.   In addition, by considering the natural extension of a superanalytic function over $\bigwedge_L$ to a superanalytic function over $\bigwedge_{L'}$ for $L< L'$ following Remark \ref{H_L-remark}, we have the inclusions $\mathcal{G}_{L}(1) \subset \mathcal{G}_{L'}(1) \subseteq \mathcal{G}_\infty(1)$. 

Similarly, a superanalytic  $(1,2)$-superfunction $H(z,\theta_1, \theta_2) = (\tilde{z}, \tilde{\theta}_1, \tilde{\theta}_2)$ from a DeWitt open neighborhood in $\bigwedge_{*>1}^0 \oplus (\bigwedge_{*>1}^1)^2$ to $\bigwedge_{*>1}^0 \oplus (\bigwedge_{*>1}^1)^2$ is of the form
\begin{eqnarray}
\tilde{z} &=& f(z) + \theta_1 \xi_1(z) + \theta_2 \xi_2(z) + \theta_1 \theta_2 g(z) \label{N=2-superanalytic1} \\
\tilde{\theta}_j &=& \psi_j (z) + \theta_1 g_j(z) + \theta_2 h_j(z) + \theta_1 \theta_2 \varphi_j(z)  \label{N=2-superanalytic2}
\end{eqnarray}
for $j = 1,2$, and $f,g, g_j, h_j$ even and $\xi_j, \psi_j, \varphi_j$ odd superanalytic $(1,0)$-super-functions in $z$. We will call this class of functions {\it $N=2$ superanalytic functions}, and denote the class of such functions by $\mathcal{G}_{*>1} (2)$.  Again, since in the ringed-space approach, one does not consider odd functions of an even variable, we let $\mathcal{H}_{*>1}(2)$, be the subclass of $\mathcal{G}_{*>1} (2)$ consisting of $(1,2)$-superfunctions $H(z,\theta_1, \theta_2) = (\tilde{z}, \tilde{\theta}_1, \tilde{\theta}_2)$ of the form (\ref{N=2-superanalytic1})-(\ref{N=2-superanalytic2}), where $\psi_j(z)$, $\xi_j(z)$ and $\phi_j$, for $j =1,2$, are identically zero and $f(z_B)$, is a complex analytic function on an open subset of $\mathbb{C}$.  Finally, let $\mathcal{C}_{*>1}(2)$ be the subclass of $\mathcal{H}_{*>1}(2)$ consisting of $\mathcal{H}_{*>1}(2)$ functions such that $g(z_B)$, and $g_j(z_B)$ and $h_j(z_B)$, for $j =1,2$, are complex analytic functions on an open subset of $\mathbb{C}$.   Thus we have
$\mathcal{C}_{*>1}(2) \subseteq \mathcal{H}_{*>1}(2) \subseteq \mathcal{G}_{*>1}(2)$, and the only time equality holds is in the case when $* = 2$,  i.e., when we are working over the Grassmann algebra $\bigwedge_2$.  In this case, $\mathcal{C}_2 (2) = \mathcal{H}_2(2) = \mathcal {G}_2(2)$.   In addition, we have the inclusions $\mathcal{G}_{L}(2) \subset \mathcal{G}_{L'}(2) \subseteq \mathcal{G}_\infty(2)$, for $0< L < L'$.

An N$= \! n$ superconformal function is a superanalytic $(1,n)$-superfunction 
on a DeWitt open subset of $\bigwedge_{*>n-1}^0 \oplus (\bigwedge_{*>n-1}^1)^n$ to $\bigwedge_{*>n-1}^0 \oplus (\bigwedge_{*>n-1}^1)^n$ that transforms the superderivations $D_j =  \frac{\partial}{\partial \theta_j} + \theta_j \frac{\partial}{\partial z}$, for $j = 1,\dots n$, in a certain ``homogeneous" way.   

In particular, an N=1 superanalytic  function $H(z,\theta) = (\tilde{z}, \tilde{\theta})$ of the form (\ref{N=1-superanalytic1})-(\ref{N=1-superanalytic2}) transforms $D = \frac{\partial}{\partial \theta} + \theta \frac{\partial}{\partial z}$  to $\tilde{D} = \frac{\partial}{\partial \tilde{\theta}} + \tilde{\theta} \frac{\partial}{\partial \tilde{z}}$ by  $D = (D \tilde{\theta}) \tilde{D} + (D\tilde{z} - \tilde{\theta} D \tilde{\theta}) \tilde{D}^2$.   We define $H$ to be {\it N=1 superconformal} if $H$ transforms $D$ homogeneously of degree one. That is, if $H(z, \theta) = (\tilde{z}, \tilde{\theta})$ satisfies $D\tilde{z} - \tilde{\theta} D \tilde{\theta} = 0$.  This is equivalent to $H$ having the form 
\begin{eqnarray}
\tilde{z} &=& f(z) + \theta g(z) \psi(z)  \label{N=1-superconformal-condition1} \\
\tilde{\theta} &=& \psi (z) + \theta g(z) \label{N=1-superconformal-condition2}
\end{eqnarray}
for $f,g$ even and $\psi$ odd superanalytic $(1,0)$-superfunctions in $z$, satisfying the condition 
\begin{equation}\label{N=1superconformal-condition4}
f'(z)  = g(z) g(z)-  \psi(z)\psi'(z), \qquad \mbox{i.e., \ \ $g^2(z) = f'(z) + \psi(z) \psi'(z)$.}
\end{equation}
Thus  an N=1 superconformal function $H$ is uniquely determined by the superanalytic functions $f(z)$ and $\psi(z)$ and a choice of square root for (\ref{N=1superconformal-condition4}).

In the N=2 superconformal setting there are generally two different coordinate systems commonly used.  We will first work in what we call  the ``nonhomogeneous" coordinate setting (see \cite{B-n2moduli} and \cite{B-axiomatic}), and then translate to the ``homogeneous" coordinate system. 

An N=2 superanalytic function $H(z,\theta_1, \theta_2) = (\tilde{z}, \tilde{\theta}_1, \tilde{\theta}_2)$  of the form (\ref{N=2-superanalytic1})-(\ref{N=2-superanalytic2}) transforms $D_1$ and $D_2$ by
\begin{eqnarray*}
D_1  &=& (D_1 \tilde{\theta}_1)\tilde{D}_1 + (D_1 \tilde{\theta}_2 ) \tilde{D}_2 + \left( D_1 \tilde{z} - \tilde{\theta}_1D_1\tilde{\theta}_1 -   \tilde{\theta}_2D_1 \tilde{\theta}_2\right) \tilde{D}_1^2 \\
D_2 &=& (D_2 \tilde{\theta}_1 )\tilde{D}_1  + (D_2 \tilde{\theta}_2 ) \tilde{D}_2 + \left( D_2 \tilde{z} - \tilde{\theta}_1 D_2 \tilde{\theta}_1  -   \tilde{\theta}_2  D_2 \tilde{\theta}_2  \right) \tilde{D}_2^2.
\end{eqnarray*}
We define  $H$ to be {\it N=2 superconformal} if it transforms $D_1$, respectively $D_2$, as $D_1 = (D_1 \tilde{\theta}_1) \tilde{D}_1 + (D_1\tilde{\theta}_2 ) \tilde{D}_2 =  (D_1 \tilde{\theta}_1) \tilde{D}_1 - (D_2\tilde{\theta}_1 ) \tilde{D}_2$, respectively 
$D_2 = (D_2 \tilde{\theta}_1 )\tilde{D}_1  + (D_2 \tilde{\theta}_2 ) \tilde{D}_2 = (D_2 \tilde{\theta}_1 )\tilde{D}_1  + (D_1 \tilde{\theta}_1 ) \tilde{D}_2$.   That is $H$ must satisfy
\begin{eqnarray}
D_1 \tilde{\theta}_1  - D_2 \tilde{\theta}_2 \ = \ D_1 \tilde{\theta}_2 +  D_2 \tilde{\theta}_1  &=& 0 \label{inhomo-condition1} \\
D_1 \tilde{z} -  \tilde{\theta}_1 D_1 \tilde{\theta}_1 -  \tilde{\theta}_2D_1 \tilde{\theta}_2  &=& 0 \\
D_2 \tilde{z} -    \tilde{\theta}_1 D_2 \tilde{\theta}_1 -  \tilde{\theta}_2 D_2 \tilde{\theta}_2  &=& 0 . \label{inhomo-condition4}
\end{eqnarray}
These conditions (\ref{inhomo-condition1})--(\ref{inhomo-condition4}) imply that an N=2 superconformal function  in the nonhomogeneous coordinate system $H(z, \theta_1, \theta_2) = (\tilde{z}, \tilde{\theta}_1, \tilde{\theta}_2)$ is of the form
\begin{eqnarray}
 \tilde{z} \! \! &=& \! \! f(z) + \theta_1 (g_1(z) \psi_1(z) +   g_2(z)\psi_2(z))  + \theta_2 (g_1(z)\psi_2(z) - g_2(z) \psi_1(z) )  \label{nonhomo-superconformal-condition1}  \\
& &  \quad  -  \theta_1 \theta_2 ( \psi_1(z)  \psi_2(z))'   \nonumber \\
\tilde{\theta}_1 \! \! &=& \! \!  \psi_1 (z) + \theta_1  g_1(z) - \theta_2 g_2(z) + \theta_1 \theta_2  (\psi_2)'(z)  \\
\qquad \ \ \tilde{\theta}_2 \! \! &=& \! \!  \psi_2(z) + \theta_1 g_2(z) + \theta_2 g_1(z) - \theta_1 \theta_2 (\psi_1)'(z) ,
\end{eqnarray}
satisfying
\begin{equation}
f'(z) =   g_1^2(z) +   g_2^2(z)  -  \psi_1(z) (\psi_1)'(z)  -    \psi_2(z) (\psi_2)'(z)  , \label{nonhomo-superconformal-condition4} 
\end{equation}
for even superanalytic $(1,0)$-superfunctions $f, g_1$ and $g_2$ and odd superanalytic $(1,0)$-superfunctions $\psi_1, \psi_2$.

It will be convenient for us to work in the ``homogeneous" coordinate system denoted by even variable $z$ and odd variables $\theta^+$ and $\theta^-$, where 
\begin{equation}\label{transform-nonhomo-homo}
\theta^\pm  = \frac{1}{\sqrt{2}} \left( \theta_1 \pm  i\theta_2 \right),
\end{equation}
or equivalently
\begin{equation}\label{transform-homo-inhomo}
\theta_1 = \frac{1}{\sqrt{2}} \left( \theta^+ + \theta^- \right) \qquad \mbox{and} \qquad
\theta_2 = - \frac{i}{\sqrt{2}} \left( \theta^+ - \theta^- \right).
\end{equation}
This is a standard transformation in N=2 superconformal field theory (cf. \cite{DRS}, \cite{B-n2moduli}, \cite{B-axiomatic}), however the nomenclature ``homogeneous" for the $(z, \theta^+, \theta^-)$ coordinate system and ``nonhomogeneous" for the $(z, \theta_1, \theta_2)$ coordinate system was first introduced by the author in \cite{B-n2moduli}.  In Remark \ref{transform-D-remark} we give some reasons for this terminology.   Other reasons for this nomenclature involve the algebra of infinitesimal N=2 superconformal transformations and are discussed in Remark \ref{homo-u(1)-remark} below as well as in \cite{B-n2moduli}.

We have that
\begin{equation*}
\frac{\partial}{\partial \theta_1} =  \frac{1}{\sqrt{2}} \Bigl(\frac{\partial}{\partial \theta^+} + \frac{\partial}{\partial \theta^-} \Bigr) \qquad \mathrm{and} \qquad 
\frac{\partial}{\partial \theta_2} = \frac{i}{\sqrt{2}} \Bigl( \frac{\partial}{\partial \theta^+} - \frac{\partial}{\partial \theta^-} \Bigr) 
\end{equation*}
or equivalently
\begin{equation*}
\frac{\partial}{\partial \theta^\pm} =  \frac{1}{\sqrt{2}} \Bigl(\frac{\partial}{\partial \theta_1} \mp i \frac{\partial}{\partial \theta_2} \Bigr).
\end{equation*}

Define
\begin{equation}
D^\pm = \frac{\partial}{\partial \theta^\pm} + \theta^\mp \frac{\partial}{\partial z} = \frac{1}{\sqrt{2}} ( D_1 \mp i  D_2 ).
\end{equation}
Note that
\begin{eqnarray}
[D^\pm, D^\pm] & = & 2(D^\pm)^2 \ = \ 0\\ 
\left[D^+,D^-\right] &=& D^+D^- + D^-D^+ \ = \ 2 \frac{\partial}{\partial z} .
\end{eqnarray}

Let $H(z,\theta^+,\theta^-) = (\tilde{z},\tilde{\theta}^+,\tilde{\theta}^-)$ be an  N=2 superanalytic function in the homogeneous coordinate system from a DeWitt open neighborhood in $\bigwedge_{*>1}^0 \oplus (\bigwedge_{*>1}^1)^2$ to $\bigwedge_{*>1}^0 \oplus (\bigwedge_{*>1}^1)^2$,  i.e., $\tilde{z}$ is an even superanalytic $(1,2)$-superfunction and $\tilde{\theta}^\pm$ are odd superanalytic $(1,2)$-superfunctions.   Then $D^+$ and $D^-$ transform under $H(z,\theta^+,\theta^-)$ by
\begin{equation}\label{transform-Dplus}
D^\pm = (D^\pm\tilde{\theta}^\pm)\tilde{D}^\pm + (D^\pm\tilde{\theta}^\mp) \frac{\partial}{\partial \tilde{\theta}^\mp} + (D^\pm\tilde{z} - \tilde{\theta}^\mp D^\pm\tilde{\theta}^\pm) \frac{\partial}{\partial \tilde{z}} .
\end{equation}

An N=2  superconformal function $H$ in the homogeneous coordinate system transforms $D^+$ and $D^-$ homogeneously of degree one.   That is, $H$ transforms $D^\pm$ by non-zero superanalytic functions times $\tilde{D}^\pm$, respectively.   Since such a superanalytic function $H(z,\theta^+,\theta^-) = (\tilde{z}, \tilde{\theta}^+,\tilde{\theta}^-)$  transforms $D^+$ and $D^-$ according to  (\ref{transform-Dplus}), $H$ is superconformal if and only if, in addition to being superanalytic, $H$ satisfies
\begin{eqnarray}
D^\pm \tilde{\theta}^\mp &=& 0, \label{basic-superconformal-condition1} \\
D^\pm\tilde{z} - \tilde{\theta}^\mp D^\pm\tilde{\theta}^\pm &=& 0, \label{basic-superconformal-condition4}
\end{eqnarray}
for $D^\pm \tilde{\theta}^\pm$ not identically zero, thus transforming $D^\pm$ by $D^\pm = (D^\pm \tilde{\theta}^\pm)\tilde{D}^\pm$.   These conditions imply that we can write $H(z,\theta^+,\theta^-) = (\tilde{z}, \tilde{\theta}^+, \tilde{\theta}^-)$ as
\begin{eqnarray}
\qquad \qquad \tilde{z} &=& f(z) + \theta^+ g^+(z) \psi^-(z) + \theta^- g^-(z) \psi^+(z) + \theta^+ \theta^- (\psi^+(z)\psi^-(z))' \label{superconformal-condition1} \\
\tilde{\theta}^\pm &=& \psi^\pm(z) + \theta^\pm g^\pm(z) \pm \theta^+ \theta^-(\psi^\pm)'(z) \label{superconformal-condition3}
\end{eqnarray}
for $f$, $g^\pm$ even and $\psi^\pm$ odd superanalytic $(1,0)$-superfunctions in $z$, satisfying the condition 
\begin{equation}\label{superconformal-condition4}
f'(z) \; = \; (\psi^+)'(z)\psi^-(z) - \psi^+(z) (\psi^-)'(z) + g^+(z) g^-(z) ,
\end{equation}
and we also require that $D^+ \tilde{\theta}^+$ and $D^- \tilde{\theta}^-$ not be identically zero.   Thus  an N=2 superconformal function $H$ is uniquely determined by the superanalytic functions $f(z)$, $\psi^\pm(z)$,  and $g^\pm(z)$ satisfying the condition (\ref{superconformal-condition4}). 

Note that, transforming between homogeneous and nonhomogeneous coordinate systems for an N=2 superconformal functions given by (\ref{superconformal-condition1})--(\ref{superconformal-condition4}), or equivalently  (\ref{nonhomo-superconformal-condition1})--(\ref{nonhomo-superconformal-condition4}),  we have that $\psi^\pm (z) = \frac{1}{\sqrt{2}} (\psi_1(z) \pm i \psi_2(z))$ and $g^\pm(z) = g_1(z) \pm i g_2(z)$.

\begin{rema}\label{transform-D-remark}
{\em {}From the properties derived above for an N=2 superconformal function in the nonhomogeneous coordinate system, we see one of the reasons for our terminology.  Namely,  that  in the nonhomogeneous coordinate system an N=2 superconformal function does not transform the superderivations $D_1$ and $D_2$, respectively, homogeneously of degree one.  Instead it transforms them as  $D_1 = (D_1 \tilde{\theta}_1) \tilde{D}_1 + (D_1 \tilde{\theta}_2 ) \tilde{D}_2$ and $D_2 = (D_2 \tilde{\theta}_1 )\tilde{D}_1  + (D_2 \tilde{\theta}_2 ) \tilde{D}_2$, respectively -- unlike the homogeneous nature of the transformation of $D^\pm$ under an N=2 superconformal function in the homogeneous coordinates.  In the latter case the superderivations transform homogeneously as $D^\pm = (D^\pm \tilde{\theta}^\pm) \tilde{D}^\pm$.}
\end{rema}

\subsection{Complex DeWitt supermanifolds and N=2 superconformal DeWitt super-Riemann surfaces}\label{supermanifolds-section}

A {\em DeWitt $(m,n)$-dimensional supermanifold over $\bigwedge_*$} is a topological space $X$ with a countable basis which is locally homeomorphic to an open subset of $(\bigwedge_*^0)^m \oplus (\bigwedge_*^1)^n$ in the DeWitt topology.  A {\em DeWitt $(m,n)$-chart on $X$ over $\bigwedge_*$} is 
a pair $(U, \Omega)$ such that $U$ is an open subset of $X$ and $\Omega$ is a homeomorphism of $U$ onto an open subset of $(\bigwedge_*^0)^m \oplus (\bigwedge_*^1)^n$ in the DeWitt topology.  A {\em superanalytic atlas of DeWitt $(m,n)$-charts on $X$ over $\bigwedge_{* > n-1}$} is a family of charts $\{(U_{\alpha}, \Omega_{\alpha})\}_{\alpha \in A}$ satisfying 
 
(i) Each $U_{\alpha}$ is open in $X$, and $\bigcup_{\alpha \in A} U_{\alpha} = X$. 

(ii) Each $\Omega_{\alpha}$ is a homeomorphism {}from $U_{\alpha}$ to a (DeWitt) open set in $(\bigwedge_{* > n-1}^0)^m \oplus (\bigwedge_{* > n-1}^1)^n$, such that $\Omega_{\alpha} \circ \Omega_{\beta}^{-1}: \Omega_{\beta}(U_\alpha \cap U_\beta) \longrightarrow \Omega_{\alpha}(U_\alpha \cap U_\beta)$ is superanalytic for all non-empty $U_{\alpha} \cap U_{\beta}$, i.e., $\Omega_{\alpha} \circ \Omega_{\beta}^{-1} = (\tilde{z}_1,\dots, \tilde{z}_m, \tilde{\theta}_1,\dots,\tilde{\theta}_n)$ where $\tilde{z}_i$ 
is an even superanalytic $\bigwedge_{* > n-1}$-superfunction in $(m,n)$-variables for $i = 1,\dots,m$, and $\tilde{\theta}_j$ is an odd superanalytic $\bigwedge_{* >n-1}$-superfunction in $(m,n)$-variables 
for $j = 1,\dots,n$.

Such an atlas is called {\em maximal} if, given any chart $(U, \Omega)$ such that
\[\Omega \circ \Omega_{\beta}^{-1} : \Omega_{\beta} (U \cap U_\beta) \longrightarrow \Omega (U \cap U_\beta)\] 
is a superanalytic homeomorphism for all $\beta$, then $(U, \Omega) \in \{(U_{\alpha}, \Omega_{\alpha})\}_{\alpha \in A}$.
 
A {\em DeWitt $(m,n)$-superanalytic supermanifold over $\bigwedge_{* > n-1}$} is a DeWitt $(m,n)$-dimensional supermanifold $M$ together with a maximal superanalytic atlas of DeWitt $(m,n)$-charts over $\bigwedge_{* > n-1}$.  

Given a DeWitt $(m,n)$-superanalytic supermanifold $M$ over $\bigwedge_{* > n-1}$, define an equivalence relation $\sim$ on M by letting $p \sim q$ if and only if there exists $\alpha \in A$ such that $p,q \in U_\alpha$ and $\pi_B^{(m,n)} (\Omega_\alpha (p)) = \pi_B^{(m,n)} (\Omega_\alpha (q))$,
where $\pi_B^{(m,n)}$ is the projection given by (\ref{projection-onto-body}).  Let $p_B$ denote the equivalence class of $p$ under this equivalence relation.  Define the {\it body} $M_B$ of $M$ to be the 
$m$-dimensional complex manifold with analytic structure given by the coordinate charts $\{((U_\alpha)_B, (\Omega_\alpha)_B) \}_{\alpha \in A}$ where $(U_\alpha)_B = \{ p_B \; | \; p \in U_\alpha \}$, and $(\Omega_\alpha)_B : (U_\alpha)_B \longrightarrow \mathbb{C}^m$ is given by $(\Omega_\alpha)_B (p_B) = \pi_B^{(m,n)} \circ \Omega_\alpha (p)$. We define the genus of $M$ to be the genus of $M_B$.

Note that $M$ is a complex fiber bundle over the complex manifold $M_{B}$; the fiber is the complex vector space $(\bigwedge_{* > n-1}^0)_S^m \oplus (\bigwedge_{* > n-1}^1)^n$.  This bundle is not in general a vector bundle since the transition functions are in general nonlinear.  

\begin{rema}\label{functorial-remark2}
{\em  Just as a single superanalytic function over a certain Grassmann algebra can be thought of as a functor from a (sub)category of Grassmann algebras to superanalytic functions over any one of these Grassmann algebras (see Remark \ref{H_L-remark}), so can a DeWitt $(m,n)$-superanalytic supermanifold over a certain Grassmann algebra be thought of as a functor from a (sub)category of Grassmann algebras to DeWitt $(m,n)$-superanalytic supermanifolds over any one of these Grassmann algebras.  Let $M$ be an $(m,n)$-superanalytic supermanifold over $\bigwedge_{L'}$ for $L'\geq n$ with coordinate atlas given by $\{(U_\alpha, \Omega_\alpha)\}_{\alpha \in A}$.  If the coordinate transition functions for $M$ are such that the coefficient functions $f_{(j), (k)} \equiv 0$ for $(k) \notin J_{L-n}$ for some $L'>L \geq n$, then the submanifold of $M$ given by $\bigcup_{\alpha \in A} \Omega_\alpha^{-1} ((\Omega_\alpha (U_\alpha))_B \times (\bigwedge_L^0)_S^m \times (\bigwedge^1_L)^n )$ 
is naturally a $(m,n)$-superanalytic supermanifold over $\bigwedge_L$.  Moreover, if $M_1$ and $M_2$ are $(m,n)$-superanalytic supermanifolds over $\bigwedge_{L'}$ which result in the same submanifold under this restriction from $\bigwedge_{L'}$ to $\bigwedge_L$, then $M_1 = M_2$.  Thus there is a natural and unique extension of any $(m,n)$-superanalytic supermanifold over $\bigwedge_L$ to a $(m,n)$-superanalytic supermanifold over $\bigwedge_*$ for $*>L$. 
}
\end{rema}

For any DeWitt $(1,n)$-superanalytic supermanifold $M$, its body $M_{B}$ is a Riemann surface.   An {\em N$= \! n$ superconformal DeWitt super-Riemann surface over $\bigwedge_{*>n-1}$, for $n=1,2$,} is a DeWitt $(1,n)$-superanalytic supermanifold over $\bigwedge_{*>n-1}$ with coordinate atlas $\{(U_{\alpha}, \Omega_{\alpha})\}_{\alpha \in A}$ such that the coordinate transition functions $\Omega_{\alpha} \circ \Omega_{\beta}^{-1}$ in addition to being superanalytic are also N$= \! n$ superconformal for all non-empty $U_{\alpha} \cap U_{\beta}$.

Since the condition that the coordinate transition functions be N$= \! n$ superconformal instead of merely superanalytic is such a strong condition (unlike in the nonsuper case), we again stress the distinction between an N$= \! n$ superanalytic DeWitt super-Riemann surface which has {\it superanalytic} transition functions versus an N$= \! n$ superconformal DeWitt super-Riemann surface which has N$= \! n$ {\it superconformal} transition functions.  In the literature one will find the term ``super-Riemann surface" or ``Riemannian supermanifold" used for both merely superanalytic structures (cf. \cite{D}) and for superconformal structures (cf. \cite{Fd}, \cite{CR}).   Throughout this paper, we will be mainly dealing with N=2 superconformal super-Riemann surfaces and N=1 superanalytic super-Riemann surfaces. 

\begin{rema}\label{DeWitt-v-ringed-space-remark} {\em In general, the transition functions for an N$= \! n$ superanalytic DeWitt super-Riemann surface, for $n = 1,2$, are $\mathcal{G}_{*>n-1}(n)$ functions.  If however, the functions are in the subclass of $\mathcal{H}_{*>n-1}(n)$ functions, or $\mathcal{C} _{*>n-1}(n)$ functions, then we will call such a super-Riemann surface a $\mathcal{H}_{*>n-1}(n)$-supermanifold or $\mathcal{C} _{*>n-1}(n)$-supermanifold, respectively.   These subclasses of DeWitt supermanifolds, $\mathcal{H}_{*>n-1}(n)$-supermanifold or $\mathcal{C} _{*>n-1}(n)$-supermanifold, are those that are equivalent to the supermanifolds studied in the ringed-space approach in, for instance, \cite{Manin1}, \cite{Rothstein}, \cite{Vaintrob}.  }
\end{rema} 

Let $M_1$ and $M_2$ be N$= \! n$ superanalytic DeWitt super-Riemann surfaces, for $n = 1,2$, with coordinate atlases $\{(U_{\alpha},$ $\Omega_{\alpha})\}_{\alpha \in A}$ and $\{(V_{\beta}, \Xi_{\beta})\}_{\beta \in B}$, respectively.  A map $F: M_1 \longrightarrow M_2$ is said to be {\it N$= \! n$ superanalytic}  if $\Xi_\beta \circ F \circ \Omega_\alpha^{-1}: \Omega_\alpha (U_\alpha \cap F^{-1} (V_\beta)) \longrightarrow \Xi_\beta(V_\beta)$ is N$= \! n$ superanalytic for all $\alpha \in A$ and $\beta \in B$ with $U_\alpha \cap F^{-1} (V_\beta) \neq \emptyset$.  If in addition, $F$ is bijective, then we say that $M_1$ and $M_2$ are {\it N$= \! n$ superanalytically equivalent}.  By {\it N$= \! n$ superanalytic structure} over a Riemann surface $M_B$, we mean an equivalence class of N$= \! n$ superanalytic equivalent atlases on a DeWitt $(1,n)$-supermanifold $M$ whose body is $M_B$. 

Now, let $M_1$ and $M_2$ be N$= \! n$ superconformal DeWitt super-Riemann surfaces, for $n=1,2$, with coordinate atlases $\{(U_{\alpha},$ $\Omega_{\alpha})\}_{\alpha \in A}$ and $\{(V_{\beta}, \Xi_{\beta})\}_{\beta \in B}$, respectively.  A map $F: M_1 \longrightarrow M_2$ is said to be {\it N$= \! n$ superconformal}  if $\Xi_\beta \circ F \circ \Omega_\alpha^{-1}: \Omega_\alpha (U_\alpha \cap F^{-1} (V_\beta)) \longrightarrow \Xi_\beta(V_\beta)$ is N$= \! n$ superconformal for all $\alpha \in A$ and $\beta \in B$ with $U_\alpha \cap F^{-1} (V_\beta) \neq \emptyset$.  If in addition, $F$ is bijective, then we say that $M_1$ and $M_2$ are {\it N$= \! n$ superconformally equivalent}.  By {\it N$= \! n$ superconformal structure} over a Riemann surface $M_B$, we mean an equivalence class of N$= \! n$ superconformally equivalent atlases on a DeWitt $(1,n)$-supermanifold $M$ whose body is $M_B$.

We define the {\it N$= \! n$ super complex plane} over $\bigwedge_{*>n-1}$, denoted by $S^n\mathbb{C}$, to be $\mathbb{C} \times (\bigwedge_{*>n-1}^0)_S \times (\bigwedge_{*>n-1}^1)^n = \bigwedge_{*>n-1}^0 \times (\bigwedge_{*>n-1}^1)^n$ with the usual topology on $\mathbb{C}$ dictating the DeWitt topology on $S^n\mathbb{C}$.   (In the notation of \cite{Ro}, this space is denoted by $\mathbb{C}_{S[L]}^{1,n}$ if $\bigwedge_* = \bigwedge_L$, and by $\mathbb{C}_S^{1,n}$ if $\bigwedge_* = \bigwedge_\infty$.)
We define the {\it N$= \! n$ super upper half-plane} over $\bigwedge_{*>n-1}$, denoted by $S^n\mathbb{H}$, to be $\mathbb{H} \times (\bigwedge_{*>n-1}^0)_S \times (\bigwedge_{*>n-1}^1)^n$ with the usual topology on $\mathbb{H}$ dictating the DeWitt topology on $S^n\mathbb{H}$. 

Note that for N$= \! n$ with $n = 1,2$, both the superplane $S^n\mathbb{C}$ and the super upper half-plane $S^n\mathbb{H}$ are not only superanalytic as supermanifolds, but also N$= \! n$ superconformal.  In addition, these are examples of $\mathcal{C} _{*>n-1}(n)$-supermanifolds.

\section{The equivalence of N=2 superconformal and N=1 superanalytic DeWitt super-Riemann surfaces}\label{equivalence-section}

In this section, we recall some results from \cite{DRS} establishing an equivalence between N=1 superanalytic  super-Riemann surfaces and N=2 superconformal super-Riemann surfaces.  Our main result in this paper, Theorem \ref{uniformization-prop} in Section \ref{uniformization-section} is formulated and proved for N=2 superconformal DeWitt super-Riemann surfaces.  In Section \ref{uniformization-section}, we will use the results from \cite{DRS} stated in this section to formulate Corollary \ref{uniformization-cor} to Theorem \ref{uniformization-prop} which gives a uniformization theorem for N=1 superanalytic DeWitt super-Riemann surfaces.  

Although we follow \cite{DRS}, there are discrepancies between some of our formulas and those given in \cite{DRS}.  For instance, there is a typo in \cite{DRS} in the transformation from the nonhomogeneous coordinate system $(z, \theta_1, \theta_2)$ to the homogeneous coordinate system $(z, \theta^+, \theta^-)$; this typo is a factor of $1/2$ erroneously introduced into the $D^\pm$ superderivations after the transformation of coordinates, and this factor is carried throughout their calculations. 

Let $U_B$ be an open set in $\mathbb{C}$.  Let  $\mathcal{SC}_{*>1}(2, U_B)$ be the set of invertible N=2 superconformal functions defined on the DeWitt open set $U_B \times ((\bigwedge_{*>1}^0)_S \oplus (\bigwedge_{*>1}^1)^2)$ in $\bigwedge_{*>1}^0 \oplus (\bigwedge_{*>1}^1)^2$.  Let $\mathcal{SA}_{*>1}(1, U_B)$ be the set of invertible N=1 superanalytic functions $H$ defined on the DeWitt open set $U_B \times (\bigwedge_{*>1})_S$ in $\bigwedge_{*>1}$ such that  the coefficients of the functions defining $H$ are restricted to lie in $\bigwedge_{*-2}$ rather than just in $\bigwedge_{*-1}$; that is in (\ref{restrict-coefficients}), we take $(k) \in J_{*-2}$ rather than $(k) \in J_{*-1}$.  

Define the map
\begin{eqnarray}
\mathcal{F}_1 : \mathcal{SC}_{*>1}(2, U_B) & \longrightarrow & \mathcal{SA}_{*>1}(1, U_B)\\
H & \mapsto & \mathcal{F}_1(H) \nonumber
\end{eqnarray}
as follows: For $H \in \mathcal{SC}_{*>1}(2, U_B)$, then in particular $H(z, \theta^+, \theta^-) = (\tilde{z}, \tilde{\theta}^+, \tilde{\theta}^-)$ is of the form (\ref{superconformal-condition1})-(\ref{superconformal-condition4}) for even functions $f$ and $g^\pm$ and odd functions $\psi^\pm$.    Define 
\begin{equation}\label{F1-definition}
\mathcal{F}_1(H) (z, \theta) = (f(z) + \psi^+(z) \psi^-(z) +  2\theta g^+(z) \psi^-(z), \   \psi^+(z) + \theta g^+(z)) .
\end{equation}
It is a straightforward calculation to show that $\mathcal{F}_1$ is in fact a homomorphism of pseudogroups. 

The invertible N=1 superanalytic function $\mathcal{F}_1(H)$ can be thought of as arising from performing the N=2 superanalytic coordinate transformation 
\begin{equation}
(z, \ \theta^+, \ \theta^-) \mapsto (u, \ \eta,  \ \alpha) = (z + \theta^+ \theta^-, \ \theta^+, \ \theta^-).
\end{equation}
Under this transformation, we obtain the N=2 superanalytic function in the even variable $u$ and the two odd variables $\eta$ and $\alpha$ given by 
\begin{eqnarray}
\qquad \qquad \tilde{u} &=& f(u) + \psi^+(u) \psi^-(u) +  2\eta g^+(u) \psi^-(u) \\
\tilde{\eta} &=& \psi^+(u) + \eta g^+(u) \\
\tilde{\alpha} &=& \psi^-(u) + \alpha g^-(u) -2 \eta \alpha (\psi^-)'(u) .
\end{eqnarray}

Conversely, define the map
\begin{eqnarray}
\mathcal{F}_2 : \mathcal{SA}_{*>1}(1, U_B) & \longrightarrow & \mathcal{SC}_{*>1}(2, U_B)\\
H & \mapsto & \mathcal{F}_2(H)  \nonumber
\end{eqnarray}
as follows: For $H \in \mathcal{SA}_{*>1}(1, U_B)$, then $H(z, \theta) = (f_1(z) + \theta \xi(z), \ \psi(z) + \theta g(z))$ for even functions $f_1(z)$ and $g(z)$ and odd functions $\xi(z)$ and $\psi(z)$, and with $g(z)$ nonvanishing.  Define $\mathcal{F}_2(H) (z, \theta^+, \theta^-) =  (\tilde{z}, \tilde{\theta}^+, \tilde{\theta}^-)$ to be of the form (\ref{superconformal-condition1})-(\ref{superconformal-condition3}) where
\begin{eqnarray}
f(z) &=& f_1(z) - \frac{\psi(z)\xi(z)}{2 g(z)}, \\
g^+(z) &=& g(z), \qquad \ \ \mbox{and} \qquad \ g^-(z)\ \  = \ \ \frac{f_1'(z)}{g(z)} - \frac{\psi'(z) \xi(z)}{g(z)^2}\\
\psi^+(z) &=& \psi(z), \qquad\ \  \mbox{and} \qquad \psi^-(z) \ \ = \ \ \frac{\xi(z)}{2g(z)}.
\end{eqnarray}
One can easily check that condition (\ref{superconformal-condition4}) is satisfied, and thus $\mathcal{F}_2(H)$ is indeed N=2 superconformal. 

We have that $\mathcal{F}_1$ and $\mathcal{F}_2$ are bijections and 
\begin{equation}\label{inverse-functors}
\mathcal{F}_1 \circ \mathcal{F}_2 = id_{\mathcal{SA}_{*>1}(1, U_B)} \qquad  \mbox{and} \qquad \mathcal{F}_2 \circ \mathcal{F}_1 = id_{\mathcal{SC}_{*>1}(2, U_B)}.
\end{equation}

Let $\mathcal{SCM}_{*>1}(2)$ be the category of N=2 superconformal DeWitt super-Riemann surfaces over the Grassmann algebra $\bigwedge_{*>1}$, and let $\mathcal{SAM}_{*>1}(1)$ be the category of N=1 superanalytic DeWitt super-Riemann surfaces $M$ over the Grassmann algebra $\bigwedge_{*>1}$ such that the transition functions for $M$ are in $\mathcal{SA}_{*>1}(1, U_B)$ for some $U_B \in \mathbb{C}$. 

Define the functor 
\begin{eqnarray}
\mathcal{F}: \mathcal{SCM}_{*>1}(2) &\longrightarrow& \mathcal{SAM}_{*>1}(1)\\
M & \mapsto & \mathcal{F}(M) \nonumber
\end{eqnarray}
as follows:  Let $M$ be an N=2 superconformal DeWitt super-Riemann surface over the Grassmann algebra $\bigwedge_{*>1}$ with coordinate atlas $\{ (U_\alpha, \Omega_\alpha)\}_{\alpha \in A}$.   Let $\mathcal{F}(M)$ be the N=1 superanalytic DeWitt super-Riemann surface with body $M_B$ obtained by patching together  DeWitt open domains in $\bigwedge_{*>1}$ with local coordinates $(z, \theta)$ by means of the transition functions $\mathcal{F}_1 ( \Omega_\alpha \circ \Omega_\beta^{-1}):  (\Omega_\beta(U_\alpha \cap U_\beta))_B \times (\bigwedge_{*>1})_S \longrightarrow (\Omega_\alpha (U_\alpha \cap U_\beta))_B \times (\bigwedge_{*>1})_S$.

Since $\mathcal{F}_1$ is a homomorphism of pseudogroups, the functor $\mathcal{F}$ is well defined.  
From (\ref{inverse-functors}), it follows that $\mathcal{F}$ is an isomorphism of categories.  Thus we have the following proposition:
\begin{prop}\label{DRS-prop}
The category $\mathcal{SCM}_{*>1}(2)$  of N=2 superconformal DeWitt super-Riemann surfaces over the Grassmann algebra $\bigwedge_{*>1}$ is isomorphic to the category $\mathcal{SAM}_{*>1}(1)$ of N=1 superanalytic DeWitt super-Riemann surfaces such that  the coefficients of the coordinate transition functions are restricted to lie in $\bigwedge_{*-2}$.
\end{prop}

\begin{rema}
{\em  In N=2 superconformal field theory, the supermanifolds that arise from superstrings propagating through space time, are N=2 superconformal DeWitt super-Riemann surfaces with half-infinite tubes attached.  These half-infinite tubes are N=2 superconformally equivalent to punctures on the N=2 superconformal DeWitt super-Riemann surface with N=2 superconformal local coordinates vanishing at the punctures.  Although, there is a bijection between N=2 superconformal super-Riemann surfaces and N=1 superanalytic super-Riemann surfaces, there is no such bijection when the extra data of punctures and local coordinates vanishing at the punctures is added.    Even for a marked point such as the origin on the corresponding superplanes, there is not a bijection between the N=2 superconformal local coordinates vanishing at the origin of $\bigwedge_{*>1}^0 \oplus (\bigwedge_{*>1}^1)^2$ and N=1 superanalytic local coordinates vanishing at the origin of $\bigwedge_{*>1}$.   For example, the N=1 superanalytic functions $H(z, \theta) = ( z + \theta \xi(z), \theta)$ vanish at the origin $(0,0)$ of the N=1 superplane $\bigwedge_{*>1}$.  However the corresponding N=2 superconformal functions $\mathcal{F}_1(H) (z, \theta^+, \theta^-) = (z + \frac{1}{2} \theta^+, \theta^+, \frac{1}{2}\xi(z) + \theta^-)$ vanish at $(0,0, -\frac{1}{2} \xi(0))$.   
   Thus one cannot simply replace N=2 superconformal worldsheets swept out by propagating superstrings by N=1 superanalytic worldsheets when the full data of the propagating strings---punctures and local coordinates vanishing at the punctures---is included.  One must either work in the N=2 superconformal setting, or take into account the discrepancies that arise by using the N=1 superanalytic setting when modeling the incoming and outgoing tubes for the superstrings.  See, for example, \cite{B-axiomatic-deformations} for further discussion of this fact. }
\end{rema}

\section{A general uniformization theorem for N=2 superconformal and N=1 superanalytic DeWitt super-Riemann surfaces}\label{uniformization-section}

We now prove the main theorem of this paper:

\begin{thm}\label{uniformization-prop} 
Every N=2 superconformal DeWitt super-Riemann surface $M$ with body $M_B$ is N=2 superconformally equivalent to a $\mathcal{H}_{*>1}(2)$-supermanifold if and only if the first \v  Cech cohomology group of $M_B$ with coefficients in $\mathcal{L}^{-1} \otimes TM_B$ is trivial, where $\mathcal{L}$ is a holomorphic line bundle over $M_B$ and $TM_B$ is the tangent bundle of $M_B$.  In other words, if $\check{H}^1(M_B, \mathcal{L}^{-1} \otimes TM_B) = 0$, then any N=2 superconformal DeWitt supermanifold $M$ with body $M_B$ is N=2 superconformally equivalent to a supermanifold with transition functions of the form 
\begin{equation}\label{reduced-transition}
H(z, \theta^+, \theta^-) = ( f(z) , \  \theta^+ g^+(z), \ \theta^- g^-(z) ) 
\end{equation}
for $f(z)$ and $g^\pm(z)$ even superanalytic $(1,0)$-superfunctions in $z$, satisfying the condition 
\begin{equation}\label{superconformal-condition-reduced}
f'(z) \; = \; g^+(z) g^-(z) .
\end{equation}

Moreover, if $M$ has transition functions of the form  (\ref{reduced-transition}), satisfying (\ref{superconformal-condition-reduced}), and $\check{H}^1(M_B,  TM_B) = 0$, then $M$ is N=2 superconformally equivalent to a supermanifold with transition functions of the form (\ref{reduced-transition}), satisfying (\ref{superconformal-condition-reduced}), such that $f(z_B)$ takes values in $\mathbb{C}$ instead of more generally in $\bigwedge_*^0$. 
\end{thm}

\begin{proof} We follow the spirit of the proof of the uniformization theorem in the N=1 superconformal genus-zero case given by Crane and Rabin in \cite{CR}.   Let $M$ be an N=2 superconformal DeWitt super-Riemann surface with coordinate atlas $\{(U_{\alpha}, \Omega_{\alpha})\}_{\alpha \in A}$.    For any $\alpha, \beta \in A$ with $U_\alpha \cap U_\beta \neq \emptyset$, we have transition function $H_{\alpha \beta} =  \Omega_\alpha \circ \Omega_\beta^{-1} : \Omega_\beta (U_\alpha \cap U_\beta) \longrightarrow \Omega_\alpha(U_\alpha \cap U_\beta)$.  We will write $H_{\alpha \beta}(z, \theta^+, \theta^-) = (\tilde{z}_{\alpha \beta}, \tilde{\theta}^+_{\alpha \beta}, \tilde{\theta}^-_{\alpha \beta})$ and denote the three even superfunctions in $z$ and two odd superfunctions in $z$ that uniquely determine $H_{\alpha \beta}$ according to (\ref{superconformal-condition1})--(\ref{superconformal-condition4}) by $f_{\alpha \beta}$, $g^\pm_{\alpha \beta}$ and $\psi^\pm_{\alpha \beta}$, respectively. 

On each triple intersection $U_\alpha \cap U_\beta \cap U_\gamma \neq \emptyset$, for $\alpha, \beta, \gamma \in A$, we have the consistency condition $H_{\alpha \gamma} = H_{\alpha \beta} \circ H_{\beta \gamma}$, which when expanded in terms of component functions give the conditions
\begin{eqnarray}
f_{\alpha \gamma}  (z) \! \! \! \! &=&\! \! \! \!  f_{\alpha \beta} (f_{\beta \gamma} (z)) + g^+_{\alpha \beta} (f_{\beta \gamma}(z)) \psi^+_{\beta \gamma} (z)   \psi^-_{\alpha \beta} (f_{\beta \gamma}(z)) -  g^-_{\alpha \beta} ( f_{\beta \gamma}(z))\label{consistency-f}\\
& & \quad \cdot \psi^+_{\alpha \beta} ( f_{\beta \gamma} (z)) \psi^-_{\beta \gamma} (z) +  (\psi^+_{\alpha \beta} \psi^-_{\alpha \beta})' (f_{\beta \gamma} (z)) \psi^+_{\beta \gamma} (z) \psi^-_{\beta \gamma}  (z)\nonumber\\
\qquad \psi^\pm_{\alpha \gamma} (z) \! \! \! \! &=& \! \!  \! \! \psi^\pm_{\alpha \beta} (f_{\beta \gamma} (z))+ g^\pm_{\alpha \beta} ( f_{\beta \gamma}(z))\psi^\pm_{\beta \gamma} (z)  \pm(\psi^\pm_{\alpha \beta})' (f_{\beta \gamma}(z))  \psi^+_{\beta \gamma}(z) \psi^-_{\beta \gamma}(z) \label{consistency-psi}\\
g^\pm_{\alpha \gamma} (z) \! \! \! \!  &=& \! \! \! \! g^\pm_{\alpha \beta}(f_{\beta \gamma}(z))g^\pm_{\beta \gamma}(z)  - 2(\psi^\pm_{\alpha \beta})'(f_{\beta \gamma}(z))g^\pm_{\beta \gamma} (z) \psi^\mp_{\beta \gamma}(z)  \label{consistency-g}\\
& & \quad -  (g^\pm_{\alpha \beta})'(f_{\beta \gamma}(z)) g^\pm_{\beta \gamma}(z) \psi^\pm_{\beta \gamma}(z) \psi^\mp_{\beta \gamma} (z) .\nonumber
\end{eqnarray}

We will use the equations above to first show that, in general, there exist N=2 superconformal changes of coordinates in each coordinate chart that give a new atlas for which the $\psi_{\alpha \beta}^\pm$ terms in the coordinate transition functions are equal to zero for all $\alpha,  \beta \in A$ if and only if $\check{H}^1 (M_B, \mathcal{L}^{-1} \otimes TM_B) = 0$.  Then we will show that similarly, we can in general set the soul part of the $f_{\alpha \beta}$ terms equal to zero if and only if $\check{H}^1 (M_B, TM_B) = 0$.  

We first expand each $\psi^\pm_{\alpha \beta}$ into its component functions, writing
\begin{equation}
\psi^\pm_{\alpha \beta}(z) = \sum_{(j) \in J^1_{*-2} } (\psi^\pm_{\alpha \beta})_{(j)} (z) \zeta_{j_1} \zeta_{j_2} \cdots \zeta_{j_{2n+1}} .
\end{equation}
We will show by induction on $n \in \mathbb{N}$, that by N=2 superconformal change of coordinates in each chart, we can set any nonzero $(\psi^\pm_{\alpha \beta})_{(j)}$ equal to zero for $(j) = (j_1, j_2, \dots, j_{2n+1}) \in J_{*-2}^1$ if and only if $\check{H}^1 (M_B, \mathcal{L}^{-1} \otimes TM_B) = 0$.  Let $z^\alpha$ denote the even coordinate on $\Omega_\alpha(U_\alpha)$ for $\alpha \in A$.

For $n=1$, letting $(j) = (j_1) \in J^1_{*-2}$ for $j_1 \in \{1,\dots, *-2\}$, the equations (\ref{consistency-psi}) reduce to
\begin{equation}\label{pre-cocycle}
( \psi^\pm_{\alpha \gamma} )_{(j)} (z^\gamma_{(\emptyset)}) = ( \psi^\pm_{\alpha \beta})_{(j)} (z^\beta_{(\emptyset)})+ (\psi^\pm_{\beta \gamma})_{(j)} (z^\gamma_{(\emptyset)}) (g^\pm_{\alpha \beta})_{(\emptyset)} ( z^\beta_{(\emptyset)}) .
\end{equation}
By the N=2 superconformal condition (\ref{superconformal-condition4}), we have
\begin{equation}\label{using-superconformal-condition4}
(f'_{\alpha \beta })_{(\emptyset)} (z^\beta_{(\emptyset)}) = (g^+_{\alpha \beta})_{(\emptyset)}  (z^\beta_{(\emptyset)}) (g^-_{\alpha \beta})_{(\emptyset)}  (z^\beta_{(\emptyset)}) .
\end{equation}

Since $(f_{\alpha \beta})_{(\emptyset)}$ is a local homeomorphism from an open set in the complex plane to an open set in the complex plane, we have that $(f_{\alpha \beta})_{(\emptyset)}' (z_{(\emptyset)}^\beta)$ is nonzero for all $z_{(\emptyset)}^\beta \in (U_\beta)_B$, which along with (\ref{using-superconformal-condition4}) implies that the $(g^\pm_{\alpha \beta})_{(\emptyset)} (z_{(\emptyset)}^\beta)$ are nonzero for all $z_{(\emptyset)}^\beta \in (U_\beta)_B$.  Thus we can write
\begin{equation}\label{section-transition}
 (g^\pm_{\alpha \beta})_{(\emptyset)} ( z^\beta_{(\emptyset)}) = \frac{1}{ (g^\mp_{\alpha \beta})_{(\emptyset)} ( z^\beta_{(\emptyset)}) } (f'_{\alpha \beta })_{(\emptyset)} (z^\beta_{(\emptyset)}) .
\end{equation}

Let $\sigma^\pm$ each be a section of a holomorphic bundle $\mathcal{L}_\pm^{-1} \otimes TM_B$ where $\mathcal{L}_\pm$ are the line bundles defined by the transition functions $(g^\pm_{\alpha \beta})_{(\emptyset)}$, respectively, for $\alpha, \beta \in A$, and the holomorphic vector field is defined by the transition functions $(f'_{\alpha \beta})_{(\emptyset)}$.  Then the transition function from $U_\alpha$ to $U_\beta$ for this section is the righthand side of Eq. (\ref{section-transition}); that is
\begin{equation}
\sigma^\pm_\beta = \frac{1}{ (g^\mp_{\alpha \beta})_{(\emptyset)} ( z^\beta_{(\emptyset)}) } (f'_{\alpha \beta })_{(\emptyset)} (z^\beta_{(\emptyset)})  \sigma^\pm_\alpha,
\end{equation}
on $U_\alpha \cap U_\beta$.
Then multiplying equation (\ref{pre-cocycle}) by these holomorphic vector fields, respectively, we have
\begin{equation}\label{cocycle}
\sigma^\pm_\alpha ( \psi^\pm_{\alpha \gamma} )_{(j)} (z^\gamma_{(\emptyset)}) = \sigma^\pm_\alpha  ( \psi^\pm_{\alpha \beta})_{(j)} (z^\beta_{(\emptyset)})+ \sigma^\pm_\beta (\psi^\pm_{\beta \gamma})_{(j)} (z^\gamma_{(\emptyset)}).
\end{equation}
This implies that $\{\sigma^\pm_\alpha  ( \psi^\pm_{\alpha \beta})_{(j)} \; | \; (\alpha, \beta) \in A \times A \}$ are the local representatives of a cocycle in the first \v Cech cohomology group of $M_B$ with coefficients in the sheaf $\mathcal{L}^{-1} \otimes TM_B$.  Thus these cocycles are, in general, coboundaries if and only if this cohomology is trivial.  If this cohomology is trivial, there exist elements $b^\pm_{(j)}$ in the zeroth cohomology $\check{H}^0(M_B, \mathcal{L}^{-1} \otimes TM_B)$, i.e. global sections, such that 
\begin{equation}\label{boundary1}
\sigma^\pm_\alpha  ( \psi^\pm_{\alpha \beta})_{(j)}  =  \sigma^\pm_\beta (b^\pm_{(j)})_\beta - \sigma^\pm _\alpha (b^\pm_{(j)})_\alpha,
\end{equation}
or in other words
\begin{equation}\label{boundary2}
( \psi^\pm_{\alpha \beta})_{(j)}  =  (g^\pm_{\alpha \beta})_{(\emptyset)} (b^\pm_{(j)})_\beta - (b^\pm_{(j)})_\alpha .
\end{equation}

For all $\alpha \in A$, define the N=2 superconformal transformation $H^\alpha_{(j)}$ by $f_\alpha(z^\alpha) = z^\alpha$, $g_\alpha^\pm(z^\alpha) = 1$ and $\psi_\alpha^\pm (z^\alpha) =  (b_{(j)}^\pm)_\alpha(z^\alpha) \zeta_{j_1}$, where $\zeta_{j_1}$ is the $j_1$-th basis element for the underlying Grassmann algebra. That is  $H^\alpha_{(j)} (z^\alpha, ( \theta^+)^\alpha, (\theta^-)^\alpha) = (\tilde{z}^\alpha, (\tilde{\theta}^+)^\alpha, (\tilde{\theta}^-)^\alpha)$ is given by 
\begin{eqnarray}
\tilde{z}^\alpha &=&  z^\alpha + (\theta^+)^\alpha (b^-_{(j)})_\alpha(z^\alpha)\zeta_{j_1} + (\theta^-)^\alpha  (b^+_{(j)})_\alpha(z^\alpha)\zeta_{j_1}\\
(\tilde{\theta}^\pm)^\alpha &=&  (b_{(j)}^\pm)_\alpha (z^\alpha) \zeta_{j_1} + (\theta^\pm)^\alpha \pm (\theta^+)^\alpha (\theta^-)^\alpha ((b_{(j)}^\pm)_\alpha)'(z^\alpha) \zeta_{j_1} .
\end{eqnarray}
Redefining the coordinates for the charts $(U_\alpha, \Omega_\alpha)$ and $(U_\beta, \Omega_\beta)$ by the N=2 superconformal transformations $H^\alpha_{(j)}$ and $H^\beta_{(j)}$, respectively, such that the new coordinate charts are $(U_\alpha, \tilde{\Omega}_\alpha = H^\alpha_{(j)} \circ \Omega_\alpha)$ and $(U_\beta, \tilde{\Omega}_\beta = H^\beta_{(j)} \circ \Omega_\beta)$, we see that the new coordinate transformation $\tilde{H}_{\alpha \beta} = \tilde{\Omega}_\alpha \circ \tilde{\Omega}_\beta^{-1} = H^\alpha_{(j)} \circ \Omega_\alpha \circ \Omega_\beta^{-1} \circ (H^\beta_{(j)})^{-1}$, has
\begin{eqnarray}
( \tilde{\psi}^\pm_{\alpha \beta})_{(j)} &=& ( \psi^\pm_{\alpha \beta})_{(j)} -   (g^\pm_{\alpha \beta})_{(\emptyset)} (b^\pm_{(j)})_\beta + (b^\pm_{(j)})_\alpha \\
&=& 0. \nonumber
\end{eqnarray}

The new atlas $\{(U_\alpha, \tilde{\Omega}_\alpha) \: | \; \alpha \in A \}$ now has the $\tilde{\psi}_{(j)}^\pm$ terms for $(j) = (j_1) \in J^1_{*-2}$ equal to zero.  For this new atlas, we again have the compatibility condition on the triple overlaps, and thus the consistency conditions (\ref{consistency-f})--(\ref{consistency-g}) again hold for the components of the new coordinate transformation functions, and we can perform the above procedure again for a new $j_1 \in \{1,\dots, *-2\}$ without changing the fact that the previous $(\psi^\pm_{\alpha \beta})_{(j)}$ terms have been set equal to zero.  Doing this repeatedly, we can, in general, set the $(\psi^\pm_{\alpha \beta})_{(j)}$ terms equal to zero for $(j) = (j_1)$, for all $j_1 = 1,\dots, *-2$ if and only if $\check{H}^1 (M_B, \mathcal{L}^{-1} \otimes TM_B) = 0$.  

Now we make the inductive assumption: assume coordinate transformations have been made such that for the new $\psi^\pm$ terms, the components $(\psi^\pm_{\alpha \beta})_{(j)}$ have been set equal to zero for $(j) = (j_1, j_2, \dots, j_{2k +1}) \in J^1_{*-2}$, and for $k = 1,\dots, n-1$.  For this new atlas, we again have the compatibility condition on the triple overlaps.  Thus for the $\zeta_{j_1} \cdots \zeta_{j_{2n + 1}}$ level of equation (\ref{consistency-psi}) applied to $z_{(\emptyset)}$, we have that equation (\ref{pre-cocycle}) applies to these new $\psi^\pm_{(j)}$'s, and thus, so does equation (\ref{cocycle}).  Thus again, in general, there exist $b^\pm_{(j)}$, each in $\check{H}^0 (M_B, \mathcal{L}^{-1} \otimes TM_B)$ such that equations (\ref{boundary1}) and (\ref{boundary2}) hold for these terms if and only if $\check{H}^1(M_B, \mathcal{L}^{-1} \otimes TM_B)$ is trivial.   

For $\alpha \in A$, define the N=2 superconformal transformation $H^\alpha_{(j)}$ by $f_\alpha(z^\alpha) = z^\alpha$, $g_\alpha^\pm(z^\alpha) = 1$ and $\psi_\alpha^\pm (z^\alpha) =  (b_{(j)}^\pm)_\alpha(z^\alpha) \zeta_{j_1} \cdots \zeta_{j_{2n+1}}$. Redefining the coordinates for the charts $(U_\alpha, \Omega_\alpha)$ by the N=2 superconformal transformations $H^\alpha_{(j)}$, for $\alpha \in A$, it is a straightforward calculation to see that for the new coordinate coordinate charts $(U_\alpha, \tilde{\Omega}_\alpha = H^\alpha_{(j)} \circ \Omega_\alpha)$ for $\alpha \in A$ the new coordinate transformations $\tilde{H}_{\alpha \beta} = \tilde{\Omega}_\alpha \circ \tilde{\Omega}_\beta^{-1} = H^\alpha_{(j)} \circ \Omega_\alpha \circ \Omega_\beta^{-1} \circ (H^\beta_{(j)})^{-1}$, have $( \tilde{\psi}^\pm_{\alpha \beta})_{(j)} = 0$.  We perform this procedure of redefining the coordinate charts for each $(j) = (j_1, \dots, j_{2k+1}) \in J^1_{*-2}$ with $k=n$, resulting in at atlas with coordinate transition functions that have all the $\psi^\pm_{(j)}$ components equal to zero for $(j)$ of length $2n+1$, while keeping the transition functions $\psi^\pm_{(j')}$ for $(j') \in J^1_{*-2}$ of length $2k+1$ for $k<n$ equal to zero.  

This proves that, in general, there exist N=2 superconformal coordinate transformations which result in an atlas of charts with coordinate transitions functions for which the $\psi^\pm$ terms are all zero if and only if $\check{H}^1(M_B, \mathcal{L}^{-1} \otimes TM_B)$ is trivial; that is, every N=2 superconformal DeWitt super-Riemann surface $M$ with body $M_B$ is N=2 superconformally equivalent to a $\mathcal{H}_{*>1}(2)$-supermanifold if and only if $\check{H}^1(M_B, \mathcal{L}^{-1} \otimes TM_B)=0$.   

The remaining nontrivial consistency conditions (\ref{consistency-f}) and (\ref{consistency-g}) for the new coordinate atlas reduce to 
\begin{eqnarray}
f_{\alpha \gamma}  (z) &=&  f_{\alpha \beta} (f_{\beta \gamma} (z))  \label{consistency-f2}\\
g^\pm_{\alpha \gamma} (z) &=&  g^\pm_{\alpha \beta}(f_{\beta \gamma}(z)) g^\pm_{\beta \gamma}(z) . \label{consistency-g2}
\end{eqnarray}

We now expand the $f_{\alpha \beta}$ terms into component functions, writing 
\begin{equation}
f_{\alpha \beta} (z) = \sum_{(j) \in J^0_{*-2} }(f_{\alpha \beta})_{(j)} (z) \zeta_{j_1} \zeta_{j_2} \cdots \zeta_{j_{2n}}. 
\end{equation}
For the $n=0$ terms, i.e., the $(\emptyset) \in J^0_{*-2}$ terms, equation (\ref{consistency-f2}) becomes
\begin{equation}\label{f-zero}
(f_{\alpha \gamma})_{(\emptyset)}  (z^\gamma_{(\emptyset)}) \; = \;   (f_{\alpha \beta})_{(\emptyset)} ((f_{\beta \gamma})_{(\emptyset)} (z^\gamma_{(\emptyset)})) \;   =  \; (f_{\alpha \beta})_{(\emptyset)} (z^\beta_{(\emptyset)}),
\end{equation}
which is just the usual cocycle condition on the Riemann surface $M_B$.   

We will show by induction on $n \in \mathbb{N}$ that we can, in general, set the soul components of these $f$ terms equal to zero if and only if $\check{H}^1(M_B, TM_B) = 0$, assuming the $\psi^\pm$ have already been set to zero.  For $n=1$, let $(j) = (j_1, j_2) \in J^0_{*-2}$.  The consistency equation (\ref{consistency-f2})  implies
\begin{equation}\label{f-two}
(f_{\alpha \gamma})_{(j)}  (z^\gamma_{(\emptyset)}) = (f_{\alpha \beta})_{(j)} (z_{(\emptyset)}^\beta) +  (f_{\alpha \beta})_{(\emptyset)}'(z^\beta_{(\emptyset)})  (f_{\beta \gamma})_{(j)} (z^\gamma_{(\emptyset)}) .
\end{equation}

Let $\sigma$ be a section of the holomorphic tangent bundle such that on coordinate charts we have
\begin{equation}\label{pre-cocycle-f}
\sigma_\beta = (f_{\alpha \beta })_{(\emptyset)}' (z^\beta_{(\emptyset)})  \sigma_\alpha.
\end{equation}
Multiplying equation (\ref{f-two}) by this vector field, we have
\begin{equation}\label{cocycle-f}
\sigma_\alpha (f_{\alpha \gamma})_{(i)}  (z^\gamma_{(\emptyset)}) =  \sigma_\alpha (f_{\alpha \beta})_{(i)} (z_{(\emptyset)}^\beta) + \sigma_\beta (f_{\beta \gamma})_{(i)} (z^\gamma_{(\emptyset)}) .
\end{equation}
This implies that $\{\sigma_\alpha  ( f_{\alpha \beta})_{(j)} \; | \; (\alpha, \beta) \in A \times A \}$ are the local representatives of a cocycle in $\check{H}^1(M_B, TM_B)$.  In general, these cocycles are coboundaries if and only if this cohomology group is trivial; in which case, there exists elements $b_{(j)}$ in the zeroth cohomology group of $M_B$, such that 
\begin{equation}\label{boundary1-f}
\sigma_\alpha  (f_{\alpha \beta})_{(j)}  =  \sigma_\beta (b_{(j)})_\beta - \sigma _\alpha (b_{(j)})_\alpha,
\end{equation}
i.e.,
\begin{equation}\label{boundary2-f}
( f_{\alpha \beta})_{(j)}  =  (f_{\alpha \beta})_{(\emptyset)}' (b_{(j)})_\beta - (b_{(j)})_\alpha .
\end{equation}

For $\alpha \in A$, define the N=2 superconformal transformation $H^\alpha_{(j)}$ by $f_\alpha(z^\alpha) = z^\alpha + (b_{(j)})_\alpha(z^\alpha) \zeta_{j_1} \zeta_{j_2}$,  $g_\alpha^+(z^\alpha) = (f_\alpha)'(z^\alpha)$,  $g_\alpha^-(z^\alpha) =1$ and $\psi_\alpha^\pm (z^\alpha) = 0$. That is  $H^\alpha_{(j)} (z^\alpha, ( \theta^+)^\alpha, (\theta^-)^\alpha) = (z^\alpha + (b_{(j)})_\alpha(z^\alpha) \zeta_{j_1} \zeta_{j_2}, \theta^+ (1 +  (b_{(j)})_\alpha'(z^\alpha) \zeta_{j_1} \zeta_{j_2}), \theta^-)$. Redefining the coordinates for the charts $(U_\alpha, \Omega_\alpha)$ by the N=2 superconformal transformations $H^\alpha_{(j)}$, for $\alpha \in A$, such that the new coordinate charts are $(U_\alpha, \tilde{\Omega}_\alpha = H^\alpha_{(j)} \circ \Omega_\alpha)$, we have that the new coordinate transformations $\tilde{H}_{\alpha \beta} = \tilde{\Omega}_\alpha \circ \tilde{\Omega}_\beta^{-1} = H^\alpha_{(j)} \circ \Omega_\alpha \circ \Omega_\beta^{-1} \circ (H^\beta_{(j)})^{-1}$, for $\alpha, \beta \in A$, have
\begin{eqnarray}
( \tilde{f}_{\alpha \beta})_{(j)} &=& ( f_{\alpha \beta})_{(j)} -   (f_{\alpha \beta})_{(\emptyset)}' (b_{(j)})_\beta + (b_{(j)})_\alpha \; = \; 0\\
\tilde{\psi}^\pm_{\alpha \beta} &=&  0.
\end{eqnarray}

The new atlas $\{(U_\alpha, \tilde{\Omega}_\alpha) \: | \; \alpha \in A \}$ now has the $\tilde{f}_{(j)}$ terms for $(j) = (j_1, j_2) \in J^0_{*-2}$ equal to zero and the $\tilde{\psi}^\pm$ terms equal to zero.  For this new atlas, we again have the compatibility condition on the triple overlaps, and thus the consistency condition (\ref{consistency-f2}) again holds on the components of the new coordinate transformation functions.  Therefore, we can perform the above procedure again for a different $(j) = (j_1, j_2) \in J^0_{*-2}$ without changing the fact that the previous $(f_{\alpha \beta})_{(j)}$ terms have been set equal to zero.  Doing this repeatedly, we can, in general, set the $(f_{\alpha \beta})_{(j)}$ terms equal to zero for all $(j) = (j_1, j_2) \in J^0_{*-2}$ if and only if $\check{H}^1(M_B, TM_B)=0$.  

Now we make the inductive assumption: assume coordinate transformations have been made such that for the new $f$ terms, the components $(f_{\alpha \beta})_{(j)}$ have been set equal to zero for $(j) = (j_1, j_2, \dots, j_{2k}) \in J^0_{*-2}$, for $k = 1,\dots, n-1$.  For this new atlas, we again have the compatibility condition on the triple overlaps.  Thus for the $\zeta_{j_1} \cdots \zeta_{j_{2n}}$ level of equation (\ref{consistency-f2}) applied to $z_{(\emptyset)}$, we have that equation (\ref{pre-cocycle-f}) applies to these new $f_{(j)}$'s, and thus, so does equation (\ref{cocycle-f}).  Thus, in general, there exists $b_{(j)}$ in the zeroth cohomology of $M_B$ such that equations (\ref{boundary1-f}) and (\ref{boundary2-f}) hold for these terms if and only if $\check{H}^1(M_B, TM_B)=0$.

For $\alpha \in A$, define the N=2 superconformal transformation $H^\alpha_{(j)}$ by $f_\alpha(z^\alpha) = z^\alpha + (b_{(j)})_\alpha(z^\alpha) \zeta_{j_1} \cdots \zeta_{j_{2n}}$, $g_\alpha^+(z^\alpha) = f_\alpha'(z^\alpha)$, $g^-_\alpha(z^\alpha) = 1$, and $\psi_\alpha^\pm (z^\alpha) = 0$. Redefining the coordinates for the charts $(U_\alpha, \Omega_\alpha)$ by the N=2 superconformal transformations $H^\alpha_{(j)}$, for $\alpha \in A$, it is a straightforward calculation to see that the new coordinate transformations have $( \tilde{f}_{\alpha \beta})_{(j)} = 0$, for $\alpha, \beta \in A$.  In general, we can  perform this procedure of redefining the coordinate charts for each $(j) = (j_1, \dots, j_{2n}) \in J^0_{*-2}$, resulting in an atlas with coordinate transition functions that have all the $f_{(j)}$ components equal to zero for $(j)$ of length $2n$, in addition to the $f_{(j')}$ terms for $(j') \in J^0_{*-2}$ of length $2k$ for $0<k<n$ equal to zero, while keeping the $\psi^\pm$ terms equal to zero, if and only if $\check{H}^1(M_B, TM_B)=0$.  

This proves that, in general, there exist N=2 superconformal coordinate transformations which result in an atlas of charts for $M$ with coordinate transitions functions for which the $\psi^\pm$ terms are all zero if and only if $\check{H}^1(M_B, \mathcal{L}^{-1} \otimes TM_B)=0$, and furthermore, in this case since this implies $\check{H}^1(M_B, TM_B)=0$,
the soul portion of the $f$ terms in the coordinate transition functions can also be set to zero.   
\end{proof}

Using Proposition \ref{DRS-prop}, we immediately obtain from Theorem \ref{uniformization-prop} a uniformization result for N=1 superanalytic DeWitt super-Riemann surfaces over $\bigwedge_{*>1}$ with body $M_B$ and coordinate transition functions with coefficients restricted to $\bigwedge_{*-2}$.   Namely we have that the uniformization of such a  $\mathcal{G}_{*>1}(1)$-supermanifold to an $\mathcal{H}_{*>1}(1)$-supermanifold is dependent on $\check{H}^1(M_B, \mathcal{L}^{-1} \otimes TM_B)$.  However, analyzing the compatibility conditions for coordinate transformations on triple overlaps for any N=1 superanalytic DeWitt super-Riemann surfaces over $\bigwedge_{*>1}$ with body $M_B$, that is where the coordinate transition functions have coefficients restricted to $\bigwedge_{*-1}$, rather than restricted to $\bigwedge_{*-2}$, we see that for $\bigwedge_* = \bigwedge_L$, the conditions on the $k$ levels for $k<L-1$ are exactly the same as for the $k$ levels with the restricted coefficients, and that the $L-1$ level gives the same cocycle property as that for the lower levels.  Thus we can extend our uniformization result to general N=1 superanalytic DeWitt super-Riemann surfaces over $\bigwedge_{*>0}$.  That is we have:

\begin{cor}\label{uniformization-cor}
Any N=1 superanalytic DeWitt super-Riemann surface over $\bigwedge_{*>0}$ with body $M_B$  is N=1 superanalytically equivalent to an $\mathcal{H}_{*>0}(1)$-supermanifold if and only if the first \v  Cech cohomology group of $M_B$ with coefficients in $\mathcal{L}^{-1} \otimes TM_B$ is trivial. In other words, if $\check{H}^1(M_B, \mathcal{L}^{-1} \otimes TM_B) = 0$, then any N=1 superanalytic DeWitt supermanifold $M$ with body $M_B$ is N=1 superanalytically equivalent to a supermanifold with transition functions of the form 
\begin{equation}\label{reduced-superanalytic}
H(z, \theta) = ( f(z), \  \theta g(z))
\end{equation}
where $f(z)$ and $g(z)$ are even superanalytic $(1,0)$-superfunctions in $z$.

Moreover, if $M$ has transition functions of the form (\ref{reduced-superanalytic}) and $\check{H}^1(M_B, TM_B) = 0$, then $M$ is N=1 superanalytically equivalent to a supermanifold with transition functions of the form (\ref{reduced-superanalytic}) such that $f(z_B)$ takes values in $\mathbb{C}$ instead of more generally in $\bigwedge_*^0$. 
\end{cor}

\begin{rema}\label{setting-remark}
{\em Although from Proposition \ref{DRS-prop}, we know that the N=2 superconformal and N=1 superanalytic settings are essentially equivalent, it is easier to prove our uniformization theorem, Theorem \ref{uniformization-prop}, first for N=2 superconformal DeWitt super-Riemann surfaces and then translate our result to N=1 superanalytic DeWitt super-Riemann surfaces, as we did in  Corollary \ref{uniformization-cor}.  The reason for this is that the coordinate transformation compatibility conditions on triple overlaps in the homogeneous coordinate setting for N=2 superconformal DeWitt supermanifolds clearly gives a cocycle in $\check{H}^1(M_B, \mathcal{L}^{-1} \otimes TM_B)$ or $\check{H}^1(M_B, TM_B)$, whereas for N=1 superanalytic DeWitt supermanifolds this clear dependency between  cocycles in these cohomology groups and the coordinate transformation compatibility conditions on triple overlaps is lost.  In particular, letting $H_{\alpha \beta} (z, \theta) = (f_{\alpha \beta} (z) + \theta \xi_{\alpha \beta} (z), \; \psi_{\alpha \beta}(z) + \theta g_{\alpha \beta}(z))$ be the transition function from the $(U_\beta, \Omega_\beta)$ coordinate chart to the $(U_\alpha, \Omega_\alpha)$ coordinate chart for an N=1 superanalytic super-Riemann surface, the compatibility condition $H_{\alpha \gamma} = H_{\alpha \beta} \circ H_{\beta \gamma}$ on triple overlaps implies the following conditions on the odd component functions:
\begin{eqnarray}
\xi_{\alpha \gamma} (z) &=& g_{\beta \gamma}(z) \xi_{\alpha \beta} (f_{\beta \gamma}(z)) + f'_{\alpha \beta} (f_{\beta \gamma}(z)) \xi_{\beta \gamma} (z) \label{first-condition} \\
& & \quad +\, \psi_{\beta \gamma}(z) \xi_{\alpha \beta}' (f_{\beta \gamma}(z)) \xi_{\beta \gamma}(z) \nonumber \\
\psi_{\alpha \gamma}(z) &=& \psi_{\alpha \beta}(f_{\beta \gamma}(z)) + g_{\alpha \beta}(f_{\beta \gamma}(z)) \psi_{\beta \gamma}(z).
\end{eqnarray}
Using induction as in the proof of Theorem 4.1, we look at the coefficients of $\zeta_{j_1} \zeta_{j_2} \cdots \zeta_{j_n}$ first for $n =1$.  Looking at the coefficient of $\zeta_j$, one can see the dependency of the $\psi_{(j)}$ terms on the cohomology $\check{H}^1(M_B, \mathcal{L})$ for $\mathcal{L}$ a line bundle over $M_B$, and can set these terms equal to zero via coordinate redefinitions if this cohomology is trivial.  However, then one obtains that the coefficient of $\zeta_j$ in condition (\ref{first-condition}) is given by 
\begin{equation}
(\xi_{\alpha \gamma})_{(j)} (z_{(\emptyset)}^\gamma) = (g_{\beta \gamma})_{(\emptyset)}(z_{(\emptyset)}^\gamma) (\xi_{\alpha \beta})_{(j)} (z^\beta_{(\emptyset)}) + (f'_{\alpha \beta})_{(\emptyset)} (z^\beta_{(\emptyset)}) (\xi_{\beta \gamma})_{(j)} (z^\gamma_{(\emptyset)}) 
\end{equation}
and this is not, in general, a cocycle in the cohomology over $M_B$ with coefficients in some sheaf over $M_B$.   This more complicated dependency of terms in the coordinate transformation functions in the N=1 superanalytic case in comparison to the N=2 superconformal case is due to the composite nature of the $\xi$ terms in the correspondence between N=2 superconformal supermanifolds and N=1 superanalytic manifolds; in particular, see Eq. (3.2) where $\xi(z) = 2g^+(z) \psi^-(z)$ when transforming an N=2 superconformal supermanifold with component functions $g^\pm$ and $\psi^\pm$ to an N=1 superanalytic supermanifold.  A similar phenomenon occurs if one works in the nonhomogeneous N=2 superconformal coordinate system as discussed in Section 7.   
}
\end{rema}

\section{Uniformization for simply connected N=2 superconformal and N=1 superanalytic DeWitt super-Riemann surfaces}

Although  Corollary \ref{uniformization-cor}, along with the results of \cite{Batchelor} (see also \cite{Ro}) and \cite{Manin1} essentially allow one to conclude the classification of simply connected N=1 superanalytic (and thus also N=2 superconformal) DeWitt super-Riemann surfaces, for completeness, we give a direct proof of this classification here. 

\subsection{A family of inequivalent N=2 superconformal structures over the Riemann sphere}\label{sphere-section}

Let $\mathcal{G}$ be the set of functions $g : (\bigwedge_{*>1}^0)^\times \rightarrow (\bigwedge_{*>1}^0)^\times$, $z \mapsto g(z)$, such that $g$ is superanalytic for $z \in  (\bigwedge_{*>1}^0)^\times$.   That is, $g$ is an even superanalytic function in $z$ such that $g_{(j)}(z_{(\emptyset)})$ for $(j) \in J^0_{*-2}$ is complex analytic for all $z_{(\emptyset)} \in \mathbb{C}^\times$ and $g^+_{(\emptyset)}$ is nonvanishing on $\mathbb{C}^\times$.  Note that $\mathcal{G}$ is a group under point-wise multiplication. 

For $g \in \mathcal{G}$, define the N=2 superconformal map
\begin{eqnarray}
I_g:   \mbox{$(\bigwedge_{*>1}^0)^\times$} \oplus (\mbox{$\bigwedge_{*>1}^1$})^2  & \longrightarrow &   \mbox{$(\bigwedge_{*>1}^0)^\times$} \oplus (\mbox{$\bigwedge_{*>1}^1$})^2 \\
(z, \theta^+, \theta^-) & \mapsto & I_g(z,\theta^+, \theta^-) =  \Bigl(\frac{1}{z}, \; \frac{i\theta^+g(z)}{z}, \; \frac{i\theta^-}{z g(z)} \Bigr) . \nonumber
\end{eqnarray}

Define $S^2\hat{\mathbb{C}}(g)$, for $g \in \mathcal{G}$, to be the genus-zero N=2 superconformal super-Riemann surface over $\bigwedge_{*>1}$ with N=2 superconformal structure given by the covering of local coordinate neighborhoods $\{ U_{\sou_g}, U_{\nor_g} \}$ and the local coordinate maps
\begin{eqnarray}
\sou_g  : U_{\sou_g}  & \longrightarrow & \mbox{$\bigwedge_{*>1}^0$} \oplus (\mbox{$\bigwedge_{*>1}^1$})^2 \label{southern-chart} \\
\nor_g : U_{\nor_g} & \longrightarrow & \mbox{$\bigwedge_{*>1}^0$}  \oplus (\mbox{$\bigwedge_{*>1}^1$})^2, \label{northern-chart}
\end{eqnarray}
which are homeomorphisms of $U_{\sou_g}$ and $U_{\nor_g}$ onto $\bigwedge_{*>1}^0 \oplus (\bigwedge_{*>1}^1)^2$, respectively, such that 
\begin{eqnarray}\label{transition}
\sou_g \circ \nor_g^{-1} : \mbox{$(\bigwedge_{*>1}^0)^\times$} \oplus \mbox{$(\bigwedge_{*>1}^1)^2$} &\longrightarrow& \mbox{$(\bigwedge_{*>1}^0)^\times$} \oplus \mbox{$(\bigwedge_{*>1}^1)^2$} \label{origin-of-I}\\ 
(z, \theta^+,\theta^-) & \mapsto & I_g(z,\theta^+,\theta^-) . \nonumber
\end{eqnarray}
Thus the body of $S^2\hat{\mathbb{C}}(g)$ is the Riemann sphere, i.e., $(S^2\hat{\mathbb{C}}(g))_B = \hat{\mathbb{C}} = \mathbb{C} \cup \{\infty\}$.

The group of N=2 superconformal automorphisms from the N=2 superconformal plane $S^2\mathbb{C}$ to itself that preserve the even coordinate is comprised of transformations 
\begin{eqnarray}\label{group}
T: S^2\mathbb{C} &\longrightarrow& S^2\mathbb{C} \\
(z, \theta^+, \theta^-) & \mapsto & \left(z, \;  \theta^+ \varepsilon^+(z), \; \frac{\theta^-}{\varepsilon^+(z)}  \right) \nonumber
\end{eqnarray}
where $\varepsilon^+(z)$ is an even superanalytic function defined for all $z \in \bigwedge_{*>1}^0$ such that $\varepsilon^+_{(\emptyset)}(z_{(\emptyset)})$ is nonzero for all $z_{(\emptyset)} \in \mathbb{C}$.  The set of all such N=2 superconformal automorphisms of the N=2 superplane that have even component $z$ is a proper subgroup of the the group of N=2 superconformal automorphisms of the N=2 superplane, and in fact, is an abelian subgroup.  

Let $\mathcal{E}^0$ be the set of even superanalytic functions $\varepsilon^+: \bigwedge_{*>1}^0 \rightarrow (\bigwedge_{*>1}^0)^\times$.  That is $\varepsilon^+_{(\emptyset)} (z_{(\emptyset)}) \neq 0$ for all $z_{(\emptyset)} \in \mathbb{C}$.  Then this set $\mathcal{E}^0$ is a group under point-wise multiplication of functions and is isomorphic to the  group of transformations of the form (\ref{group}).  If we restrict the domain of elements of $\mathcal{E}^0$ to $(\bigwedge^0_{*>1})^\times$, we see that $\mathcal{E}^0$ is a subgroup of $\mathcal{G}$.  And since $\mathcal{G}$ is in fact abelian, it is a normal subgroup.   Let $\mathcal{E}^\infty$ be the subgroup of $\mathcal{G}$ given by the set of functions $\{\varepsilon^+(1/z) :  (\bigwedge^0_{*>1})^\times \rightarrow (\bigwedge^0_{*>1})^\times  \; | \; \varepsilon^+(z)  \in \mathcal{E}^0 \}$.  Let $\mathcal{E} = \mathcal{E}^0 \mathcal{E}^\infty$.  Then $\mathcal{E}$ is a proper normal subgroup of  $\mathcal{G}$. 

\begin{lem}\label{equivalent-lemma}
Let $g \in \mathcal{G}$.  There exists $n \in \mathbb{Z}$ such that $S^2 \hat{\mathbb{C}} (g)$ is N=2 superconformally equivalent to $S^2 \hat{\mathbb{C}} (z^n)$.
\end{lem}

\begin{proof}
Let $I^0(z) = 1/z$.  The inverse of $\varepsilon^0$ in $\mathcal{E}^0$, given by $1/\varepsilon^0 (z) = I^0 \circ \varepsilon^0(z)$, is in $\mathcal{E}^0$ and the function $\varepsilon^\infty (1/z) = \varepsilon^\infty \circ I^0(z)$ is also in $\mathcal{E}^0$.  Thus changing coordinates in the chart $(U_{\nor_g}, \nor_g)$ of $S^2 \hat{\mathbb{C}} (g)$ by $T_{I^0 \circ \varepsilon^0}$ and changing coordinates in the chart $(U_{\sou_g}, \sou_g)$ of $S^2 \hat{\mathbb{C}} (g)$ by $T_{\varepsilon^\infty \circ I^0}$, the new change of coordinates from $(U_{\nor_g}, T_{I^0 \circ \varepsilon^0} \circ \nor_g)$ to $(U_{\sou_g}, T_{\varepsilon^\infty \circ I^0} \circ \sou_g)$ is given by 
\begin{eqnarray}
T_{\varepsilon^\infty \circ I^0} \circ I_g \circ T_{I^0 \circ \varepsilon^0}^{-1} (z, \theta^+ ,\theta^-) &=& \left( \frac{1}{z}, \frac{i\theta^+ g(z) \varepsilon(z)}{z}, \frac{i \theta^- }{ z g(z) \varepsilon (z)} \right) \\
&=& I_{h} (z, \theta^+, \theta^-).
\end{eqnarray}
Therefore for $g,h \in \mathcal{G}$, if $h(z) = g(z) \epsilon (z)$ for some $\epsilon \in \mathcal{E}$, then $S^2\hat{\mathbb{C}}(g)$ is N=2 superconformally equivalent to $S^2\hat{\mathbb{C}}(h)$.  That is, there is a surjection from the set $\mathcal{G}/\mathcal{E}$ to the set of N=2 superconformal equivalence classes of $S^2\hat{\mathbb{C}} (g)$, for $g \in \mathcal{G}$. 

Let $g \in \mathcal{G}$ and let $g(z) = g_B(z) + g_S(z)$ be the decomposition of $g$ into body and soul components.  Then defining $f(z) = \log (1 + g_S(z)/g_B(z))$ for $z \in (\bigwedge_{*>1}^0)^\times$ we have that $f$ is a well-defined function from $(\bigwedge_{*>1}^0)^\times$ to $ \bigwedge_{*>1}^0$.  Thus $1 + g_S(z)/g_B(z) = e^{f(z)} \in \mathcal{E}$, and since $g(z) = g_B(z) (1 + g_S(z)/g_B(z))$, we have that $S^2\hat{\mathbb{C}}(g)$ is N=2 superconformally equivalent to $S^2 \hat{\mathbb{C}}(g_B(z))$.  

But now the function $(1/z_B, i\theta^+_1 g_B(z_B)/z_B)$ gives the structure of a holomorphic line bundle over the Riemann sphere $\hat{\mathbb{C}}$.   These holomorphic line bundles are classified by the first Chern class or equivalently by a positive integer degree of the tautological line bundle or its dual.  Thus there is a holomorphic function $\epsilon_B : \mathbb{C}^\times \rightarrow \mathbb{C}^\times$ such that $g_B(z_B) = \epsilon_B(z_B) z_B^n$ for some $n \in \mathbb{Z}$, and extending the domain of $\epsilon_B$ to $(\bigwedge_{*>1}^0)^\times$, we have $\epsilon_B \in \mathcal{E}$.  Thus $S^2 \hat{\mathbb{C}} (g_B)$ is N=2 superconformally equivalent to $S^2 \hat{\mathbb{C}} (z^n)$ for some $n \in \mathbb{Z}$.
\end{proof}

Let 
\begin{equation}
SL(2, \mbox{$\bigwedge_{*-2}^0$}) = \left\{ \left( \begin{array}{cc} a & b \\
c & d 
\end{array}
\right) \; | \; a,b,c,d \in \mbox{$\bigwedge_{*-2}^0$}, \; ad-bc = 1 \right\},
\end{equation}
and let $GL(1, \bigwedge_{*-2}^0) = (\bigwedge_{*-2}^0)^\times$.  

For each $n \in \mathbb{Z}$, and 
\begin{equation}
\alpha = \left( \begin{array}{ccc} a & b & 0 \\
c & d & 0\\
0 & 0& \epsilon
\end{array}
\right) \in SL(2, \mbox{$\bigwedge_{*-2}^0$}) \times GL(1, \mbox{$\bigwedge_{*-2}^0$}),
\end{equation}
define 
\begin{equation} 
\alpha \cdot_n (z, \theta^+, \theta^-) = \left( \frac{az + b}{cz + d} , \ \theta^+ \epsilon (cz + d)^{n-1}, \ \theta^- \epsilon^{-1} (cz + d)^{-n-1} \right) 
\end{equation}
for $(z, \theta^+, \theta^-) \in (\bigwedge_{*>1}^0 \smallsetminus \{-d_B/c_B \} \times (\bigwedge_{*>1}^0)_S) \times (\bigwedge_{*>1}^1)^2$.

The group $SL(2, \mbox{$\bigwedge_{*-2}^0$}) \times GL(1, \bigwedge_{*-2}^0)$ acts on $S^2 \hat{\mathbb{C}}(z^n)$ for $n \in \mathbb{Z}$ as global N=2 superconformal transformations as follows:  For each $n \in \mathbb{Z}$, and $\alpha \in SL(2, \mbox{$\bigwedge_{*-2}^0$}) \times GL(1, \bigwedge_{*-2}^0)$, define the map
\begin{multline}
T^{n,\alpha}_\sou : \bigl(\mbox{$\bigwedge_{*>1}^0$} \smallsetminus \bigl(\{- d_B/c_B \} \times (\mbox{$\bigwedge_{*>1}^0$})_S \bigr) \bigr) \oplus\mbox{$(\bigwedge_{*>1}^1)^2$}  \longrightarrow \\
\bigl( \mbox{$\bigwedge_{*>1}^0$} \smallsetminus \bigl(\{a_B/c_B \} \times
(\mbox{$\bigwedge_{*>1}^0$})_S\bigr) \bigr) \oplus \mbox{$(\bigwedge_{*>1}^1)^2$} 
\end{multline}
by 
\begin{equation}
T_\sou^{n,\alpha}(z, \theta^+, \theta^-) = \alpha \cdot_n (z, \theta^+, \theta^-)  .
\end{equation}
In addition, define 
\begin{multline}
T^{n,\alpha}_\nor: \bigl(\mbox{$\bigwedge_{*>1}^0$} \smallsetminus \bigl(\{-a_B/b_B \} \times
(\mbox{$\bigwedge_{*>1}^0$})_S \bigr)\bigr) \oplus \mbox{$(\bigwedge_{*>1}^1)^2$} \longrightarrow
\\
\bigl( \mbox{$\bigwedge_{*>1}^0$} \smallsetminus \bigr(\{d_B/b_B \} \times
(\mbox{$\bigwedge_{*>1}^0$})_S \bigr) \bigr) \oplus \mbox{$(\bigwedge_{*>1}^1)^2$}
\end{multline}
by 
\begin{equation}
T_\nor^{n,\alpha}(z, \theta^+, \theta^-) =  \left( \frac{c+ dz}{a+bz} , \ \theta^+ \epsilon (a + bz)^{n-1}, \ \theta^- \epsilon^{-1} (a + bz)^{-n-1} \right); 
\end{equation}
that is $T_\nor^{n,\alpha} (z, \theta^+, \theta^-) =I_{z^n}^{-1} \circ T^{n,\alpha}_\sou \circ I_{z^n} (z, \theta^+, \theta^-)$ for $(z,\theta^+, \theta^-) \in ((\bigwedge_{*>1}^0)^\times \smallsetminus (\{-a_B/b_B \} \times (\mbox{$\bigwedge_{*>1}^0$})_S) )\oplus  \mbox{$(\bigwedge_{*>1}^1)^2$}$.  

Let $\{ (U_\sou, \sou), (U_\nor, \nor) \}$ be the atlas for $S^2\hat{\mathbb{C}}(z^n)$ given by (\ref{southern-chart})--(\ref{transition}) with $g(z) = z^n$.   We define $T^{n,\alpha}: S^2\hat{\mathbb{C}}(z^n) \longrightarrow S^2 \hat{\mathbb{C}}(z^n)$ by 
\begin{equation}\label{T1}
T^{n,\alpha}(p) = \left\{
  \begin{array}{ll} 
      \sou^{-1} \circ T^{n,\alpha}_\sou \circ \sou (p) & \mbox{if $p \in U_\sou \smallsetminus X_1$}, \\  
      \nor^{-1} \circ T^{n,\alpha}_\nor \circ \nor (p) & \mbox{if $p \in
           U_\nor \smallsetminus X_2$},
\end{array} \right.
\end{equation}
where $X_1 = \sou^{-1}( (\{- d_B/c_B \} \times (\bigwedge_{*>1}^0)_S ) \oplus (\bigwedge_{*>1}^1)^2 )$ and $X_2 = \nor^{-1}(( \{-a_B/b_B \} \times (\bigwedge_{*>1}^0)_S ) \oplus  (\bigwedge_{*>1}^1)^2 )$.
This defines $T^{n,\alpha}$ for all $p \in S^2\hat{\mathbb{C}}(z^n)$ unless: \\
(i) $a_B = 0$ and  $p \in \nor^{-1}( (\{0 \} \times (\bigwedge_{*>1}^0)_S) \oplus  \mbox{$(\bigwedge_{*>1}^1)^2$}
)$; or \\
(ii) $d_B = 0$ and $p \in \sou^{-1}((\{0 \} \times (\bigwedge_{*>1}^0)_S) \oplus
\mbox{$(\bigwedge_{*>1}^1)^2$})$. \\
In case (i), we define
\begin{equation}\label{T2}
T^{n,\alpha} (p) = \sou^{-1} \left(  \frac{a + bz}{c + dz}  , \ i\theta^+ \epsilon (c+dz)^{n-1} ,\  i\theta^- \epsilon^{-1} (c+dz)^{-n-1} \right),
\end{equation} 
for $\nor (p) = (z,\theta^+, \theta^-) = (z_S, \theta^+, \theta^-)$.  

In case (ii), we define 
\begin{equation}\label{T3}
T^{n,\alpha}(p) = \nor^{-1} \left(   \frac{cz + d}{az + b} , \  -i  \theta^+ \epsilon (az + b)^{n-1} , \ 
-i \theta^- \epsilon^{-1} (az + b)^{-n-1}   \right) 
\end{equation} 
for $\sou (p) = (z,\theta^+, \theta^-) = (z_S, \theta^+, \theta^-)$. 

Note that with this definition, $T^{n,\alpha}$ is uniquely determined by $T^{n,\alpha}_\sou$, i.e., by its value on $\sou (U_\sou)$.  Or equivalently, $T^{n,\alpha}$ is uniquely determined by $T^{n,\alpha}_\nor$, i.e., by its value on $\nor (U_\nor)$.

The group of transformations determined by this action of $SL(2, \mbox{$\bigwedge_{*-2}^0$}) \times GL(1$, $\bigwedge_{*-2}^0)$ is a subgroup of the group of automorphisms of $S^2 \hat{\mathbb{C}} (z^n)$ for each $n \in \mathbb{Z}$.  In fact it is a proper subgroup, but we do not need this fact here; see \cite{B-n2moduli} for a proof of this fact in the case $n = 0$  and \cite{B-autogroups} for the cases $n \neq 0$.  By Lemma \ref{equivalent-lemma}, $SL(2, \mbox{$\bigwedge_{*-2}^0$}) \times GL(1, \bigwedge_{*-2}^0)$ also acts via N=2 superconformal automorphisms on $S^2 \hat{\mathbb{C}}(g)$, for $g \in \mathcal{G}$, in the obvious way. 

\begin{lem}\label{inequivalent-lemma}
Let $m,n \in \mathbb{Z}$. If $m \neq n$, then $S^2 \hat{\mathbb{C}} (z^m)$ and $S^2 \hat{\mathbb{C}} (z^n)$ are N=2 superconformally inequivalent. 
\end{lem}

\begin{proof} 
By first acting on $S^2 \hat{\mathbb{C}}(z^n)$ by a global N=2 superconformal transformation $T^{n,\alpha}$ for $\alpha \in SL(2, \mbox{$\bigwedge_{*-2}^0$})$, we can assume without loss of generality that $F$ sends the even component of the points $(0, \theta^+, \theta^-)$, $(1, \theta^+, \theta^-)$ and $(\infty, \theta^+, \theta^-)$ to the even points $0, 1$, and infinity, respectively.   That is, in terms of the local coordinate charts $\{(U_{\sou_m}, \sou_m), (U_{\nor_m}, \nor_m)\}$ and $\{(U_{\sou_n}, \sou_n), (U_{\nor_n}, \nor_n)\}$ for $S^2 \hat{\mathbb{C}}(z^m)$ and $S^2 \hat{\mathbb{C}}(z^n)$, respectively, we have $F(\sou_m^{-1}(0,$ $\theta^+, \theta^-)) =\sou_n^{-1}(0,\rho^+,\rho^-)$,  $F(\sou_m^{-1}(1,\theta^+, \theta^-)) =\sou_n^{-1}(1,\rho^+,\rho^-)$, and $F(\nor_m^{-1}(0,\theta^+,$ $\theta^-)) = \nor_n^{-1}(0,\rho^+,\rho^-)$.

Any N=2 superconformal equivalence that fixes these points is equivalent to a redefinition of the coordinates in the local coordinate charts $(U_{\nor_m}, \nor_m)$ and $(U_{\sou_m}, \sou_m)$ by automorphisms of the two copies of the N=2 superconformal plane $\nor_m(U_{\nor_m})$ and  $\sou_m(U_{\sou_m})$  that preserve the even coordinate.  

The automorphisms of the N=2 superconformal plane $S^2\mathbb{C}$ that preserve the even coordinate are of the form $T_{\varepsilon^+}(z, \theta^+, \theta^-) = (z, \theta^+ \varepsilon^+(z), \theta^- / \varepsilon^+(z))$ for $\varepsilon^+ \in \mathcal{E}^0 \leq \mathcal{E}$.   Thus changing coordinates in the chart $(U_{\nor_m}, \nor_m)$ of $S^2 \hat{\mathbb{C}} (z^m)$ by $T_{\varepsilon^+}$ for $\varepsilon^+ \in \mathcal{E}^0$ and changing coordinates in the chart $(U_{\sou_m}, \sou_m)$ of $S^2 \hat{\mathbb{C}} (z^m)$ by $T_{\varepsilon^-}$ for $\varepsilon^- \in \mathcal{E}^0$, the new change of coordinates from $(U_{\nor_m}, T_{\varepsilon^+} \circ \nor_m)$ to $(U_{\sou_m}, T_{\varepsilon^-} \circ \sou_m)$ is given by 
\begin{equation}\label{transformation-for-lemma2}
T_{\varepsilon^-} \circ I_{z^m} \circ T_{\varepsilon^+}^{-1} (z, \theta^+ ,\theta^-) = \left( \frac{1}{z}, \frac{i\theta^+ z^{m-1} \varepsilon^-(\frac{1}{z})}{ \varepsilon^+(z)}, \frac{i \theta^- \varepsilon^+(z)}{ z^{m+1} \varepsilon^-(\frac{1}{z})} \right) = I_{z^n} (z, \theta^+, \theta^-)
\end{equation}
for 
\begin{equation}
z^n = z^m \frac{\varepsilon^-(\frac{1}{z})}{\varepsilon^+(z)}.
\end{equation}

Since $\varepsilon^-(z) \in \mathcal{E}^0$, we have $\varepsilon^-(1/z) \in \mathcal{E}^\infty$.  And since $1/\varepsilon^+(z)$ is the inverse of $\varepsilon^+(z)$ in $\mathcal{E}^0$, we have that $\varepsilon^-(1/z)/\varepsilon^+(z) \in \mathcal{E}^0 \mathcal{E}^\infty=  \mathcal{E}$.  But this implies $z^m \mathcal{E} = z^n\mathcal{E}$ which implies that $m=n$. 
\end{proof}

Lemmas \ref{equivalent-lemma} and \ref{inequivalent-lemma} imply the following:

\begin{cor}\label{bijection-cor}  
There is a bijection between the set of N=2 superconformal equivalence classes of N=2 super-Riemann spheres in $\{S^2\hat{\mathbb{C}} (g) \; | \; g \in \mathcal{G} \}$ and the set $\mathcal{G}/\mathcal{E} \cong \mathbb{Z}$.  In other words, the quotient group $\mathcal{G}/\mathcal{E}$ classifies the N=2 superconformal structures of the form (\ref{southern-chart})-(\ref{transition}) over the Riemann sphere, and the moduli space of such N=2 super-Riemann spheres is given by $\{S^2\hat{\mathbb{C}} (z^n) \; | \; n \in \mathbb{Z}\}$.
\end{cor}

\subsection{The Uniformization Theorem for simply connected N=2 superconformal and N=1 superanalytic super-Riemann surfaces}\label{uniformization-subsection}

\begin{thm}\label{uniformization-thm}  
Any N=2 superconformal DeWitt super-Riemann surface with simply connected noncompact body is N=2 superconformally equivalent to the N=2 super plane $S^2 \mathbb{C}$ or the N=2 super upper half-plane $S^2\mathbb{H}$.  Any  N=2 superconformal DeWitt super-Riemann surface with genus-zero, simply connected compact body is N=2 superconformally equivalent to one of the unique N=2 superconformal structures over the Riemann sphere $\{ S^2 \hat{\mathbb{C}}(z^n) \; | \; n \in \mathbb{Z}\}$, where $S^2\hat{\mathbb{C}}(z^n)$ is given explicitly by the covering of local coordinate neighborhoods $\{ U_{\sou_n}, U_{\nor_n} \}$ and the local coordinate maps
\begin{eqnarray}
\sou_n  : U_{\sou_n}  & \longrightarrow & \mbox{$\bigwedge_{*>1}^0$} \oplus (\mbox{$\bigwedge_{*>1}^1$})^2  \\
\nor_n : U_{\nor_n} & \longrightarrow & \mbox{$\bigwedge_{*>1}^0$}  \oplus (\mbox{$\bigwedge_{*>1}^1$})^2,
\end{eqnarray}
which are homeomorphisms of $U_{\sou_n}$ and $U_{\nor_n}$ onto $\bigwedge_{*>1}^0 \oplus (\bigwedge_{*>1}^1)^2$, respectively, such that 
\begin{eqnarray}
\sou_n \circ \nor_n^{-1} : \mbox{$(\bigwedge_{*>1}^0)^\times$} \oplus \mbox{$(\bigwedge_{*>1}^1)^2$} &\longrightarrow& \mbox{$(\bigwedge_{*>1}^0)^\times$} \oplus \mbox{$(\bigwedge_{*>1}^1)^2$} \\ 
(z, \theta^+,\theta^-) & \mapsto &  \left(\frac{1}{z}, \frac{i \theta^+}{z} z^n, \frac{i \theta^-}{z} z^{-n} \right) . \nonumber
\end{eqnarray}
In particular, the moduli space of simply connected N=2 superconformal DeWitt super-Riemann surfaces under N=2 superconformal equivalence is isomorphic to the moduli space of simply connected N=2 superconformal $\mathcal{C}_{*>1}(2)$-supermanifolds under N=2 superconformal equivalence.  
\end{thm}

\begin{proof}
Since $M_B$ is simply connected, we have that both $\check{H}^1(M_B, \mathcal{L}^{-1} \otimes TM_B)$ and $\check{H}^1(M_B,TM_B)$ are trivial; see for instance \cite{Hartshorne}.  Thus by Theorem \ref{uniformization-prop}, $M$ is given by a coordinate atlas with coordinate transition functions of the form (\ref{reduced-transition}).   By the uniformization theorem for Riemann surfaces, we know that $M_B$ is conformally equivalent to $\mathbb{C}$, $\mathbb{H}$ or $\hat{\mathbb{C}}$.   This implies that for the body there exist coordinate redefinitions $f^\alpha_{(\emptyset)} : (\Omega_\alpha)_{(\emptyset)} ((U_\alpha)_B) \longrightarrow \mathbb{C}$ for $\alpha \in A$ such that the new atlas under these coordinate transformations is equivalent to the standard coordinate atlas for $\mathbb{C}$, $\mathbb{H}$, or $\hat{\mathbb{C}}$.  

It remains to show that there exist N=2 superconformal coordinate redefinitions that uniformize the body of $M$.  Furthermore, we must show that we can further reduce, under N=2 superconformal transformations, the coordinate atlas on $M$ to be of the form (\ref{southern-chart})--(\ref{transition}) for the compact case, and to be the usual coordinate atlases on $\mathbb{C}$ and $\mathbb{H}$ with trivial transition functions in the soul directions in the noncompact case.

Letting $f^\alpha_{(\emptyset)} : (\Omega_\alpha)_{(\emptyset)} ((U_\alpha)_B) \longrightarrow \mathbb{C}$ for $\alpha \in A$ be the coordinate redefinitions of $M_B$ taking $M_B$ to $\mathbb{C}$, $\mathbb{H}$, or $\hat{\mathbb{C}}$, let $H^\alpha : (\Omega_\alpha)_{(\emptyset)} ((U_\alpha)_B) \longrightarrow \bigwedge_*^0 \times (\bigwedge_*^1)^2$ be given by $H^\alpha(z, \theta^+, \theta^-) = (f^\alpha_{(\emptyset)} (z), \theta^+ (f^\alpha_{(\emptyset)})'(z), \theta^-)$, for $\alpha \in A$.   Under these N=2 superconformal coordinate redefinitions, we see that $M$ is N=2 superconformally equivalent to an N=2 super-Riemann surface whose body is $\mathbb{C}$, $\mathbb{H}$, or $\hat{\mathbb{C}}$, respectively, and with atlas $\{(U_\alpha, \Omega_\alpha) \; | \: \alpha \in A\}$ with transition functions given by 
\begin{equation}\label{reduce-to-uniformized-body}
H_{\alpha \beta}(z, \theta^+, \theta^-)  = \Omega_\alpha \circ \Omega_\beta^{-1} (z, \theta^+, \theta^-) = 
((f_{\alpha \beta})_{(\emptyset)} (z), \theta^+ g^+_{\alpha \beta}(z), \theta^- g^-_{\alpha \beta}(z))
\end{equation} 
for $\alpha, \beta \in A$, and where each $(f_{\alpha \beta})_{(\emptyset)} (z)$ is given by $z$ in the case of $M_B = \mathbb{C}$ and $M_B = \mathbb{H}$, or by $z$ or $z^{-1}$ in the case of $M_B = \hat{\mathbb{C}}$, and where the $g_{\alpha \beta}^\pm$ satisfy the cocycle condition
\begin{equation}\label{consistency-g3}
g^\pm_{\alpha \gamma} (z^\gamma) =  g^\pm_{\alpha \beta}(z^\beta) g^\pm_{\beta \gamma}(z^\gamma) ,
\end{equation}
on the coordinate atlas $\{(U_\alpha, \Omega_\alpha)\}_{\alpha \in A}$, for all $\alpha, \beta, \gamma \in A$ with $U_\alpha \cap U_\beta \cap U_\gamma \neq \emptyset$.  Note that by the superconformal condition (\ref{superconformal-condition-reduced}) the $g^-_{\alpha \beta}$ are completely determined by $g^+_{\alpha \beta}$.  

Define retractible submanifold(s) of $M_B$, each denoted by $M_{ret}$, to be $M_B$ itself if $M_B = \mathbb{C}$ or $\mathbb{H}$, to be $\hat{\mathbb{C}} \smallsetminus  \{\infty \}$ or $\mathbb{C}^\times \cup \{\infty \}$ if $M_B = \hat{\mathbb{C}}$.   Writing the $g_{\alpha \beta}^\pm$ for $\alpha, \beta \in A$ in component form
\begin{equation}
g_{\alpha \beta}^+ (z) =  \sum_{(j) \in J^0_{*-2}} (g^+_{\alpha \beta})_{(j)} \zeta_{j_1} \zeta_{j_2} \cdots \zeta_{j_{2n}},
\end{equation}
we note that the cocycle condition (\ref{consistency-g3}) at the zero level reduces to 
\begin{equation}
(g^+_{\alpha \gamma})_{(\emptyset)} (z^\gamma_{(\emptyset)}) =  (g^+_{\alpha \beta})_{(\emptyset)}(z^\beta_{(\emptyset)}) (g^+_{\beta \gamma})_{(\emptyset)} (z^\gamma_{(\emptyset)}) ,
\end{equation}
on the body of $M_B$.  On $M_{ret} \subseteq M_B$, this cohomology is trivial; that is, there exists a global section which gives a trivialization of the line bundle associated to these transition maps.  Let  $h^\alpha : (U_\alpha)_B \times \mathbb{C} \longrightarrow \mathbb{C}^2$ be the coordinate redefinitions which trivialize this line bundle.  Let $H^\alpha_{(\emptyset)} : \Omega_\alpha (U_\alpha) \longrightarrow \bigwedge_{*>1}^0 \oplus (\bigwedge_{*>1}^1)^2$ be given by $H^\alpha(z, \theta^+, \theta^-) = (z,\theta^+h^\alpha (z) , \theta^- (h^\alpha(z))^{-1} )$.  With these coordinate redefinitions,  the new coordinate transition functions over any retractible submanifold $M_{ret}$ of $M_B$ are of the form (\ref{reduced-transition}) with $g^\pm_{\alpha \beta}$ of the form
\begin{equation}
g_{\alpha \beta}^\pm (z) =  1+ \sum_{(j) \in J^0_{*-2} \smallsetminus \{(\emptyset)\}} (g^\pm_{\alpha \beta})_{(j)} \zeta_{j_1} \zeta_{j_2} \cdots \zeta_{j_{2n}}.
\end{equation}

We will show by induction on $n \in \mathbb{Z}_+$ that we can set the soul components of these $g^\pm$ terms equal to zero by N=2 superconformal coordinate redefinitions.  For $n=1$, let $(j) = (j_1, j_2) \in J^0_{*-2}$.  The consistency equation (\ref{consistency-g3})  implies
\begin{equation}\label{g-two}
(g^+_{\alpha \gamma})_{(j)}  (z^\gamma_{(\emptyset)}) = (g^+_{\beta \gamma})_{(j)} (z_{(\emptyset)}^\gamma) +  (g^+_{\alpha \beta})_{(j)} (z^\beta_{(\emptyset)}) .
\end{equation}
This implies that $( g^+_{\beta \gamma})_{(j)}$ is a cocycle in the first \v Cech cohomology of $M_
{ret}$.  Since $M_{ret}$ is retractible this cohomology is trivial, and there exists elements $h_{(j)}$ in the zeroth cohomology of $M_{ret}$, such that 
\begin{equation}\label{boundary1-g}
(g^+_{\beta \gamma})_{(j)}  =  (h_{(j)})_\gamma - (h_{(j)})_\beta.
\end{equation}

For $\alpha \in A$, define the N=2 superconformal coordinate transformation $H^\alpha_{(j)}$ by $f_\alpha(z^\alpha) = z^\alpha$, $g^+_\alpha(z^\alpha) = 1 + (h_{(j)})_\alpha(z^\alpha) \zeta_{j_1} \zeta_{j_2}$,  $g_\alpha^-(z^\alpha) =1 -  (h_{(j)})_\alpha(z^\alpha) \zeta_{j_1} \zeta_{j_2}$ and $\psi_\alpha^\pm (z^\alpha) = 0$. That is  $H^\alpha_{(j)} (z^\alpha, ( \theta^+)^\alpha, (\theta^-)^\alpha) = (z^\alpha,  \; \theta^+ (1 +  (h_{(j)})_\alpha(z^\alpha) \zeta_{j_1} \zeta_{j_2}), \; \theta^-(1 -  (h_{(j)})_\alpha(z^\alpha) \zeta_{j_1} \zeta_{j_2}))$. Redefining the coordinates for the charts $(U_\alpha, \Omega_\alpha)$ by the N=2 superconformal transformations $H^\alpha_{(j)}$, for $\alpha \in A$, such that the new coordinate charts are $(U_\alpha, \tilde{\Omega}_\alpha = H^\alpha_{(j)} \circ \Omega_\alpha)$, we have that the new coordinate transformations $\tilde{H}_{\alpha \beta} = \tilde{\Omega}_\alpha \circ \tilde{\Omega}_\beta^{-1} = H^\alpha_{(j)} \circ \Omega_\alpha \circ \Omega_\beta^{-1} \circ (H^\beta_{(j)})^{-1}$, for $\alpha, \beta \in A$, are now of the form $(z, \theta^+ \tilde{g}_{\alpha \beta}^+(z), \theta^- \tilde{g}_{\alpha \beta}^-(z))$ with 
\begin{eqnarray}
(\tilde{g}^\pm_{\alpha \beta})_{(\emptyset)} &=& 1\\
( \tilde{g}^+_{\alpha \beta})_{(j)} &=& ( g^+_{\alpha \beta})_{(j)} -   (h_{(j)})_\beta + (h_{(j)})_\alpha \; = \; 0\\
(\tilde{g}^-_{\alpha \beta})_{(j)} &=& 0.
\end{eqnarray}

For this new atlas of $M_{ret}$, we again have the compatibility condition on the triple overlaps, and thus the consistency condition (\ref{g-two}) again holds on the components of the new coordinate transformation functions, and we can perform the above procedure again for a different $(j) = (j_1, j_2) \in J^0_{*-2}$ without changing the fact that the previous $(g^\pm_{\alpha \beta})_{(j)}$ terms have been set equal to zero.  Doing this repeatedly, we can set the $(g^\pm_{\alpha \beta})_{(j)}$ terms equal to zero for all $(j) = (j_1, j_2) \in J^0_{*-2}$.  

Now we make the inductive assumption: assume coordinate transformations have been made such that the for the new $g^\pm$ terms, the components $(g^\pm_{\alpha \beta})_{(j)}$ have been set equal to zero for $(j) = (j_1, j_2, \dots, j_{2k}) \in J^0_{*-2}$, for $k = 1,\dots, n-1$.  For this new atlas, we again have the compatibility condition on the triple overlaps.  Thus for the $\zeta_{j_1} \cdots \zeta_{j_{2n}}$ level of equation (\ref{consistency-g3}) applied to $z_{(\emptyset)}$, we have that equation (\ref{g-two}) applies to these new $g^+_{(j)}$'s.  Thus there exists $h_{(j)}$ in the zeroth cohomology of $M_{ret}$ such that equation (\ref{boundary1-g}) holds for these terms.

For $\alpha \in A$, define the N=2 superconformal transformation $H^\alpha_{(j)}$ by $f_\alpha(z^\alpha) = z^\alpha$, $g_\alpha^+(z^\alpha) = 1 +  (h_{(j)})_\alpha(z^\alpha) \zeta_{j_1} \cdots \zeta_{j_{2n}}$, $g^-_\alpha(z^\alpha) = 1 - (h_{(j)})_\alpha(z^\alpha) \zeta_{j_1} \cdots \zeta_{j_{2n}}$, and $\psi_\alpha^\pm (z^\alpha) = 0$. Redefining the coordinates for the charts $(U_\alpha, \Omega_\alpha)$ by the N=2 superconformal transformations $H^\alpha_{(j)}$, for $\alpha \in A$, it is a straightforward calculation to see that the new coordinate transformations are of the form $(z, \theta^+ g^+(z), \theta^- g^-(z))$ with $g^\pm_{(\emptyset)} = 1$ and $g^\pm_{(j)} = 0$.  We perform this procedure of redefining the coordinate charts for each $(j) = (j_1, \dots, j_{2n}) \in J^0_{*-2}$, resulting in at atlas with coordinate transition functions that have all the $g^\pm_{(j)}$ components of length $2k$ for $k=1,\dots, n$ equal to zero.   By induction, we have that the $g^\pm_{\alpha \beta}$ terms of the transition functions on $M_{ret}$ are equal to the constant function 1. 

In the case where $M_B = M_{ret}$, i.e., $M_B = \mathbb{C}$ or $\mathbb{H}$, this proves that $M$ is N=2 superconformally equivalent to $S^2\mathbb{C}$ or $S^2\mathbb{H}$, respectively.  In the case where $M_{ret} = \hat{\mathbb{C}} \smallsetminus \{0\}$ or $\hat{\mathbb{C}} \smallsetminus \{\infty\}$, we obtain an N=2 super-Riemann surface whose body is $\hat{\mathbb{C}}$, and such that $M$ is covered by two coordinate charts $(U_\sou, \sou)$ and $(U_\nor, \nor)$ with coordinate transition function  $\sou \circ \nor^{-1} (z, \theta^+, \theta^-) = (z^{-1}, \theta^+ g^+(z), \theta^- g^-(z))$, for $g^+$ an even superanalytic $(1,0)$-superfunction defined and nonzero for all $z \in (\bigwedge_{*-2}^0)^\times$; in other words $g^+ \in \mathcal{G}$.  The condition that the transition functions be N=2 superconformal, and thus in particular, satisfy (\ref{superconformal-condition4}), implies that
\begin{equation}
\frac{-1}{z^2} = g^+(z) g^-(z).
\end{equation}
Letting $g(z) = -izg^+(z)$, we have 
\begin{equation}
\sou \circ \nor^{-1} (z, \theta^+, \theta^-) = (z^{-1},\  i\theta^+ g(z)/z,\  i\theta^- /(zg^-(z))) = I_g(z, \theta^+, \theta^-).
\end{equation}  
Thus $M$ is N=2 superconformally equivalent to $S^2 \hat{\mathbb{C}}(g)$ for some $g \in \mathcal{G}$.  By Corollary \ref{bijection-cor}, the result follows.
\end{proof}

It is easy to see that the classification of simply connected N=2 superconformal DeWitt super-Riemann surfaces up to N=2 superconformal equivalence in fact coincides with the classification of holomorphic line bundles over the underlying Riemann surface up to  holomorphic equivalence.   We can see this by observing that holomorphic line bundles are classified by $\check{H}^1(M_B, \mathbb{C}^\times)$.  For $M_B$ simply connected, this cohomology is trivial if $M_B$ is noncompact and $\mathbb{Z}$ if $M_B \cong \hat{\mathbb{C}}$. 

One can see this bijective correspondence between simply connected N=2 super-Riemann surfaces and holomorphic line bundles over $M_B \cong \hat{\mathbb{C}}$ explicitly, by noting that the N=2 super-Riemann sphere $S^2 \hat{\mathbb{C}}(z^n)$ for $n \in \mathbb{Z}$ has, as a substructure, the $GL(1, \mathbb{C})$-bundle over $\hat{\mathbb{C}}$ given by the transition function $iz_{(\emptyset)}^{n-1} : \mathbb{C}^\times \longrightarrow \mathbb{C}^\times$, corresponding to the transition function for the first fermionic component of $S^2\hat{\mathbb{C}}(z^n)$ restricted to the fiber in the first component of $\theta^+ = \theta_{(1)}^+ \zeta_1 + \theta_{(2)}^+ \zeta_2 + \cdots$.  (Or equivalently, one can restrict to any $(j)$-th component, for $(j) \in J^1_{*>1}$.)  Moreover, the $GL(1, \mathbb{C})$-bundle over $\hat{\mathbb{C}}$ with transition function $iz_{(\emptyset)}^{n-1} : \mathbb{C}^\times \longrightarrow \mathbb{C}^\times$, for $n \in \mathbb{Z}$, picks out a unique $S^2 \hat{\mathbb{C}}(z^n)$.    Under this bijection between equivalence classes of N=2 super-Riemann surfaces and equivalence classes of holomorphic line bundles over the body, the N=2 super-Riemann surface $S^2\hat{\mathbb{C}} (z^n)$ corresponds to the holomorphic line bundle over $\hat{\mathbb{C}}$ of degree $-n+1$.  Thus we have the following corollary:

\begin{cor}\label{U(1)-cor}  
N=2 superconformal DeWitt super-Riemann surfaces with simply connected body $M_B$ are classified up to N=2 superconformal equivalence by holomorphic line bundles over $M_B$ up to conformal equivalence.  
\end{cor}

Using Proposition \ref{DRS-prop} and Corollary \ref{uniformization-cor} we have the following corollary to Theorem \ref{uniformization-thm} and Corollary \ref{U(1)-cor}, which gives the uniformization for simply connected N=1 superanalytic DeWitt super-Riemann surfaces:
\begin{cor}\label{uniformization-cor-genus-zero}  
Any N=1 superanalytic DeWitt super-Riemann surface with simply connected noncompact body is N=1 superanalytically equivalent to the N=1 super plane $S^1 \mathbb{C} = \bigwedge_{*>0}$ or the N=1 super upper half-plane $S^1\mathbb{H}$.  Any  N=1 superanalytic DeWitt super-Riemann surface with genus-zero, simply connected compact body is N=1 superanalytically equivalent to one of the unique N=1 superanalytic structures over the Riemann sphere given explicitly by the covering of local coordinate neighborhoods $\{ U_{\sou_n}, U_{\nor_n} \}$ and the local coordinate maps 
$\sou_n  : U_{\sou_n}   \longrightarrow \bigwedge_{*>0}$ and $\nor_n : U_{\nor_n}  \longrightarrow  \bigwedge_{*>0}$, which are homeomorphisms of $U_{\sou_n}$ and $U_{\nor_n}$ onto $\bigwedge_{*>0}$, respectively, such that 
\begin{eqnarray}
\sou_n \circ \nor_n^{-1} : \mbox{$\bigwedge_{*>0}^\times$} &\longrightarrow& \mbox{$\bigwedge_{*>0}^\times$}\\ 
(z, \theta) & \mapsto &  \left(\frac{1}{z}, \frac{i \theta}{z} z^n \right) . \nonumber
\end{eqnarray}
In particular, the moduli space of simply connected N=1 superanalytic DeWitt super-Riemann surfaces under N=1 superanalytic equivalence is isomorphic to the moduli space of simply connected $\mathcal{C}_{*>0}(1)$-supermanifolds under N=1 superanalytic equivalence, and thus, is isomorphic to the moduli space of holomorphic line bundles over a simply connected Riemann surface under holomorphic equivalence.  
\end{cor}

\section{Uniformization for certain supermanifolds in the genus-one case}\label{torus-section}

We note that $\check{H}^1(M_B, \mathcal{L}^{-1} \otimes TM_B)$ and $\check{H}^1(M_B, TM_B)$ are both non-trivial if $M_B$ is a complex torus.  Thus by Theorem \ref{uniformization-prop} there are, in general, obstructions to uniformizing an N=2 superconformal DeWitt supertorus to one with transition functions of the form (\ref{reduced-transition}), i.e., to a $\mathcal{H}_{*>1}(2)$-supermanifold.  In addition, there are obstructions to further uniformizing an N=2 superconformal  $\mathcal{H}_{*>1}(2)$-supertorus to a $\mathcal{C}_{*>1}(2)$-supermanifold.   However, we can analyze the moduli space of genus-one N=2 superconformal DeWitt super-Riemann surfaces which have coordinate transition functions that correspond to the trivial cocycles in $\check{H}^1(M_B, \mathcal{L}^{-1} \otimes TM_B)$ and $\check{H}^1(M_B, TM_B)$.  In this section, we show that these supertori are classified up to N=2 superconformal equivalence, by holomorphic line bundles over the underlying complex torus, up to holomorphic equivalence.  Using the correspondence between N=2 superconformal and N=1 superanalytic super-Riemann surfaces given by Proposition \ref{DRS-prop}, this implies that N=1 superanalytic supertori corresponding to the trivial cocycles are also classified by holomorphic line bundles over the underlying complex torus, up to holomorphic equivalence.

\subsection{The moduli space of complex tori, automorphisms, and theta functions}

In this section we present some standard facts about complex tori and theta functions, following for example \cite{Br}, \cite{Miranda}, \cite{Debarre}.  Throughout this section and this section only, $z$ denotes a complex variable rather than an even supervariable. 

Let $\tau \in \mathbb{H}$, and let $\Gamma_\tau = \mathbb{Z} \oplus \tau \mathbb{Z}$.  The group $\Gamma_\tau$ acts on $\mathbb{C}$ by translation, and the quotient $\mathbb{C}/\Gamma_\tau$ defines a complex torus, also known as an elliptic curve.  The moduli space (up to conformal equivalence) of complex tori is given by $PSL(2,\mathbb{Z})\backslash \mathbb{H}$, where the action of $PSL(2,\mathbb{Z})$ on $\mathbb{H}$ is given by 
\begin{equation}
\left( \begin{array}{cc} a & b  \\
c & d \\
\end{array} \right) \cdot \tau = \frac{a \tau + b}{c\tau + d}
\end{equation}
for $a,b,c,d \in \mathbb{Z}$ satisfying $ad-bc = 1$.   From now on, when we refer to a lattice $\Gamma_\tau$ for some $\tau \in \mathbb{H}$, it is implied that we mean the equivalence class of $\tau$ in the moduli space $PSL(2,\mathbb{Z})\backslash \mathbb{H}$.

Let $\omega_n$ be a primitive $n$-th root of unity for $n \in \mathbb{Z}_+$.  The group of automorphisms for a complex torus $\mathbb{C}/\Gamma_\tau$ are given by translations by elements of $\Gamma_\tau$ along with the following groups written multiplicatively and acting on $\mathbb{C}/\Gamma_\tau$ via multiplication:\\
(i) $\langle \omega_4 \rangle =   \{ \omega_4^k \; | \; 1\leq k \leq 4 \}$ if $\tau = i$;\\
(ii)  $\langle \omega_6 \rangle = \{ \omega_6^k \; | \; 1\leq k \leq 6 \}$ if $\tau = e^{2\pi i/3}$;\\
(iii) $\langle \omega_2 \rangle = \{1, -1\}$ otherwise.

Now fix $\tau \in \mathbb{H}$.  Let $\pi_\tau : \mathbb{C} \longrightarrow \mathbb{C}/\Gamma_\tau$ be the canonical projection map.  Let $\{ ((U_\alpha)_B, (\Omega_\alpha)_B) \}_{\alpha \in A}$ be a coordinate atlas on $\mathbb{C}/\Gamma_\tau$ given by taking $\pi_\tau^{-1}((U_\alpha)_B)$ to be an open set in $\mathbb{C}$ such that $\gamma_1(\pi_\tau^{-1}((U_\alpha)_B)) \cap \gamma_2(\pi_\tau^{-1}((U_\alpha)_B)) = \emptyset$ for distinct $\gamma_1, \gamma_2 \in \Gamma_\tau$.   Then the coordinate transition functions for the chart on $\mathbb{C}/\Gamma_\tau$ are given by translation by elements of the lattice $\Gamma_\tau$.  

A theta function associated to $\Gamma_\tau$, denoted $\vartheta_\tau$,  is an entire function on $\mathbb{C}$ that is not identically zero such that for each $\gamma \in \Gamma_\tau$ there exist constants $a_\gamma, b_\gamma \in \mathbb{C}$ satisfying
\begin{equation}
\qquad \vartheta_\tau (z + \gamma) = e^{2 \pi i(a_\gamma z + b_\gamma)} \vartheta_\tau(z)  \qquad \mbox{for all $\gamma \in \Gamma_\tau$ and $z \in \mathbb{C}$.}
\end{equation}
The constants $\{a_\gamma, b_\gamma\}_{\gamma \in \Gamma_\tau}$ are called the {\it type} of the theta function $\vartheta_\tau$, and are equivalently defined as maps
\begin{equation}\label{ab-complex}
\begin{array}{rclrcl}
a : \Gamma_\tau  & \longrightarrow  &\mathbb{C} &\qquad \quad b : \Gamma_\tau  & \longrightarrow  &\mathbb{C}\\ 
\gamma & \mapsto  &a_\gamma &  \qquad  \gamma & \mapsto &  b_\gamma 
\end{array}
\end{equation}
satisfying
\begin{eqnarray}
a_{\gamma_1 + \gamma_2}   &=& a_{\gamma_1} + a_{\gamma_2} \label{ab1-C} \\
b_{\gamma_1 + \gamma_2}   &=& (b_{\gamma_1} + b_{\gamma_2} + a_{\gamma_1} \gamma_2) \, \mathrm{mod}\; \mathbb{Z}, \label{ab2-C}
\end{eqnarray}
for all $\gamma_1, \gamma_2 \in \Gamma_\tau$.   Thus we can also refer to the type of the theta function as the pair of maps $(a,b)$.

Let $\{g_\gamma : \mathbb{C} \longrightarrow \mathbb{C}^\times\}_{\gamma \in \Gamma_\tau}$ be a set of holomorphic functions which satisfy the condition
\begin{equation}\label{g-cocycle}
\; g_{\gamma_1 + \gamma_2}(z) =  g_{\gamma_2}(z) g_{\gamma_1}(z + \gamma_2) \qquad \mbox{for all $\gamma_1, \gamma_2 \in \Gamma_\tau$ and $z \in \mathbb{C}$}.
\end{equation}
Then
\begin{equation}\label{theta-ration}
g_\gamma (z) = \frac{\vartheta_\tau (z+ \gamma)}{\vartheta_\tau(z)} = e^{2\pi i (a_\gamma z + b_\gamma)}
\end{equation}
for some theta function $\vartheta_\tau$ of type $(a,b)$.  If $\vartheta_\tau$ and $\tilde{\vartheta}_\tau$ are of the same type, then they define the same function $g_\gamma$ via (\ref{theta-ration}).

Let $\Theta_\tau$ denote the space of theta functions associated to $\Gamma_\tau = \mathbb{Z} \oplus \tau \mathbb{Z}$ modulo equivalence up to type.  Note that $\Theta_\tau$ is a group under point-wise multiplication.  

A theta function $\vartheta_\tau$ associated to $\Gamma_\tau$ is called {\it trivial} if it is of the form 
\begin{equation}\label{trivial-theta}
\vartheta_\tau (z) = e^{az^2 +bz +c} \qquad \mbox{for some $a,b,c \in \mathbb{C}$.}
\end{equation}  
These are the theta functions that never vanish, and thus for such a $\vartheta_\tau$, we have that $\vartheta_\tau(z+\gamma)/\vartheta_\tau(z) = e^{2az\gamma + a\gamma^2 + b\gamma}$ is a nonvanishing entire function on $\mathbb{C}$.  A trivial theta function given by (\ref{trivial-theta}) is of type $\{a_\gamma, b_\gamma\}_{\gamma \in \Gamma}$ where $a_\gamma = \frac{-i}{\pi} a \gamma$ and $b_\gamma = \frac{-i}{2\pi} (b\gamma + a \gamma^2)$. 

Let $\mathcal{T}_\tau$ denote the set of trivial theta functions associated to $\Gamma_\tau$ modulo equivalence up to type.  Then $\mathcal{T}_\tau$ is in fact a subgroup of $\Theta_\tau$.  Two theta functions are said to be {\it equivalent} if their quotient is a trivial theta function.   Thus the quotient group $\Theta_\tau/\mathcal{T}_\tau$ is the moduli space of theta functions up to type modulo equivalence with respect to the trivial theta functions.   This space is often referred to as the set of {\it normalized} theta functions up to type, and such theta functions are often expressed uniquely via a pair consisting of a Hermitian form associated to the Riemann form of the theta function and a map $\alpha : \Gamma_\tau \longrightarrow U(1, \mathbb{C})$; this pair is called the {\it Appell-Humbert data} of the theta function \cite{Br}, \cite{Debarre}.

\subsection{A family of inequivalent N=2 superconformal structures over $\mathbb{C}/\Gamma_\tau$}\label{torus-superconformal-section}

An N=2 superconformal structure over $\mathbb{C}/\Gamma_\tau$ is an N=2 superconformal DeWitt super-Riemann surface $M$ with body $M_B = \mathbb{C}/\Gamma_\tau$, such that if $\{((U_\alpha)_B, (\Omega_\alpha)_B \}_{\alpha \in A}$ is the coordinate atlas for $\mathbb{C}/\Gamma_\tau$ as given in the previous section, then $M$ is covered by coordinate charts $\{(U_\alpha, \Omega_\alpha)\}_{\alpha \in A}$ where $\Omega_\alpha : U_\alpha \longrightarrow \bigwedge_{*>1}^0 \oplus (\bigwedge_{*>1}^1)^2$ maps $U_\alpha$ onto $(\Omega_\alpha)_B((U_\alpha)_B) \times  (\bigwedge_{*>1}^0)_S \oplus (\bigwedge_{*>1}^1)^2$  and the coordinate transition functions $\Omega_\alpha \circ \Omega_\beta^{-1}$ are N=2 superconformal for all $\alpha, \beta \in A$. 

\begin{rema}\label{torus-automorphisms-remark}
{\em The group of automorphisms of $\mathbb{C}/\Gamma_\tau$ extend to $M$ by acting via
$H(z, \theta^+, \theta^-) = (z+ \gamma, \theta^+, \theta^-)$ on translations, for $\gamma \in \Gamma_\tau$ and by $H(z, \theta^+, \theta^-) = (\omega z, \omega \theta^+, \theta^-)$ for multiplicative automorphisms with $\omega \in \mathbb{C}^\times$.  }
\end{rema}

\begin{defn}\label{define-super-theta}{\em
A {\it super-theta function} on $\bigwedge^0_*$ associated to $\Gamma_\tau$, denoted, $S\vartheta_\tau$ is an even superanalytic function $S\vartheta_\tau : \bigwedge^0_* \longrightarrow \bigwedge^0_*$ satisfying 
\begin{equation}
\quad \ S\vartheta_\tau (z + \gamma) = e^{2 \pi i(a_\gamma z + b_\gamma)} S\vartheta_\tau(z) \qquad \mbox{for all $\gamma \in \Gamma_\tau$ and $z \in \bigwedge_*^0$}.
\end{equation}
where 
\begin{equation}\label{ab}
\begin{array}{rclrcl}
a : \Gamma_\tau  & \longrightarrow  & \mbox{$\bigwedge_*^0$} &\qquad \quad b : \Gamma_\tau  & \longrightarrow  &\mbox{$\bigwedge_*^0$}\\ 
\gamma & \mapsto  &a_\gamma &  \qquad  \gamma & \mapsto &  b_\gamma 
\end{array}
\end{equation}
are even superfunctions on $\Gamma_\tau$ satisfying
\begin{eqnarray}
a_{\gamma_1 + \gamma_2}   &=& a_{\gamma_1} + a_{\gamma_2} \label{ab1} \\
b_{\gamma_1 + \gamma_2}   &=& (b_{\gamma_1} + b_{\gamma_2} + a_{\gamma_1} \gamma_2) \, \mathrm{mod}\; \mathbb{Z}, \label{ab2}
\end{eqnarray}
for all $\gamma_1, \gamma_2 \in \Gamma_\tau$.  The pair of maps $(a,b)$ on $\Gamma_\tau$ with values in $\bigwedge_*^0$ is called the {\it type} of the super-theta function. 

A {\it trivial} super-theta function on $\bigwedge^0_*$ associated to $\Gamma_\tau$ is one whose image lies in $(\bigwedge_*^0)^\times$. 
}
\end{defn}

In \cite{FR}, Freund and Rabin introduce a notion of super-theta function.  However, their notion is different than ours given here, in that they are interested in extending the notion of a classical theta function $\vartheta_\tau(z)$ to $\vartheta_{\tau + \theta \delta}(z)$ where $z$ is a complex variable $\theta$ is an odd super variable and $\delta$ is an odd parameter.  Rather, we are interested in extending the classical notion of theta function to one whose range is in $\bigwedge_*^0$, by extending the maps $(a,b)$ which define the type of the theta function to be even functions on the lattice rather than just complex-valued functions on the lattice.  (Note that in addition, we extend the domain to $\bigwedge_{*>1}^0$, but this is trivially done via (\ref{more-than-one-variable})).

\begin{rema}\label{super-theta-remark} 
{\em  Let $S\vartheta_\tau$ be a super-theta function on $\bigwedge_*^0$ of type $(a,b)$ as in Definition \ref{define-super-theta}.  Writing $a$ and $b$ as
\begin{equation}
a = \sum_{(j) \in J_*^0} a_{(j)} \zeta_{j_1}\cdots \zeta_{j_n} \quad \mbox{and} \quad  b = \sum_{(j) \in J_*^0} b_{(j)} \zeta_{j_1}\cdots \zeta_{j_n}
\end{equation}
we have from (\ref{ab1}) and (\ref{ab2}) that 
\begin{eqnarray}
a_{\gamma_1 + \gamma_2, (j)}   &=& a_{\gamma_1, (j)} + a_{\gamma_2, (j)} \label{ab1-component} \\
b_{\gamma_1 + \gamma_2, (j)}   &=& (b_{\gamma_1, (j)} + b_{\gamma_2, (j)} + a_{\gamma_1, (j)} \gamma_2) \, \mathrm{mod}\; \mathbb{Z}, \label{ab2-component}
\end{eqnarray}
for all $\gamma_1, \gamma_2 \in \Gamma_\tau$ and $(j) \in J_*^0$.  Thus a super-theta function on $\bigwedge_*^0$ of type $(a,b)$, is equivalent to a $2^{*-1}$-tuple $(\theta_{\tau, (\emptyset)}, \theta_{\tau, (12)}, \theta_{\tau, (13)}, \dots, \theta_{\tau, (1234)}, \dots )$ where $\theta_{\tau, (j)}$ for $(j) \in J^0_*$ is an ordinary theta function of type $(a_{(j)}, b_{(j)})$.
 }
\end{rema}

\begin{rema}\label{trivial-remark}
{\em
Let $\pi_B : \bigwedge_*^0 \longrightarrow \mathbb{C}$ denote the canonical projection onto the body of $\bigwedge_*^0 = (\bigwedge_*^0)_B \oplus (\bigwedge_*^0)_S$.   The trivial super-theta functions on $\bigwedge^0_*$ are exactly those super-theta functions such that the theta function $(S\vartheta_\tau)_B (z_B) = \pi_B \circ S\vartheta_\tau (z_B)$ is a trivial theta function. 
}
\end{rema}

Let $M$ be an N=2 superconformal DeWitt super-Riemann surface with $M_B = \mathbb{C}/\Gamma_\tau$.  For the remainder of this section, we restrict to the case when the coordinate transition functions for $M$ have their $\psi^\pm$ components and the soul part of their $f$ components equal to zero; that is, they are of the form 
\begin{equation}\label{torus-coordinate-transformations}
H_\gamma (z,\theta^+, \theta^-) = (z+ \gamma, \, \theta^+ g_\gamma (z), \, \theta^- (g_\gamma(z))^{-1})
\end{equation}
for $\gamma \in \Gamma_\tau$, and thus $M$ is a $\mathcal{H}_{*>1}(2)$-supermanifold. 
Then the compatibility condition on triple overlaps imposes the following condition on the transition functions
\begin{eqnarray}
\lefteqn{H_{\gamma_1 + \gamma_2} (z, \theta^+, \theta^-) }\\
&=& H_{\gamma_1} \circ H_{\gamma_2}(z, \theta^+, \theta^-) \nonumber \\
%&=& H_{\gamma_1} (z+ \gamma_2,\,  \theta^+ g_{\gamma_2}(z), \, \theta^- (g_{\gamma_2}(z))^{-1}) 
%\nonumber \\
&=& (z + \gamma_1 + \gamma_2, \, \theta^+ g_{\gamma_2}(z) g_{\gamma_1}(z + \gamma_2), \, \theta^- (g_{\gamma_2}(z) g_{\gamma_1}(z + \gamma_2))^{-1}), \nonumber
\end{eqnarray}
which implies that the even superanalytic functions $g_\gamma(z)$ for $\gamma \in \Gamma_\tau$ must satisfy
\begin{equation}\label{torus-cocycle}
 g_{\gamma_1 + \gamma_2}(z) =  g_{\gamma_2}(z) g_{\gamma_1}(z + \gamma_2) \qquad \mbox{for all $\gamma_1, \gamma_2 \in \Gamma_\tau$ and $z \in \bigwedge_{*>1}^0$}.
\end{equation}
Conversely, any family of even superanalytic functions $g_\gamma : \bigwedge_{*>1}^0 \longrightarrow (\bigwedge_{*>1}^0)^\times$ for $\gamma \in \Gamma_\tau$ satisfying (\ref{torus-cocycle}) defines an N=2 superconformal structure on $\mathbb{C}/\Gamma_\tau$ in this way.

\begin{lem}\label{theta-lemma}  
The solutions to equation (\ref{torus-cocycle}) are given by
\begin{equation}\label{g-from-theta}
g_\gamma (z) = \frac{S\vartheta_\tau (z+ \gamma)}{S\vartheta_\tau(z)} = e^{2 \pi i (a_\gamma z + b_\gamma)} 
\end{equation}
where $S\vartheta_\tau$ is a super-theta function on $\bigwedge_{*-2}^0$ associated to $\Gamma_\tau$ of type $(a,b)$. 
\end{lem}

\begin{proof}
Writing $g_\gamma$ in component form, we have
\begin{equation}
g_\gamma (z) = \sum_{(j) \in J^0_{*-2} }(g_\gamma)_{(j)} (z) \zeta_{j_1} \zeta_{j_2} \cdots \zeta_{j_{2n}}. 
\end{equation}
Restricting to $z_{(\emptyset)} \in \mathbb{C}$, this imposes the condition
\begin{equation}\label{torus-cocycle-level0}
(g_{\gamma_1 + \gamma_2})_{(\emptyset)}(z_{(\emptyset)}) =  (g_{\gamma_2})_{(\emptyset)}(z_{(\emptyset)}) (g_{\gamma_1})_{(\emptyset)}(z_{(\emptyset)} + \gamma_2) \quad \mbox{for all $\gamma_1, \gamma_2 \in \Gamma_\tau$}
\end{equation}
on the complex analytic functions $(g_\gamma)_{(\emptyset)} : \mathbb{C} \longrightarrow \mathbb{C}^\times$.
This implies that there is a theta function $\vartheta_\tau : \mathbb{C} \longrightarrow \mathbb{C}$ associated to $\Gamma_\tau$ such that 
\begin{equation}\label{level0-solutions}
(g_\gamma)_{(\emptyset)} (z_{(\emptyset)}) = \frac{\vartheta_\tau ( z_{(\emptyset)} + \gamma)}{\vartheta_\tau(z_{(\emptyset)})} = e^{2\pi i (a_{\gamma, (\emptyset)} z_{(\emptyset)} + b_{\gamma, (\emptyset)})}
\end{equation}
where $(a_{\gamma, (\emptyset)}, b_{\gamma, (\emptyset)})_{\gamma \in \Gamma_\tau}  = (a_{(\emptyset)}, b_{(\emptyset)})$ is the type of $\vartheta_\tau$.

Writing $g_\gamma(z) = e^{2 \pi i h_\gamma (z)}$, we have that the condition (\ref{torus-cocycle})  is equivalent to the condition 
\begin{equation}
h_{\gamma_1 + \gamma_2}(z) =  h_{\gamma_2}(z) + h_{\gamma_1}(z + \gamma_2) \qquad \mbox{for all $\gamma_1, \gamma_2 \in \Gamma_\tau$},
\end{equation}
and expanding $h_\gamma$ in component form, we have that for each $(j) \in J^0_{*-2}$
\begin{equation}\label{cocycle-condition-level-j}
h_{\gamma_1 + \gamma_2, (j)} (z_{(\emptyset)})  = h_{\gamma_2, (j)}(z_{(\emptyset)}) + h_{\gamma_1, (j)}(z_{(\emptyset)} + \gamma_2) \qquad \mbox{for all $\gamma_1, \gamma_2 \in \Gamma_\tau$}.
\end{equation}
Then the fact that any solution to (\ref{torus-cocycle-level0}) is of the form (\ref{level0-solutions}), implies that any solution to  (\ref{cocycle-condition-level-j}) is of the form 
\begin{equation}
h_{\gamma, (j)} (z_{(\emptyset)}) = a_{\gamma, (j)} z_{(\emptyset)} + b_{\gamma, (\emptyset)}
\end{equation}
for some 
\begin{equation}
\begin{array}{rclrcl}
a_{(j)} : \Gamma_\tau  & \longrightarrow  & \mathbb{C} &\qquad \quad b_{(j)} : \Gamma  & \longrightarrow  & \mathbb{C} \\
\gamma & \mapsto  &a_{\gamma, (j)}  &  \qquad  \gamma & \mapsto &  b_{\gamma, (j)}
\end{array}
\end{equation}
maps satisfying
\begin{eqnarray}
a_{\gamma_1 + \gamma_2, (j)}   &=& a_{\gamma_1, (j)} + a_{\gamma_2, (j)} \\
b_{\gamma_1 + \gamma_2, (j)}   &=& (b_{\gamma_1, (j)} + b_{\gamma_2, (j)} + a_{\gamma_1, (j)} \gamma_2) \, \mathrm{mod}\; \mathbb{Z}.
\end{eqnarray}
for all $\gamma_1, \gamma_2 \in \Gamma_\tau$.   Letting $a = \sum_{(j) \in J_{*-2}^0} a_{(j)} \zeta_{j_1}\cdots \zeta_{j_n}$ and $b = \sum_{(j) \in J_{*-2}^0} b_{(j)} \zeta_{j_1}\cdots$ $\cdot \zeta_{j_n}$, the result follows.  
\end{proof}

Let $S_*\Theta_\tau$ denote the set of super-theta functions on $\bigwedge^0_{*-2}$ associated to $\Gamma_\tau$ modulo equivalence up to type.  This set is a group under point-wise multiplication.  Let $S_*\mathcal{T}_\tau$ denote the subgroup of $S_*\Theta_\tau$ consisting of the trivial super-theta functions.  

If $M$ is an N=2 superconformal DeWitt super-Riemann surface over $\bigwedge_{*>1}$, with body $\mathbb{C}/\Gamma_\tau$ and with transition functions of the form (\ref{torus-coordinate-transformations}), then this structure defines a super-theta function over $\bigwedge_{*-2}^0$ associated to $\Gamma_\tau$ which is unique up to type.  Conversely, given a super-theta function $S\vartheta_\tau$ over $\bigwedge_{*-2}^0$ associated to $\Gamma_\tau$, then (\ref{torus-coordinate-transformations}) and (\ref{g-from-theta}) define the transition functions for an N=2 superconformal DeWitt super-Riemann surface over $\bigwedge_{*>1}$, with body $\mathbb{C}/\Gamma_\tau$. 

Denote by $S^2\mathbb{T}_\tau (S\vartheta_\tau)$ the N=2 superconformal DeWitt super-Riemann surface over $\bigwedge_{*>1}$, with body $\mathbb{C}/\Gamma_\tau$ and with transition functions determined by the super-theta function $S\vartheta_\tau \in S_*\Theta_\tau$. That is, the transition functions for $S^2\mathbb{T}_\tau (S\vartheta_\tau)$ are given by 
\begin{equation}
H_\gamma (z, \theta^+, \theta^-) = \left(z + \gamma, \; \theta^+ \frac{S\vartheta_\tau(z + \gamma)}{S\vartheta_\tau (z)}, \; \theta^- \frac{S\vartheta_\tau(z)}{S\vartheta_\tau(z + \gamma)} \right)
\end{equation}
for $\gamma \in \Gamma_\tau$. 

\begin{lem}\label{equivalent-torus-lemma}  
Let $S\vartheta^{(1)}_\tau, S\vartheta^{(2)}_\tau \in S_*\Theta_\tau$ be two super-theta functions on $\bigwedge_{*-2}^0$ associated to $\Gamma_\tau$.  If $S\vartheta^{(1)}_\tau (z)= S\vartheta^{(2)}_\tau(z) S\vartheta^T_\tau (z)$ for some trivial super-theta function $S\vartheta^T_\tau \in S_*\mathcal{T}_\tau$, then the N=2 superconformal DeWitt super-Riemann surfaces over $\mathbb{C}/\Gamma_\tau$ uniquely determined by $S\vartheta^{(1)}_\tau$ and $S\vartheta^{(2)}_\tau$, respectively, and denoted by $S^2\mathbb{T}_\tau(S\vartheta^{(1)}_\tau)$ and $S^2\mathbb{T}_\tau(S\vartheta^{(2)}_\tau)$, respectively, 
are  N=2 superconformally equivalent. 
\end{lem}

\begin{proof}
If $S\vartheta^T_\tau$ is a trivial super-theta function, then $S\vartheta^T_\tau \in \mathcal{E}^0$, and we have an automorphism of the N=2 superconformal plane $S^2\mathbb{C}$ that preserves the even coordinate given by $T(z, \theta^+, \theta^-) = (z, \theta^+ S\vartheta^T_\tau(z), \theta^- / S\vartheta^T_\tau(z))$.   

The transition functions for the N=2 superconformal DeWitt super-Riemann surface $S^2\mathbb{T}_\tau(S\vartheta^{(j)}_\tau)$ are given by $H^{(j)}_\gamma (z, \theta^+, \theta^-) = (z + \gamma, \theta^+ g_\gamma^{(j)}(z), \theta^- (g^{(j)}_\gamma(z))^{-1})$ where $g^{(j)}_\gamma (z) = S\vartheta^{(j)}_\tau (z + \gamma)/S\vartheta^{(j)}_\tau(z)$, for $j = 1,2,$ respectively.  Thus if $S\vartheta^{(1)}_\tau (z)= S\vartheta^T_\tau (z) S\vartheta^{(2)}_\tau(z)$, then $g^{(1)}_\gamma (z)=g^{(2)}_\gamma (z) S\vartheta^T_\tau (z + \gamma)/S\vartheta^T_\tau(z)$, and we have
\begin{eqnarray}
\qquad \ \ H_\gamma^{(1)}(z,\theta^+, \theta^-) &=&  (z + \gamma, \theta^+ g_\gamma^{(1)}(z), \theta^- (g^{(1)}_\gamma(z))^{-1}) \\
&=& \left(z + \gamma, \theta^+ \frac{g_\gamma^{(2)}(z) S\vartheta^T_\tau (z + \gamma)}{S\vartheta^T_\tau(z)}, \theta^- \frac{S\vartheta^T_\tau (z)}{g_\gamma^{(2)}(z) S\vartheta^T_\tau(z+ \gamma)} \right) \nonumber \\
&=& T^{-1} \circ H^{(2)}_\gamma \circ T(z, \theta^+, \theta^-).\nonumber
\end{eqnarray}
Therefore the automorphism of the N=2 superplane $T$ lifts to an N=2 superconformal bijection from  $S^2\mathbb{T}_\tau(S\vartheta^{(1)}_\tau)$ to $S^2\mathbb{T}_\tau(S\vartheta^{(2)}_\tau)$.
\end{proof}

\begin{lem}\label{theta-quotient-lemma} We have that $S_*\Theta_\tau/S_*\mathcal{T}_\tau \cong \Theta_\tau/\mathcal{T}_\tau$.  That is, the group of super-theta functions on $\bigwedge_{*-2}^0$ associated to $\Gamma_\tau$ up to type modulo the subgroup of trivial super-theta functions is isomorphic to the group of theta functions associated to $\Gamma_\tau$ up to type modulo the subgroup of trivial theta functions.
\end{lem}

\begin{proof} Let $S_*\vartheta_\tau \in S_*\Theta_\tau$ be of type $(a,b)$ for functions $(a,b)$ of the form (\ref{ab}) satisfying (\ref{ab1}) and (\ref{ab2}).   Let $a= a_B + a_S$ and $b= b_B + b_S$ be the body and soul decompositions of the functions $a$ and $b$, respectively.  That is $a_B = \pi_B \circ a$ and $a_S = \pi_S \circ a$, for $\pi_B : \bigwedge_*^0 \longrightarrow \mathbb{C}$ and $\pi_S : \bigwedge_*^0 \longrightarrow (\bigwedge_*^0)_S$ the canonical projections onto the body and soul, respectively, of $\bigwedge_*^0 = (\bigwedge_*^0)_B \oplus (\bigwedge_*^0)_S$, and similarly for $b$.

From Remark \ref{trivial-remark}, we have that the super-theta function of type $(a_S, b_S)$ is trivial.  Denoting this trivial super-theta function by $S\vartheta^T_\tau$, we have
\begin{equation}
\frac{S\vartheta_\tau(z + \gamma)}{S\vartheta^T_\tau(z + \gamma)} = e^{2 \pi i (a_{B,\gamma} z  + b_{B,\gamma})} \frac{ S\vartheta_\tau (z)}{S\vartheta^T_\tau (z)}
\end{equation}
implying that $S\vartheta_\tau/S\vartheta^T_\tau$ is a theta function.
\end{proof}

Lemmas \ref{equivalent-torus-lemma} and \ref{theta-quotient-lemma} imply that any N=2 superconformal supertorus with transition functions restricted to be of the form (\ref{torus-coordinate-transformations}) (or equivalently transition functions whose $\psi^\pm$ and soul $f$ components correspond to the zero cocycles in $\check H(M_B, \mathcal{L}^{-1} \otimes TM_B)$ and $\check H(M_B, TM_B)$, respectively) is N=2 superconformally equivalent to $S^2\mathbb{T}_\tau(\vartheta_\tau)$ for some $\vartheta_\tau \in \Theta_\tau/\mathcal{T}_\tau$.   The following lemma, will allow us to conclude that these are in fact representatives of distinct N=2 superconformal equivalence classes.

\begin{lem}\label{inequivalent-torus-lemma} 
Let $\vartheta^{(1)}_\tau, \vartheta^{(2)}_\tau \in \Theta_\tau$ be two theta functions associated to $\Gamma_\tau$ up to type.  If $\vartheta^{(1)}_\tau \mathcal{T}_\tau \neq \vartheta^{(2)}_\tau \mathcal{T}_\tau$ (that is, the ratio of $\vartheta^{(1)}_\tau$  to  $\vartheta^{(2)}_\tau$ is not a trivial theta function) then the N=2 super-Riemann surfaces over $\mathbb{C}/\Gamma_\tau$ uniquely determined by $\vartheta^{(1)}_\tau$ and $\vartheta^{(2)}_\tau$, respectively, and denoted by 
$S^2\mathbb{T}_\tau(\vartheta^{(1)}_\tau)$ and $S^2\mathbb{T}_\tau(\vartheta^{(2)}_\tau)$, respectively, are N=2 superconformally inequivalent.
\end{lem}

\begin{proof} Let $F: S^2\mathbb{T}_\tau(\vartheta^{(1)}_\tau) \longrightarrow S^2\mathbb{T}_\tau(\vartheta^{(2)}_\tau)$ be an N=2 superconformal equivalence.  By Remark \ref{torus-automorphisms-remark}, acting on $S^2\mathbb{T}_\tau(\vartheta^{(1)}_\tau)$ by global automorphisms, we can assume without loss of generality that $F$ restricted to $\mathbb{C}/\Gamma_\tau$ is the identity.   Thus $F$ is a transformation from the fiber bundle over $\mathbb{C}/\Gamma_\tau$ defined by $S^2\mathbb{T}_\tau(\vartheta^{(1)}_\tau)$ to the fiber bundle over $\mathbb{C}/\Gamma_\tau$ defined by $S^2\mathbb{T}_\tau(\vartheta^{(2)}_\tau)$.   Then to keep the even component of the transition functions equal to $z + \gamma$ and the $\psi^\pm$ components equal to zero, $F$ must be of the form $F(z, \theta^+, \theta^-) = (z, \theta^+ \varepsilon^+(z),\theta^- (\varepsilon^+(z))^{-1})$, for some $\varepsilon^+ \in \mathcal{E}^0 \cap S_*\mathcal{T}_\tau$.   

Moreover, the transition functions for the N=2 super-Riemann surface $S^2\mathbb{T}_\tau(\vartheta^{(j)}_\tau)$ are given by $H^{(j)}_\gamma (z, \theta^+, \theta^-) = (z + \gamma, \theta^+ g_\gamma^{(j)}(z), \theta^- (g^{(j)}_\gamma(z))^{-1})$ where $g^{(j)}_\gamma (z) = \vartheta^{(j)}_\tau (z + \gamma)/\vartheta^{(j)}_\tau(z)$, for $j = 1,2$, respectively, and $F$ must satisfy
\begin{equation}
H_\gamma^{(1)}(z,\theta^+, \theta^-) = F^{-1} \circ H^{(2)}_\gamma \circ F (z, \theta^+, \theta^-),
\end{equation}
i.e.,
\begin{multline}
(z + \gamma, \, \theta^+ g_\gamma^{(1)}(z), \, \theta^- (g^{(1)}_\gamma(z))^{-1}) \\
 = \left(z + \gamma, \, \theta^+ \frac{g_\gamma^{(2)}(z) \varepsilon^+ (z + \gamma)}{\varepsilon^+(z)}, \, \theta^- \frac{\varepsilon^+ (z)}{g_\gamma^{(2)}(z) \varepsilon^+(z+ \gamma)} \right).
\end{multline}
Thus $\varepsilon^+ \in \mathcal{T}_\tau$, and $\vartheta^{(1)}_\tau (z) \mathcal{T}_\tau = \vartheta^{(2)}_\tau(z) \mathcal{T}_\tau$.
\end{proof}

In the case of $M_B$ a complex torus, the Appell-Humbert Theorem states that the holomorphic line bundles over $\mathbb{C}/\Gamma_\tau$ are classified by $\Theta_\tau/\mathcal{T}_\tau$, that is by theta functions associated to $\Gamma_\tau$ up to type and equivalence by a trivial theta function cf. \cite{Br}.  The bijection between equivalence classes of N=2 superconformal DeWitt super-Riemann surfaces over a given torus with coordinate transition functions of the form (\ref{torus-coordinate-transformations}) and equivalence classes of holomorphic line bundles over this torus is given explicitly, for instance, by restricting to one of the $(j)$-th components, for $(j) \in J^1_{*>1}$, of the first fermionic component, i.e., the $\theta^+$ term.   Thus from Lemmas \ref{equivalent-torus-lemma}, \ref{theta-quotient-lemma}, and \ref{inequivalent-torus-lemma} we have the following theorem:

\begin{thm}\label{torus-theorem} 
The N=2 superconformal equivalence classes of  N=2 superconformal DeWitt super-Riemann surfaces with body $\mathbb{C}/\Gamma_\tau$, and which have coordinate transition functions of the form  (\ref{torus-coordinate-transformations}), are in one-to-one correspondence with theta functions associated to $\Gamma_\tau$ of a given type up to equivalence by the trivial theta functions.  That is, these N=2 superconformal super-Riemann surfaces are classified up to N=2 superconformal equivalence by $\Theta_\tau/\mathcal{T}_\tau \cong \{ S^2 \mathbb{T}_\tau (\vartheta_\tau) \; | \; \vartheta_\tau \in \Theta_\tau /\mathcal{T}_\tau\}$, or equivalently by holomorphic line bundles over $\mathbb{C}/\Gamma_\tau$ up to holomorphic equivalence.  

Similarly, N=1 superanalytic DeWitt super-Riemann surfaces with body $\mathbb{C}/\Gamma_\tau$ and which are also $\mathcal{H}_{*>0}(1)$-supermanifolds and such that the even transition function is just $z + \gamma$,  are in one-to-one correspondence with theta functions associated to $\Gamma_\tau$ of a given type up to equivalence by the trivial theta functions.  That is, these N=1 superanalytic super-Riemann surfaces are classified up to N=1 superanalytic equivalence by holomorphic line bundles over $\mathbb{C}/\Gamma_\tau$ up to holomorphic equivalence.  
\end{thm}

\section{The nonhomogeneous N=2 superconformal coordinates and an interpretation of uniformization in terms of loop groups}\label{nonhomo-section} 

In this section, we transfer some of our results to the nonhomogeneous N=2 supercoordinates since there are many results in N=2 superconformal field theory that employ this coordinate system.  This setting is also convenient for giving an interpretation of the Uniformization Theorems \ref{uniformization-thm} and \ref{torus-theorem}  in terms of $GL(1)$ loop groups over $\bigwedge_{*>1}^0$ as we do in Section \ref{line-bundle-subsection}.

\subsection{N=2 superconformal structures over $\hat{\mathbb{C}}$ and over $\mathbb{C}/\Gamma_\tau$ in nonhomogeneous coordinates}

In the nonhomogeneous coordinate system, the uniformized N=2 superconformal superspheres $S^2\hat{\mathbb{C}}(z^n)$, for $n \in \mathbb{Z}$, are given by the covering of local coordinate neighborhoods $\{ U_{\sou_n}, U_{\nor_n} \}$ and the local coordinate maps $\sou_n$ and $\nor_n$ which are homeomorphisms of $U_{\sou_n}$ and $U_{\nor_n}$ onto $\bigwedge_{*>1}^0 \oplus (\bigwedge_{*>1}^1)^2$, respectively, such that $\sou_n \circ \nor_n^{-1} : (\bigwedge_{*>1}^0)^\times \oplus (\bigwedge_{*>1}^1)^2 \longrightarrow (\bigwedge_{*>1}^0)^\times \oplus (\bigwedge_{*>1}^1)^2$ is given by 
\begin{multline}\label{transition-nonhomo}
\sou_n \circ \nor_n^{-1}(z, \theta_1, \theta_2) =  \left(\frac{1}{z}, \; \frac{i\theta_1}{2z} \left(z^n + z^{-n} \right) - \frac{\theta_2}{2z} \left(z^n - z^{-n} \right)  , \right. \\
\left. \frac{\theta_1}{2z} \left(z^n - z^{-n} \right) + \frac{i\theta_2}{2z} \left(z^n + z^{-n} \right)  \right) .
\end{multline}

For an N=2 superconformal super-Riemann surface with body $\mathbb{C}/\Gamma_\tau$ and with transition functions in the homogeneous coordinate system given by $H_\gamma (z, \theta^+, \theta^-) = (z + \gamma, \theta^+ g_\gamma (z),$ $\theta^- (g_\gamma (z))^{-1})$ for $g_\gamma(z) = e^{2 \pi i (a_\gamma z + b_\gamma)} = \vartheta_\tau(z + \gamma)/\vartheta_\tau(z)$ with $\vartheta_\tau \in \Theta_\tau/\mathcal{T}_\tau$, we have that in the nonhomogeneous coordinate system the transition functions are of the form 
\begin{eqnarray}\label{transition-nonhomo-torus}
\lefteqn{H_\gamma(z, \theta_1, \theta_2) }\\
&=& \left(z + \gamma, \; \frac{\theta_1}{2} \left(g_\gamma(z) + \frac{1}{g_\gamma(z)} \right) + \frac{i \theta_2}{2} \left(g_\gamma(z) - \frac{1}{g_\gamma(z)} \right), \right. \nonumber \\
& & \quad  \left.  \frac{-i\theta_1}{2} \left(g_\gamma(z) - \frac{1}{g_\gamma(z)} \right) + \frac{ \theta_2}{2} \left(g_\gamma(z) + \frac{1}{g_\gamma(z)} \right)\right) \nonumber \\
&=& \left(z + \gamma, \; \theta_1 \cosh (2\pi i(a_\gamma z + b_\gamma )) +i \theta_2 \sinh (2\pi i(a_\gamma z + b_\gamma )), \right. \nonumber\\
& & \quad  \left.  \-i\theta_1 \sinh (2\pi i(a_\gamma z + b_\gamma )) + \theta_2 \cosh (2\pi i(a_\gamma z + b_\gamma )) \right) \nonumber
\end{eqnarray}

From  (\ref{transition-nonhomo}) and (\ref{transition-nonhomo-torus}), it seems that the view from the homogeneous coordinate system is far less opaque than that from the nonhomogeneous system.  However, as we shall see in Section \ref{line-bundle-subsection}, the nonhomogeneous coordinate setting does give us an intuitive explanation for the classification of genus-zero N=2 superconformal DeWitt super-Riemann surfaces, and genus-one N=2 superconformal DeWitt super-Riemann surfaces corresponding to the trivial cocycles in $\check{H}^1(M_B, \mathcal{L}^{-1} \otimes TM_B)$ and $\check{H}^1(M_B,TM_B)$,  that we obtained in Sections \ref{uniformization-subsection} and \ref{torus-superconformal-section}.

\subsection{N=1 superconformal DeWitt super-Riemann surfaces, affine $\mathfrak{u}(1)$ and the $GL(1)$ loop group over $\bigwedge_{*>1}^0$}\label{line-bundle-subsection}

Much of the interpretation below was inspired in part by discussions the author had with Yi-Zhi Huang.

Recall that an N=1 superconformal DeWitt super-Riemann surface is a DeWitt $(1,1)$-dimensional supermanifold for which the transition functions, in addition to being superanalytic are N=1 superconformal.   As proved in \cite{B-memoir}, the Lie superalgebra of infinitesimal N=1 superconformal transformations is given by the superderivations
\begin{eqnarray}
L_n(z,\theta) &=& - \biggl( z^{n + 1} \frac{\partial}{\partial z} + (\frac{n + 1}{2})z^n   \theta \frac{\partial }{\partial \theta}  \biggr) \label{L-N=1} \\
G_{n -\frac{1}{2}} (z,\theta) &=& -  z^n \Bigl(   \frac{\partial }{\partial \theta}  -  \theta  \frac{\partial}{\partial z} \Bigr)   \label{G-N=1}
\end{eqnarray}
for $n\in \mathbb{Z}$.  Define the N=1 Neveu-Schwarz algebra to be the  Lie superalgebra with basis 
consisting of the central element $d$, even elements $L_n$, and odd elements $G_{n + 1/2}$ for $n \in \mathbb{Z}$, satisfying the supercommutation relations  
\begin{eqnarray}
\left[L_m ,L_n \right] \! \! &=& \! \! (m - n)L_{m + n} + \frac{1}{12} (m^3 - m) \delta_{m + n 
, 0} \; d , \label{Virasoro-relation1} \\
\left[ L_m, G_{n + \frac{1}{2}} \right] \! \! &=& \! \! \left(\frac{m}{2} - n - \frac{1}{2} \right) G_{m+n+\frac{1}{2}} ,\\
\left[ G_{m + \frac{1}{2}} , G_{n - \frac{1}{2}} \right] \! \! &=& \! \! 2L_{m + n}  + \frac{1}{3} (m^2 + m) \delta_{m + n , 0} \; d , \label{Neveu-Schwarz-relation-last}
\end{eqnarray}
for $m, n \in \mathbb{Z}$.   It is straightforward to check that the superderivations (\ref{L-N=1}), (\ref{G-N=1}) give a representation of the N=1 Neveu-Schwarz Lie superalgebra with central charge zero. 

As proved in \cite{B-n2moduli}, the Lie superalgebra of infinitesimal N=2 superconformal transformations in the nonhomogeneous coordinate system is given by the superderivations 
\begin{eqnarray}
L_n(z,\theta_1,\theta_2) &=& - \biggl( z^{n + 1} \frac{\partial}{\partial z} + (\frac{n + 1}{2})z^n \Bigl(  \theta_1 \frac{\partial }{\partial \theta_1}  + \theta_2 \frac{\partial }{\partial \theta_2 }  \Bigr) \biggr) \label{L-N=2} \\
J_n(z,\theta_1,\theta_2) &=& i z^n \biggl( \theta_1  \frac{\partial }{\partial \theta_2} -   \theta_2  \frac{\partial }{\partial \theta_1} \biggr)    \\
\qquad G^{(1)}_{n -\frac{1}{2}} (z,\theta_1,\theta_2) &=& - \biggl( z^n \Bigl(   \frac{\partial }{\partial \theta_1}  -  \theta_1  \frac{\partial}{\partial z} \Bigr)   - n z^{n-1} \theta_1 \theta_2   \frac{\partial }{\partial \theta_2}  \biggr)  \\
G^{(2)}_{n -\frac{1}{2}} (z,\theta_1,\theta_2) &=&   - \biggl( z^n \Bigl( \frac{\partial}{\partial \theta_2} - \theta_2 
\frac{\partial}{\partial z} \Bigr)  + n z^{n-1}  \theta_1 \theta_2 \frac{\partial}{\partial \theta_1} \biggr) \label{G-N=2}
\end{eqnarray}
for $n\in \mathbb{Z}$.  Define the N=2 Neveu-Schwarz algebra to be the  Lie superalgebra with basis 
consisting of the central element $d$, even elements $L_n$, $J_n$, and odd elements $G^{(j)}_{n + 1/2}$, for $n \in \mathbb{Z}$, $j = 1,2$, satisfying the super commutation relations  
\begin{eqnarray}
\left[L_m ,L_n \right] \! \! &=& \! \! (m - n)L_{m + n} + \frac{1}{12} (m^3 - m) \delta_{m + n 
, 0} \; d , \label{Virasoro-relation2} \\
\left[ L_m, G^{(j)}_{n + \frac{1}{2}} \right] \! \! &=& \! \! \left(\frac{m}{2} - n - \frac{1}{2} \right) G^{(j)}_{m+n+\frac{1}{2}} ,\\
\left[ G^{(j)}_{m + \frac{1}{2}} , G^{(j)}_{n - \frac{1}{2}} \right] \! \! &=& \! \! 2L_{m + n}  + \frac{1}{3} (m^2 + m) \delta_{m + n , 0} \; d , \\
\qquad \ \left[ G^{(1)}_{m + \frac{1}{2}} , G^{(2)}_{n - \frac{1}{2}} \right] \! \! &=& \! \!  - i  (m-n+1) J_{m+n} , \\
\left[ J_m, J_n \right] \! \! &=& \! \! \frac{1}{3} m \delta_{m+n,0} d, \qquad \qquad \qquad \qquad \ \  \left [L_m, J_n \right] \  = \   -n J_{m+n}, \\
\left[ J_m, G^{(1)}_{n + \frac{1}{2}} \right] \! \! &=& \! \! - i G^{(2)}_{m+n+\frac{1}{2}} , \qquad  \qquad \qquad  \ \ \ \left[ J_m, G^{(2)}_{n + \frac{1}{2}} \right]  \  = \   i G^{(1)}_{m+n+\frac{1}{2}} , \label{transformed-Neveu-Schwarz-relation-last}
\end{eqnarray}
for $m, n \in \mathbb{Z}$, $j=1,2$.   It is straightforward to check that the superderivations (\ref{L-N=2})--(\ref{G-N=2}) give a representation of the N=2 Neveu-Schwarz Lie superalgebra with central charge zero. 

Thus the N=2 Neveu-Schwarz Lie superalgebra has two copies of the N=1 Neveu-Schwarz Lie superalgebra (given by the $L_n$'s and $G^{(j)}_{n-1/2}$'s, for $j=1$ or $j=2$) along with a copy of a $\mathfrak{u}(1)$ affine Lie algebra running between them given by the $J_n$'s.  

\begin{rema}\label{homo-u(1)-remark} {\em  Defining $G^\pm_{n + \frac{1}{2}} = \frac{1}{\sqrt{2}} (G^{(1)}_{n+\frac{1}{2}} \mp i G^{(2)}_{n+\frac{1}{2}})$, and rewriting the supercommutation relations between $J_n$ and $G^\pm_{n + \frac{1}{2}}$ in terms of this new basis, we find that $[J_m, G^\pm_{n+\frac{1}{2}}] = \pm G^\pm_{m + n + \frac{1}{2}}$.   We call the basis $\{ L_n, J_n, G^\pm_{n + \frac{1}{2}}, d \; | \; n \in \mathbb{Z} \}$ the homogeneous basis for the N=2 Neveu-Schwarz algebra.  Performing the change of variables (\ref{transform-nonhomo-homo}) from nonhomogeneous coordinates to homogeneous coordinates, we have that the corresponding superderivations (\ref{L-N=2})--(\ref{G-N=2}), give a representation of the N=2 Neveu-Schwarz Lie superalgebra in the homogeneous basis.  And here we see another motivation for our terminology:  in the homogeneous basis, the $J_0$ term has homogeneous supercommutation relations with $G^\pm_{n + \frac{1}{2}}$.  Through the exponentiation of infinitesimal N=2 superconformal transformations, this results in the relative simplicity of the change of coordinate formulas in the homogeneous case versus the nonhomogeneous case for genus-zero (\ref{transition-nonhomo}) and restricted genus-one (\ref{transition-nonhomo-torus}) uniformized N=2 super-Riemann surfaces.  
}
\end{rema}

Another way of viewing the relationship between the N=1 and N=2 Neveu-Schwarz Lie superalgebras realized as superderivations, is by considering that the superderivations (\ref{L-N=1}) and (\ref{G-N=1}) along with 
\begin{equation}
J_n (z, \theta)   = z^n \theta \frac{\partial}{\partial \theta}  \qquad \mbox{and} \qquad G^*_{n-\frac{1}{2}} (z, \theta)  =  i z^n \left( \frac{\partial}{\partial \theta} + \theta \frac{\partial}{\partial z}\right) \label{recap-G*}
\end{equation}
for $n\in \mathbb{Z}$, also give a representation of the N=2 Neveu-Schwarz Lie superalgebra with central charge zero; see \cite{B-axiomatic-deformations}.  (Note that $G^*_{n-\frac{1}{2}} (z, \theta) = G_{n-\frac{1}{2}} (z, i \theta)$.)  As shown in \cite{B-axiomatic-deformations}, these superderivations generate the algebra of all N=1 superanalytic coordinate transformations, and these are generated by  (\ref{L-N=1}) and (\ref{G-N=1}) along with the zero term of the affine $\mathfrak{u}(1)$, namely $J_0(z,\theta) = \theta \frac{\partial}{\partial\theta}$.  

From the fact that the N=2 Neveu-Schwarz Lie superalgebra of infinitesimal N=2 superconformal transformations contains two subalgebras of infinitesimal N=1 superconformal transformations, which give rise to N=1 superconformal submanifold structures for certain N=2 super-Riemann surfaces, one can in some intuitive sense, think of the N=2 superconformal moduli space as arising from two copies of the N=1 superconformal moduli space with an exponentiated copy of affine $\mathfrak{u}(1)$ running between them.

More specifically, by the classification of N=1 superconformal DeWitt super-Riemann surfaces in \cite{CR}, up to N=1 superconformal equivalence, there is only one N=1 superconformal structure over $\mathbb{C}$, $\mathbb{H}$, and $\hat{\mathbb{C}}$, respectively.   The unique equivalence class of N=1 superconformal DeWitt super-Riemann surfaces with compact genus-zero body is given by two coordinate charts $\{(U_\sou, \sou), (U_\nor, \nor)\}$ and local coordinate maps $\sou: U_\sou \longrightarrow \bigwedge_{*>0}$ and $\nor:U_\nor \longrightarrow \bigwedge_{*>0}$ which are homeomorphisms of $U_\sou$ and $U_\nor$ onto $\bigwedge_{*>0}$, respectively, such that $\sou \circ \nor^{-1} : \bigwedge_{*>0}^\times \longrightarrow  \bigwedge_{*>0}^\times$ is given by $\sou \circ \nor^{-1} (z, \theta) = (1/z, i\theta/z)$.  We denote this genus-zero N=1 superconformal DeWitt super-Riemann surface by $S^1\hat{\mathbb{C}}(1)$.  The N=2 superconformal DeWitt super-Riemann surface $S^2\hat{\mathbb{C}}(1)$, has two embeddings of $S^1\hat{\mathbb{C}}(1)$, given by the identity mapping on the even subspace (the body and the even fiber component) and identifying the one fermionic component of $S^1\hat{\mathbb{C}}(1)$ with the fiber corresponding to either the first or the second fermionic component of $S^2 \hat{\mathbb{C}} (1)$.

In the genus-one case, the classification of N=1 superconformal DeWitt super-Riemann surfaces over a complex torus $\mathbb{C}/\Gamma_\tau$ was given in \cite{CR}  (see also \cite{FR}, \cite{Hodgkin}).  The N=1 superconformal transition functions $H_\gamma(z,\theta)$ are given by:
\begin{eqnarray}
\qquad H_{m+n\tau} (z,\theta) \! \! &=& \! \!  (z+ m + n b, \epsilon_1^m \epsilon_2^n \theta) \qquad  \qquad \mbox{(nontrivial spin structure})\\
H_{m +n\tau}(z, \theta) \! \! &=&  \! \! (z + m + n b +n\theta\delta, \theta  + n\delta) \qquad \mbox{(trivial spin structure})
\end{eqnarray}
for $m,n \in \mathbb{Z}$, where $b \in \bigwedge_{*-1}^0$ with $b_B = \tau$, $\delta \in \bigwedge_{*-1}^1$, and $(\epsilon_1, \epsilon_2) = (\pm 1, \mp 1)$ or $(-1,-1)$.   The first case corresponding to the nontrivial spin structure results in one distinct N=1 superconformal equivalence class for each $b \in \bigwedge_{*-1}^0$ with $b_B = \tau$.  The family of distinct N=1 superconformal equivalence classes in the case of trivial spin structure is parameterized by $b_S = b - \tau \in (\bigwedge_{*-1}^0)_S$ and $\delta \in \bigwedge_{*-1}^1/<\pm1>$.

The N=1 superconformal DeWitt super-Riemann surface with trivial spin structure and with $b = b_B = \tau$ and $\delta = 0$ is an N=1 superconformal submanifold of the N=2 superconformal DeWitt super-Riemann surface $S^2\mathbb{T}_\tau (1)$ with two unique embeddings given by mapping the N=1 fermionic component onto either the first or the second N=2 fermionic component for this trivial $(\bigwedge_{*>1}^0)_S \times (\bigwedge_{*>1}^1)^2$-bundle over $\mathbb{C}/\Gamma_\tau$.

The N=1 superconformal DeWitt super-Riemann surfaces with nontrivial spin structure and with $b = b_B = \tau$ are N=1 superconformally equivalent to an N=1 superconformal submanifold of $S^2 \mathbb{T}_\tau (\vartheta_\tau)$ with $\vartheta_\tau(z + m + n \tau)/\vartheta_\tau(z) = e^{\pi i n}$, i.e., where $\vartheta_\tau$ is  the theta function of type $(a_\gamma, b_\gamma)$ for $a_\gamma = 0$ and $b_\gamma = b_{m + n \tau} = n/2$.  Again the embedding can be done in two different ways: by embedding the N=1 fermionic component into either the first fermionic component or the second fermionic component.  This N=2 superconformal supertorus is the N=2 superconformal DeWitt super-Riemann surface over $\mathbb{C}/\Gamma_\tau$ with transition functions given by $H_{m+n\tau}(z, \theta^+, \theta^-) = (z+ m + n \tau, e^{\pi in} \theta^+, e^{-\pi in}\theta^-) = (z + m + n \tau, e^{\pi i n} \theta^+, e^{\pi i n} \theta^-)$ in the homogeneous coordinate system, which are coincidentally given by $H_{m+n\tau}(z, \theta_1, \theta_2) = (z+ m + n \tau, e^{\pi in} \theta_1, e^{\pi in}\theta_2)$ in the nonhomogeneous coordinate system.

Using the setting and results of \cite{B-n2moduli}, consider the group given by N=2 superconformal transformations of the form
\begin{multline}\label{group2}
\exp \Bigl(-  \sum_{n \in \Z} A_n J_n(z,\theta^+, \theta^-)  \Bigr)   \cdot  a_0^{-J_0(z,\theta^+, \theta^-)} \cdot (z, \theta^+, \theta^-) \\
= \Bigl( z, \; \theta^+ a_0 \exp \Bigl( \sum_{n \in \Z} A_n z^n \Bigr), \; \theta^- a_0^{-1} \exp \Bigl( -\sum_{n \in \Z} A_n z^n \Bigr) \Bigr)
\end{multline}
for $ a_0 \in (\bigwedge_{*-2}^0)^\times$, and $A_n \in \bigwedge_*^0$, for $n \in \Z$, 
where 
\begin{equation}\label{J-notation}
J_n(z,\theta^+,\theta^-) = - z^n\Bigl(\theta^+\frac{\partial}{\partial \theta^+} - \theta^- \frac{\partial}{\partial \theta^-}\Bigr) 
\end{equation}
and the series $\sum_{n \in \Z} A_n z^n$ has an infinite radius of convergence.  It follows from Theorem 6.10 in \cite{B-n2moduli}, there is a bijection between transformations of the form (\ref{group}) and of the form (\ref{group2}).  Similarly, exponentiating the $J_n(z, \theta^+, \theta^-)$ terms for $n \leq 0$ acting on $(1/z, i\theta^+/z, i\theta^-/z)$, one obtains the transformations in a neighborhood of infinity.  

This exponentiation of these infinitesimals $J_n(z, \theta^+, \theta^-)$, which represent the $\mathfrak{u}(1)$ affine Lie subalgebra of the N=2 Neveu-Schwarz algebra, over $\bigwedge_{*>1}^0$ gives us the full connected component of the identity in the $GL(1)$ loop group over $\bigwedge_{*>1}^0$ \cite{PS}, in the case of $n \in \mathbb{Z}$, which is the group $\mathcal{E}$.  And in the case $n \in \mathbb{N}$, we obtain the subgroup of the connected component of the $GL(1)$ loop group which, over $\bigwedge_{*>1}^0$, corresponds to the subgroup $\mathcal{E}^0$  of $\mathcal{E}$.  The case corresponding to $n\in \mathbb{N}$ and the subgroup $\mathcal{E}^0$ occurs when $M_B$ is noncompact, and the case corresponding to $n \in \mathbb{Z}$ and the subgroup $\mathcal{E}$ occurs when $M_B$ is compact.  In the compact genus-zero case, we have the group $\mathcal{G}$ corresponding to the full loop group, and $\mathcal{G}/\mathcal{E} \cong \mathbb{Z}$ counts the connected components, in addition to classifying the holomorphic $GL(1)$-bundles over $\hat{\mathbb{C}}$.   And similarly in the genus-one case, exponentiating the affine $\mathfrak{u}(1)$, over the two classes of genus-one N=1 superconformal DeWitt super-Riemann surfaces with trivial and non-trivial spin structure, respectively, we again arrive at the holomorphic $GL(1)$-bundles over the underlying body manifold giving rise to the moduli space of genus-one N=2 superconformal DeWitt super-Riemann surfaces with transition functions restricted to contain no odd functions of an even variable.


\begin{thebibliography}{CdGP}

\bibitem[B1]{B-announce} K. Barron, A supergeometric interpretation of vertex operator superalgebras, {\it Int. Math. Res. Notices}, 1996 No. 9, Duke University Press, 409--430.

\bibitem[B2]{B-thesis} K. Barron, {\it The supergeometric interpretation of vertex operator superalgebras}, Ph.D. Thesis, Rutgers University, October 1996.

\bibitem[B3]{B-vosas} K. Barron, N=1 Neveu-Schwarz vertex operator superalgebras over Grassmann algebras and with odd formal variables, in: {\it Proceedings of the International Conference on 
Representation Theory, 1998}, ed. by J. Wang and Z. Lin, China Higher Education Press \& Springer-Verlag, Beijing, 2000, 9--36.

\bibitem[B4]{B-memoir} K. Barron,  {\it The moduli space of $N = 1$ superspheres with tubes and the
sewing operation}, Memoirs Amer. Math. Soc., vol. 162, no. 772, (2003).

\bibitem[B5]{B-iso} K. Barron, The notion of N=1 supergeometric vertex operator superalgebra and the isomorphism theorem, {\it Commun. in Contemp. Math.}, vol. 5, no. 4, (2003), 481--567. 

\bibitem[B6]{B-change} K. Barron, Superconformal change of variables for N=1 Neveu-Schwarz vertex operator superalgebras, {\it J. of Algebra}, {\bf 277} (2004), 717--764. 

\bibitem[B7]{B-n2moduli} K. Barron, The moduli space of N=2 super-Riemann spheres with tubes, {\it Comm. in Contemp. Math.}, {\bf 9} (2007), 857--940.

\bibitem[B8]{B-axiomatic} K. Barron, Axiomatic aspects of N=2 vertex superalgebras with odd formal variables, {\it Commun. in Alg.}, {\bf 38} (2010), 1199--1268.

\bibitem[B9]{B-axiomatic-deformations} K. Barron, On axiomatic aspects of N=2 vertex superalgebras with odd formal variables, and deformations of N=1 vertex superalgebras, arXiv:0710.5755v1.

\bibitem[B10]{B-autogroups} K. Barron, Automorphism groups of N=2 superconformal super-Riemann spheres, {\it J. Pure Appl. Algebra}, {\bf 214} (2010), 1973--1987.

\bibitem[Bat]{Batchelor} M. Batchelor, Two approaches to supermanifolds, {\it Trans. of the Amer. Math. Soc.} {\bf 258} (1980), 257--270. 

\bibitem[BPZ]{BPZ} A. Belavin, A. Polyakov and A. Zamolodchikov, Infinite conformal symmetries in two-dimensional quantum field theory, {\it Nuclear Phys.}  {\bf B241} (1984), 333--380.

\bibitem[BR]{BR} M. Bergvelt, J. Rabin,  Supercurves, their Jacobians, and super KP equations, {\it Duke Math. J.} {\bf  98} (1999), no. 1, 1--57.

\bibitem[Br]{Br} V. Br\^inz\u anescu, {\it Holomorphic vector bundles over compact complex surfaces},
Lecture Notes in Math., {\bf 1624}, Springer-Verlag, Berlin, 1996.

\bibitem[CdGP]{CdGP} P. Candelas, X. de la Ossa, P. Green, and L. Parkes, A pair of Calabi-Yau manifolds as an exactly soluble superconformal field theory, {\it Nuclear Phys. B} {\bf 359} (1991) 21--74.

\bibitem[CK]{CK} D. Cox and S. Katz, {\it Mirror Symmetry and Algebraic Geometry}, Math. Surveys and Monographs Amer. Math. Soc, vol. 68, (1999).

\bibitem[CR]{CR} L. Crane and J. Rabin, Super Riemann surfaces: uniformization and Teichm\"uller theory, {\it Commun. Math. Phys.} {\bf 113} (1988), 601--623.

\bibitem[Deb] {Debarre} O. Debarre, {\it Complex Tori and Abelian Varieties}, SMF/AMS Texts and Monographs, Amer. Math. Soc., Providence, RI;  Soci\'et\'e Math\'ematique de France, Paris, {\bf 11} (2005).

\bibitem[DeW]{D} B. DeWitt, {\it Supermanifolds}, 2nd ed., Cambridge Monogr. Math. Phys., Cambridge Univ. Press, Cambridge, 1992.  

\bibitem[DPZ]{DPZ} P. Di Vecchia, J. Petersen and H. Zheng, N=2 extended superconformal theories in two dimensions, {\it Phys. Lett. } {\bf B162} (1985),  327--332.

\bibitem[DRS]{DRS} S. Dolgikh, A. Rosly and A. Schwarz, Supermoduli spaces, {\it Commun. in Math. Phys.} {\bf 135} (1990),  91--100.

\bibitem[DL]{DL} C. Dong and J. Lepowsky, Generalized Vertex Algebras and Relative Vertex Operators, {\it Progr. Math.} {\bf 112}, Birkh\"auser, Boston, 1993. 

\bibitem[FaR]{Falqui-Reina-1} G. Falqui and C. Reina, $N=2$ super Riemann surfaces and algebraic geometry, {\it J. Math. Phys.} {\bf 31}  (1990), no. 4, 948--952.

\bibitem[FFR]{FFR} A. Feingold, I. Frenkel and J. Ries, {\it Spinor Construction of Vertex Operator Algebras, Triality and $E_8^{(1)}$}, Contemp. Math. {\bf 121}, Amer. Math. Soc., Providence, 
1991.  

\bibitem[FLM]{FLM} I. B. Frenkel, J. Lepowsky and A. Meurman, {\it Vertex Operator Algebras and the Monster}, Academic Press, New York, 1988. 

\bibitem[FrR]{FR} P. Freund and J. Rabin, Supertori are Algebraic Curves, {\it Commun. in Math. Phys.} {\bf 114} (1988), 131--145.

\bibitem[F]{Fd} D. Friedan, Notes on string theory and two-dimensional conformal field theory, in: {\it Unified String Theories}, World Scientific, Singapore, 1986, 162--213.

\bibitem[FMS]{FMS}  D. Friedan, E. Martinec and S. Shenker,  Conformal invariance, supersymmetry and string theory, {\it  Nucl. Phys.} {\bf B271} (1986), no. 1, 93--165. 

\bibitem[Ge]{Ge} D. Gepner, Lectures on N=2 string theory, in: {\it Superstrings '89 (Trieste, 1989)},  World Sci. Publishing, River Edge, NJ, 1990, 238--302. 

\bibitem[Gi]{Gi} A. Givental, Elliptic Gromov-Witten invariants and the  generalized mirror conjecture, in: {\it Integrable systems and algebraic geometry (Kobe/Kyoto, 1997)}, World Sci. Publishing, River Edge, NJ, 1998, 107--155.

\bibitem[GP]{GP} B. Greene and M. Plesser, Duality in Calabi-Yau moduli space, {\it Nuclear Phys. B} {\bf 338} (1990),  15--37. 

\bibitem[Ha]{Hartshorne} R. Hartshorne, {\it Algebraic geometry}, Graduate Texts in Mathematics, No. 52. Springer-Verlag, New York-Heidelberg, 1977. 

\bibitem[Ho]{Hodgkin}  L. Hodgkin, A direct calculation of super-Teichm\"uller space, {\it Letters in Math. Phys.}, {\bf 14} (1987), 47--53. 

\bibitem[HKK]{HKKP} K. Hori, S. Katz, A. Klemm, R. Pandharipande, R. Thomas, C. Vafa, R. Vakil, and E. Zaslow, {\it Mirror symmetry}, Clay Math. Monographs, 1, Amer. Math. Soc., Providence, 2003.

\bibitem[H1]{H-book} Y.-Z. Huang, {\it Two-Dimensional Conformal Geometry and Vertex Operator Algebras}, Progress in Math. Vol. 148, Birkh\"auser, Boston, 1997.

\bibitem[H2]{H-tensor} Y.-Z. Huang, A theory of tensor products for module categories for a vertex operator algebra, IV, {\it J. Pure Appl. Alg.} {\bf 100} (1995), 173--216.

\bibitem[H3]{H-orbifold} Y.-Z. Huang, A nonmeromorphic extension of the moonshine module vertex operator algebra, in: {\it Moonshine, the Monster and related topics, Proc. Joint Summer Research Conference, Mount Holyoke, 1994}, ed. C. Dong and G. Mason, Contemp. Math., Vol. 193, Amer. Math. Soc., Providence, 1996, 123--148.

\bibitem[H4]{H1997} Y.-Z. Huang, Intertwining operator algebras, genus-zero modular functors and genus-zero conformal field theories, in {\it Operads: Proceedings of Renaissance Conferences}, eds. J.L. Loday, J. Stasheff and A. A. Voronov, Contemp. Math., Vol. 202, Amer. Math.Soc., Providence, 1997, 335--355. 

\bibitem[H5]{H1998} Y.-Z. Huang, Genus-zero modular functors and intertwining operator algebras, {\it Internat. J. Math.} {\bf  9} (1998) 845--863. 

\bibitem[H6]{H1999} Y.-Z. Huang, A functional-analytic theory of vertex (operator) algebras, I, {\it Commun. Math. Phys.} {\bf 204} (1999) 61--84. 

\bibitem[H7]{H2003} Y.-Z. Huang, A functional-analytic theory of vertex (operator) algebras, II, {\it Commun. Math. Phys.} {\bf 242} (2003) 425--444. 

\bibitem[H8]{H2005} Y.-Z. Huang, Differential equations and intertwining operators, {\it Commun. Contemp. Math.} {\bf 7} (2005) 375--400. 

\bibitem[H9]{H2005-2} Y.-Z. Huang, Differential equations, duality and modular invariance, {\it Commun. Contemp. Math.} {\bf 7} (2005),  649--706.

\bibitem[HL1]{HL1} Y.-Z. Huang and J. Lepowsky, Toward a theory of tensor products for representations of a vertex operator algebra, in: {\it Proc. 20th International Conference on Differential Geometric
Methods in Theoretical Physics, New York, 1991}, ed. S. Catto and A. Rocha, World Scientific, Singapore, 1992, 344--354.

\bibitem[HL2]{HL2} Y.-Z. Huang and J. Lepowsky, Vertex operator algebras and operads, in: {\it The Gelfand Mathematical Seminar, 1990--1992}, ed. L. Corwin, I. Gelfand and J. Lepowsky, Birkh\"auser, Boston, 1993, 145--161.

\bibitem[HL3]{HL3} Y.-Z. Huang and J. Lepowsky, Operadic formulation of the notion of vertex operator algebra, in: {\it Mathematical Aspects of Conformal and Topological Field Theories and Quantum
Groups}, ed. P. Sally, M. Flato, J. Lepowsky, N. Reshetikhin and G. Zuckerman, Contemp. Math. Vol. 175, Amer. Math. Soc., Providence, 1994, 131--148.

\bibitem[HL4]{HL4} Y.-Z. Huang and J. Lepowsky, A theory of tensor products for module categories for a vertex operator algebra, I, {\it Selecta Mathematica, New Series} {\bf 1} (1995), 699--756.

\bibitem[HL5]{HL5} Y.-Z. Huang and J. Lepowsky, A theory of tensor products for module categories for a vertex operator algebra, II, {\it Selecta Mathematica, New Series} {\bf 1} (1995), 757--786.

\bibitem[HL6]{HL6} Y.-Z. Huang and J. Lepowsky, Tensor products of modules for vertex operator algebras and vertex tensor categories, in: {\it Lie Theory and Geometry, in honor of Bertram
Kostant}, ed. R. Brylinski, J.-L. Brylinski, V. Guillemin and V. Kac, Progr. Math, Vol. 123, Birkh\"auser, Boston, 1994, 349--383.

\bibitem[HL7]{HL7} Y.-Z. Huang and J. Lepowsky, A theory of tensor products for module categories for a vertex operator algebra, III, {\it J. Pure Appl. Alg.} {\bf 100} (1995), 141--171.

\bibitem[KT]{KT} Y. Kanie and A. Tsuchiya, Vertex operators in conformal field theory on $\mathbb{P}^1$ and monodromy representations of braid groups, in: {\it Conformal Field Theory and Solvable Lattice
Models, Advanced Studies in Pure Math.} {\bf 16}, Kinokuniya Company Ltd., Tokyo, 1988, 297--372.

\bibitem[KS]{KS} A. Konechny and A. Schwarz, Theory of $(k\oplus l\vert q)$-dimensional supermanifolds, {\it Selecta Math. (N.S.)} {\bf  6}  (2000), 471--486.

\bibitem[LR]{LR} C. LeBrun and M. Rothstein, Moduli of super Riemann surfaces, {\it Commun. Math. Phys.} {\bf 117} (1988), 159--176.

\bibitem[L]{Leites}  D. Leites, Introduction to the theory of supermanifolds, {\it Uspekhi Mat. Nauk} {\bf 35} (1980), 3--57. 

\bibitem[LVW]{LVW} W. Lerche, C. Vafa, and N.  Warner, Chiral rings in N=2 superconformal theories, {\it Nuclear Phys. B} {\bf 324} (1989), 427--474.  

\bibitem[M1]{Manin1} Y. Manin, {\it Gauge field theory and complex geometry}, Springer-Verlag, 1988.  Translated by N. Koblitz and J. R. King from the 1984 Russian original. 

\bibitem[M2]{Manin2} Y. Manin, {\it Topics in noncommutative geometry}, Princeton University Press, Princeton, NJ, 1991. 

\bibitem[MR]{Manin-Radul} Y. Manin, A. Radul, A supersymmetric extension of the Kadomtsev-Petviashvili hierarchy, {\it Comm. Math. Phys.} {\bf  98} (1985), no. 1, 65--77.

\bibitem[Mil]{Milas} A. Milas,  Characters, supercharacters and Weber modular functions, {\it J. Reine Angew. Math} {\bf  608} (2007), 35--64

\bibitem[Mir]{Miranda} R. Miranda, {\it Algebraic curves and Riemann surfaces}, Grad. Studies in Math., Amer. Math. Soc., Providence, RI, {\bf 5} (1995).

\bibitem[PS]{PS} A. Pressley and G. Segal, {\it Loop Groups}, Oxford Math. Monogr., Oxford University Press, New York, 1986.

\bibitem[Ra]{Ra}  J. Rabin, Super elliptic curves, {\it J. Geom. Phys.} {\bf 15} (1995), no. 3, 252--280.
 
\bibitem[Rog]{Ro} A. Rogers, {\it Supermanifolds, Theory and Applications}, World Scientific Publishing Co. Pte. Ltd., Hackensack, NJ, 2007.

\bibitem[Rot]{Rothstein} M. Rothstein, Deformations of complex supermanifolds, {\it Proc. Amer. Math. Soc.} {\bf 95} (1985), p. 255--260.

\bibitem[S]{Schwarz1984} A. Shvarts, On the definition of superspace, {\it  Teoret. Mat. Fiz.} {\bf 60}(1984),  37--42. 

\bibitem[STT]{STT}  A. Semikhatov, A. Taormina, Y. Tipunin,  Higher-level Appell functions, modular transformations, and characters, {\it Comm. Math. Phys.} {\bf 255} (2005), no. 2, 469--512.

\bibitem[V]{Vaintrob} A. Vaintrob, Deformations of complex supermanifolds and coherent sheaves on them, {\it J. Sov. Math} {\bf 51} (1990), 2140--2188. 

\bibitem[W]{Wa} N. Warner, Lectures on N=2 superconformal theories and singularity theory, in: {\it Superstrings '89 (Trieste, 1989)}, World Sci. Publishing, River Edge, NJ, 1990, 197--237. 

\bibitem[YZ]{YZ} M. Yu and H. Zheng {\it N=2 superconformal invariance in two-dimensional quantum field theories}, {\it Nucl. Phys.} {\bf B288} (1987), 275--300.
 
\end{thebibliography}
\end{document}